\documentclass[11pt]{article}
\usepackage[english]{babel}
\usepackage{latexsym}
\usepackage{epsfig}
\usepackage{graphicx}
\usepackage{comment}
\usepackage{amsfonts,amsmath,mathtools}
\usepackage{color}

\def\dv{\mathbb}
\def\caP{{\cal P}(N)}
\def\m{m} 
\def\am{m^{\diamond}} 
\def\dual{\ast}
\def\refl{\star}
\def\0{\mbox{\sf 0}}
\def\calB{{\cal B}}
\def\calC{{\cal C}}
\def\calD{{\cal D}}
\def\calE{{\cal E}}
\def\calI{{\cal I}}
\def\calL{{\cal L}}

\def\calS{{\cal S}}
\def\calT{{\cal T}}
\def\calW{{\cal W}}
\def\calZ{{\cal Z}}
\def\cone{\mbox{cone}\,}
\def\ext{\mbox{ext}\,}
\def\conv{\mbox{conv}\,}
\def\scon{\mbox{\rm semi-cone}}

\definecolor{mygrey}{gray}{0.6}
\definecolor{myorange}{cmyk}{0,0.4,1,0}

\newtheorem{thm}{Theorem}
\newtheorem{example}{Example}
\newtheorem{remark}{Remark}

\newtheorem{defin}{Definition}
\newtheorem{lem}[thm]{Lemma}
\newtheorem{cor}[thm]{Corollary}
\newenvironment{proof}{\begin{trivlist}\item[] \mbox{\it Proof. }}
{\hfill$\Box$ \end{trivlist}}

\begin{document}

\title{Facets of the cone of exact games}
\author{Milan Studen\'{y}, V\'{a}clav Kratochv\'{\i}l}
\date{\today} 
\maketitle

\begin{abstract}
The class of {\em exact\/} transferable utility coalitional games, introduced in 1972 by Schmeidler, has been
studied both in the context of game theory and in the context of imprecise probabilities.
We characterize the cone of exact games by describing the minimal set of linear inequalities defining this cone; these facet-defining inequalities for the exact cone appear to correspond to certain
set systems (= systems of coalitions). We noticed that non-empty proper coalitions having non-zero coefficients
in these facet-defining inequalities form set systems with particular properties.

More specifically, we introduce the concept of a {\em semi-balanced\/} system of coalitions, which generalizes the classic concept of a balanced coalitional system in cooperative game theory. The semi-balanced coalitional systems provide valid inequalities for the exact cone and minimal semi-balanced systems (in the sense of inclusion of set systems) characterize this cone. We also introduce basic classification of minimal
semi-balanced systems, their pictorial representatives  and
a substantial concept of an {\em indecomposable\/} (minimal)
{\em semi-balanced system\/} of coalitions.
The main result of the paper is that indecomposable semi-balanced systems are in one-to-one correspondence with facet-defining inequalities for the exact cone.
The secondary relevant result is the rebuttal of a former conjecture
claiming that a coalitional game is exact iff
it is totally balanced and its anti-dual is also totally balanced.
We additionally characterize those inequalities which are facet-defining both
for the exact cone and the cone of totally balanced games.
\end{abstract}

\section{Introduction}\label{sec.introd}
The class of {\em exact\/} (transferable utility) games is
one of the topics of interest in cooperative game theory;
see \cite[\S\,V.1]{Rosenmuller00} or \cite[\S\,3.4]{Grabish16}.
Exact games were introduced by Schmeidler \cite{Schmeidler72} in 1972
within a wide framework of cooperative games with possibly infinite amount of players.
In this paper we, however, consider the usual game-theoretical framework
and assume that the set $N$ of players (for considered cooperative games) is a (fixed) non-empty finite set.
Mathematically equivalent concept of  a {\em coherent lower probability\/} (to the one of an exact game) has later appeared in the context of imprecise probabilities \cite{Walley91}, where $N$ has the interpretation of the sample space for considered (discrete) probability distributions.
Exact games have various applications described in detail in \cite[\S\,1]{CHK11}.
For example, it was shown in \cite{CHK09} that exact games coincide with risk allocation games with no aggregate uncertainty and in \cite{CBH05} that non-negative exact games coincide with multi-issue allocation games.

\subsection{Overview of\/ former related results}\label{ssec.overview}
It was shown already by Schmeidler \cite[\S\,2]{Schmeidler72}
that the exact games involve traditional {\em convex\/} (= {\em supermodular\/}) games \cite{Shapley72} and
form a subclass of a popular class of {\em totally balanced\/} games \cite[Chapter\,V]{Rosenmuller00}.
The latter class is then included in the class of {\em balanced\/} games \cite{Shapley67},
defined as the class of cooperative games with non-empty {\em core\/} polyhedron.
Note that all these game-theoretical concepts have also their counterparts in the context of imprecise probabilities;
see \cite{MM18} for more details about the correspondence.

In our framework of a finite set $N$ of players, one can consider the geometric point of view on the situation.
It follows from the results in \cite{CHK11,LBH12} that the set of (characteristic functions for) exact games over a fixed set $N$
of players forms a polyhedral cone; the same is true for other three classes of cooperative games mentioned above.
Thus, exact games over $N$ can be characterized by means of a finite number of linear
inequalities and a natural question of theoretical interest is what is the minimal set of such inequalities.
Note in this context that we do not distinguish between an inequality and its multiple by a positive factor and
that the uniqueness of the minimal defining set of linear inequalities for a polyhedral set is relative
to the affine (= a shifted linear) space generated by the set;
these are the so-called {\em facet-defining\/} inequalities for the polyhedral set.

Let us mention other related results. Balanced games were characterized in terms of linear
inequalities already in the 1960s \cite{Bondareva63,Shapley67} and the facet-defining inequalities
for this cone correspond to certain systems of subsets of $N$, called the {\em minimal balanced\/} set
systems (= collections). The facet-defining inequalities for the cone of supermodular games
also correspond to special set systems, namely to those consisting of two sets $S,T\subseteq N$ such that
$S\setminus T$ and $T\setminus S$ are singletons \cite[Corollary\,11]{KVV10}.
The facet-defining inequalities for the cone of totally balanced games were recently characterized
in \cite{KS19}: they also correspond to set systems, called {\em irreducible\/} (minimal) {\em balanced\/}
systems on {\em subsets\/} of $N$.

The existence of a finite system of linear inequalities characterizing exact games
follows from the results by Lohmann {\em et al.\/} \cite{LBH12}. The inequalities reported in \cite{LBH12}
also correspond to systems of subsets of $N$; nevertheless, not all
these inequalities are facet-defining for the exact cone. The problem of the characterization
of facet-defining inequalities for this cone was then discussed in \cite[\S\,6]{KS19}
where a conjecture has been raised about their form. That conjecture, confirmed in case $|N|\leq 5$,
has an equivalent formulation saying that a game is exact iff it is totally balanced and its {\em anti-dual\/} game
\cite{Oishi16} is totally balanced as well.

Recall that each polyhedral cone can be characterized by means of (finitely many)
linear inequalities, which
characterization is named the {\em outer description}; nonetheless, each polyhedral cone can
alternatively be defined as the conic hull of finitely many vectors, which characterization is named the
{\em inner description}. The latter approach leads to the task to characterize the extreme rays of (a pointed
version of) the polyhedral cone. 
The generators of the extreme rays of the cones of (suitably) standardized cooperative games are named {\em extreme games\/} \cite[\S\,V.4]{Rosenmuller00}.
Note in this context that the inner description for the cone of balanced games was presented in \cite[\S\,5.2]{KS19}.
On the other hand, the number of extreme rays for the other three (standardized) cones
seems to grow more than exponentially with $|N|$; this observation decreases the hope that they have manageable
inner description. The available results here are the criteria to recognize
extreme games: we have proposed simple and easily implementable such linear criteria
based on the (vertices of the) corresponding core polyhedron both in the supermodular case \cite{SK16} and in case of exact games \cite{SK18}.

\subsection{Main results in this paper}
This paper is devoted to the problem of characterization of facet-defining inequalities
for the cone of exact games. We follow the line of research indicated above.
More specifically, we introduce the concept of a ({\em minimal\/}) {\em semi-balanced\/} system of subsets of $N$,
which generalizes the classic concept of a (minimal) balanced set system on $N$ from \cite{Shapley67}.
Linear inequalities assigned to these set systems are shown to characterize the cone of exact games (see Corollary~\ref{cor.exa-char}).
This result is analogous to \cite[Theorem~3.4]{LBH12} in which exact games are characterized by means of
the so-called ``exact balanced collections" of subsets of $N$, but there is one important difference.
Our semi-balanced set systems technically differ from the collections of sets introduced by Lohmann {\em et al.\/} \cite{LBH12},
although the assigned sets of inequalities (as a whole) are necessarily equivalent.

The point is that our concept of a semi-balanced set system allows one to recognize easily certain hidden
symmetry. More specifically, every game over $N$ is exact iff its anti-dual game is exact and this fact is reflected
in the  linear inequalities for exact games: an inequality is valid/facet-defining for the exact cone iff
the same holds for its {\em conjugate\/} inequality (see Section~\ref{ssec.complementarity}).
Each semi-balanced set system $\calS$ has a {\em complementary\/} semi-balanced set system $\calS^{\refl}$ and the inequality
assigned to $\calS^{\refl}$ is conjugate to the inequality assigned to $\calS$.
We also introduce basic classification for minimal semi-balanced set systems, called briefly {\em min-semi-balanced\/}
systems (on $N$), into four basic classes; the class of minimal balanced (= {\em min-balanced\/})
set systems on $N$ is one of them (Section~\ref{ssec.classification}).
In addition to that we introduce pictorial representatives for (permutational types of) min-semi-balanced systems
(Section~\ref{ssec.picto-grams}), which easily encode the assigned inequalities and reflect both the classification
and complementarity relationships.

Besides that we establish a certain one-to-many correspondence between min-balanced systems on $N$ involving at least 3 sets and {\em purely\/} min-semi-balanced set systems on $N$, that is, those that are not balanced (Section~\ref{ssec.get-purely}).
This correspondence may be a basis for a procedure to generate the complete list of min-semi-balanced systems on $N$
on basis of the list of all min-balanced systems on $N$.
The point is that if $|N|\geq 3$ then every facet-defining inequality for the exact cone corresponds to a
purely min-semi-balanced system.
Note in this context that an analogous observation that the min-balanced systems are not needed for the description
of the exact cone was already made in \cite[\S\,5,\,Theorem\,5.4]{LBH12}.

Nonetheless, even purely min-semi-balanced set systems can be superfluous in the sense that the assigned inequalities
are not facet-defining. We introduce a narrower concept of an {\em indecomposable\/} min-semi-balanced
set system (see Section~\ref{sec.indecomposable-min}) and show that if $|N|\geq 3$ then any facet-defining inequality
corresponds to an indecomposable set system. In fact, our main result is that the facet-defining inequalities for
the exact cone are just those assigned to indecomposable min-semi-balanced set systems (see Theorem~\ref{thm.main-result}).

We also put more light on the relation of the cones of exact and totally balanced games.
Specifically, we first give a counter-example to a former conjecture from \cite[\S\,6]{KS19}
mentioned above in Section~\ref{ssec.overview} (Section~\ref{ssec.counter-example}).
Then we derive, as a consequence of our main result, that
those facet-defining inequalities for the cone of totally balanced games which concern {\em strict
subsets\/} of $N$ are also facet-defining for the cone of exact games (Section~\ref{ssec.sub-balanced}).
In fact, these inequalities correspond to min-semi-balanced systems from one of four basic classes
in our classification (from Section~\ref{ssec.classification}).
In particular, every irreducible min-balanced system on a strict subset of $N$ from \cite{KS19} can be extended
uniquely to a certain special indecomposable min-semi-balanced set system (on $N$).

\subsection{Structure of the paper}
In Section~\ref{sec.preliminaries} we recall elementary concepts and basic facts.
Our concept of a {\em semi-balanced\/} system is introduced in Section~\ref{sec.semi-balanced}.
We give several equivalent definitions of a minimal semi-balanced system, called a {\em min-semi-balanced\/} system,
there and introduce the linear inequalities assigned to these set systems.
The cone of exact games is characterized by means of these inequalities in Section~\ref{sec.char-exact}.
In Section~\ref{sec.properties} we then introduce the concept of a complementary set system and basic classification
of min-semi-balanced systems. We also propose special pictures to represents these set systems
and the inequalities assigned to them.
In Section~\ref{sec.purely-min} we establish the correspondence of purely min-semi-balanced systems to min-balanced ones
and shown that the inequalities assigned to min-balanced systems are superfluous if $|N|\geq 3$.
The concept of an {\em indecomposable\/} min-semi-balanced set system is then defined formally in Section~\ref{sec.indecomposable-min}.
We formulate our main result there saying that indecomposability is a necessary and sufficient condition for the assigned inequality to be facet-defining.
Section~\ref{sec.relation} is devoted to the relation of the exact cone to the cone of totally balanced games.
We first give a counter-example to the conjecture from \cite[\S\,6]{KS19}.
Then we characterize those facet-defining inequalities for the totally balanced cone which are also facet-defining
for the exact cone. In Conclusions (Section~\ref{sec.coclusion}) we summarize our findings and give a reference to our catalogue of indecomposable min-semi-balanced systems over sets of low cardinality. Some of the longer technical proofs are moved to the Appendix.

\section{Preliminaries}\label{sec.preliminaries}
Throughout the paper the symbol $N$ will denote a finite set of {\em players} and we restrict ourselves to the non-degenerate case $|N|\geq 2$.
The power set $\caP := \{ S : S\subseteq N\}$ of the set of players is the set of {\em coalitions}.
The symbol $\subset\/$ will denote strict inclusion of either sets or set systems, that is, $S\subset T$ iff $S\subseteq T$ and $S\neq T$.
The set of real numbers will be denoted by ${\dv R}$, the set of rational numbers by ${\dv Q}$.

\subsection{Basic versions of linear combinations}\label{ssec.linear}
We are going to deal with vectors in real Euclidean spaces ${\dv R}^{\calI}$, where $\calI$ is
a non-empty finite index set. Our elementary linear algebraic operations will concern the space ${\dv R}^{N}$
in which case one has $\calI =N$. But later on, some more advanced geometric considerations
will be in space ${\dv R}^{\calI}$ where $\calI$ will be a class of subsets of $N$, mostly $\calI=\caP$.

The {\em incidence vector\/} of a coalition (= set) $S\subseteq N$ will be denoted by $\chi_{S}\in {\dv R}^{N}$:
$$
\chi_{S}(i):=
\left\{
\begin{array}{cl}
1 & \mbox{if $i\in S$,}\\
0 & \mbox{if $i\not\in S$,}
\end{array}
\right.
\qquad \mbox{for $i\in N$.}
$$
By a {\em constant vector\/} in ${\dv R}^{N}$ will be called a vector whose components
equal each other, that is, a vector of the form $[r,\ldots,r]\in {\dv R}^{N}$ where $r\in {\dv R}$.
A special case of a constant vector in ${\dv R}^{N}$ is the {\em zero vector\/} in ${\dv R}^{N}$, denoted by $\0$.
A finite set ${\sf S}$ of vectors in ${\dv R}^{N}$ is {\em linearly independent} if
$\sum_{x\in\, {\sf S}} \alpha_{x}\cdot x=\0$ with $\alpha_{x}\in {\dv R}$, $x\in {\sf S}$,
implies  $[\,\forall\, x\in {\sf S} ~~\alpha_{x}=0\,]$, otherwise it is called {\em linearly dependent}.
Analogously, a finite set ${\sf S}\subseteq {\dv R}^{N}$ is
{\em affinely independent} if
$$
\forall\, \alpha_{x}\in {\dv R},\, x\in {\sf S},\qquad
[\,\,\sum_{x\in\, {\sf S}} \alpha_{x}=0 ~~\&~~ \sum_{x\in\, {\sf S}} \alpha_{x}\cdot x=\0\,\,] ~\Rightarrow~
[\,\forall\, x\in {\sf S} ~~\alpha_{x}=0\,]\,,
$$
otherwise it is called {\em affinely dependent}.
\smallskip

Other elementary concepts apply to any real Euclidean space ${\dv R}^{\calI}$, $0<|\calI|<\infty$.
We are going to use the symbol
$$
\langle \theta,x\rangle ~:=~ \sum_{\iota\in\,\calI}\, \theta(\iota)\cdot x(\iota)\qquad\mbox{for
vectors $\theta,x\in{\dv R}^{\calI}$}
$$
to denote the respective scalar product in ${\dv R}^{\calI}$.
A (finite) linear combination $\sum_{i\in I} \lambda_{i}\cdot x_{i}$
of vectors $x_{i}\in {\dv R}^{\calI}$ with real coefficients $\lambda_{i}\in {\dv R}$ is called
\begin{itemize}
\item {\em non-zero\/} if there is $i\in I$ with $\lambda_{i}\neq 0$,
\item {\em affine\/} if\/ $\sum_{i\in I} \lambda_{i}=1$,
\item {\em conic\/} if $\lambda_{i}\geq 0$ for all $i\in I$, and
\item {\em convex\/} if it is both affine and conic.
\end{itemize}
The {\em convex hull\/} of a set ${\sf S}\subseteq {\dv R}^{\calI}$
is the collection of all convex combinations of vectors from ${\sf S}$;
it will be denoted by $\conv ({\sf S})$. A set ${\sf S}\subseteq {\dv R}^{\calI}$ is
{\em convex\/} if it is closed under convex combinations: ${\sf S}=\conv ({\sf S})$.
The {\em conic hull\/} of ${\sf S}\subseteq {\dv R}^{\calI}$
is the set of all conic combinations of vectors from ${\sf S}$;
it will be denoted by $\cone ({\sf S})$.

Analogously, the {\em affine hull\/} of a set ${\sf S}\subseteq {\dv R}^{\calI}$
is the collection of all affine combinations of vectors from ${\sf S}$.
It is always an {\em affine subspace\/} of ${\dv R}^{\calI}$, that is,
a subset ${\sf A}\subseteq {\dv R}^{\calI}$ closed under affine combinations.
A non-empty affine subspace is always 
a shifted linear subspace of ${\dv R}^{\calI}$,
that is, a set of the form ${\sf A}=x+{\sf L}:=\{\, x+y\,:\ y\in {\sf L}\,\}$, where $x\in {\dv R}^{\calI}$ and
${\sf L}\subseteq {\dv R}^{\calI}$ is a linear subspace \cite[\S\,1]{Brondsted83}; ${\sf L}$
is then uniquely determined by ${\sf A}$ while $x$ not.
The {\em dimension\/} of a set ${\sf S}\subseteq {\dv R}^{\calI}$, denoted by $\dim ({\sf S})$, is the
dimension of its affine hull ${\sf A}$, defined as the dimension of ${\sf L}$.
A {\em hyperplane\/} in ${\dv R}^{\calI}$ is an affine subspace ${\sf H}$ of ${\dv R}^{\calI}$
of the dimension $|\calI|-1$. An equivalent condition is that it is the set of solutions
$x\in {\dv R}^{\calI}$ to the equation $\langle \theta,x\rangle=\beta$, where $\beta\in {\dv R}$
and $\theta\in {\dv R}^{\calI}$ is a {\sf non-zero\/} vector in ${\dv R}^{\calI}$; see \cite[\S\,1]{Brondsted83}.

\subsection{Some concepts from polyhedral geometry}\label{ssec.polyhedrals}
Throughout the paper we assume that the reader is familiar with standard concepts and basic facts from polyhedral geometry; see \cite{Brondsted83,Schrijver98,BK92,Ziegler95}, for example. Nevertheless, for reader's convenience,
we recall those of them that are used (repeatedly) in our paper.
\smallskip

Given distinct $x,y\in {\dv R}^{\calI}$, the convex hull of $\{x,y\}$
is the {\em closed segment}, denoted by $[x,y]$, while the {\em open segment\/}, denoted by $]x,y[$,
consists of convex combinations of $x$ and $y$ which have both coefficients non-zero:
$$
]x,y[ ~:=~ \{\, (1-\alpha)\cdot x+\alpha\cdot y\,:\ 0<\alpha<1\,\}\,.
$$
A {\em polyhedron\/} in ${\dv R}^{\calI}$, $0<|\calI|<\infty$, is the
set of vectors $x\in {\dv R}^{\calI}$ specified by finitely many linear
inequalities $\langle \theta,x\rangle\geq \beta$ for $x\in {\dv R}^{\calI}$, where $\theta\in {\dv R}^{\calI}$
and $\beta\in {\dv R}$. A polyhedron is called {\em rational\/} if, moreover, $\theta\in {\dv Q}^{\calI}$
and $\beta\in {\dv Q}$ for these inequalities.
A set ${\sf S}\subseteq {\dv R}^{\calI}$ is {\em bounded\/} if
there are constants $c_{0},c_{1}\in {\dv R}$ such that $c_{0}\leq x(\iota)\leq c_{1}$
for any component $x(\iota)$, $\iota\in\calI$, of any $x\in {\sf S}$.
A {\em polytope\/} (in ${\dv R}^{\calI}$) is the convex hull
of a non-empty finite set of vectors in ${\dv R}^{\calI}$.
A fundamental result in polyhedral geometry says that a subset of ${\dv R}^{\calI}$
is a polytope iff it is a non-empty {\em bounded polyhedron}; see \cite[Theorem\,9.2]{Brondsted83} or \cite[Theorem\,2.15]{Ziegler95}.

A {\em face\/} of a polyhedron ${\sf P}\subseteq {\dv R}^{\calI}$, more precisely an {\em exposed face\/} of
${\sf P}$ \cite[\S\,5]{Brondsted83}, is a subset ${\sf F}\subseteq {\sf P}$ consisting
of vectors $x\in {\sf P}$ satisfying $\langle \theta,x\rangle=\beta$ for some
$\theta\in {\dv R}^{\calI}$, $\beta\in {\dv R}$, such that $\langle \theta,y\rangle\geq\beta$ is a valid inequality for all $y\in {\sf P}$. In case of a {\sf polyhedron\/} ${\sf P}$, this is equivalent to
the condition that ${\sf F}\subseteq {\sf P}$ is a convex subset of it such that
one has $[y,z]\subseteq {\sf F}$ whenever $y,z\in{\sf P}$ and $]y,z[\,\cap\, {\sf F}\neq\emptyset$; use \cite[\S\,8]{Brondsted83} or see \cite[Theorem\,7.51]{BK92}.
The number of faces of a polyhedron is finite; see \cite[Corollary\,8.5]{Brondsted83}. The {\em face-lattice\/} of ${\sf P}\subseteq {\dv R}^{\calI}$ is the collection of its faces, ordered by inclusion relation $\subseteq$; it is indeed a lattice in usual sense \cite[\S\,5]{Brondsted83}.

{\em Facet\/} of a non-empty polyhedron ${\sf P}$ is its face of the dimension $\dim ({\sf P})-1$. An equivalent definition is that facet is a maximal face ${\sf F}$ of ${\sf P}$ distinct from ${\sf P}$; use \cite[Corollary\,8.6]{Brondsted83}.
A basic fact is that each {\em full-dimensional\/} proper polyhedron ${\sf P}\subset {\dv R}^{\calI}$, that is,
each polyhedron with $\dim ({\sf P})=|\calI|$ and ${\sf P}\neq {\dv R}^{\calI}$, is specified by those valid inequalities for ${\sf P}$
which define facets and the specification of\/ ${\sf P}$ by {\em facet-defining\/} inequalities is the
{\sf unique\/} inclusion-minimal inequality description of ${\sf P}$ (up to
positive multiples of inequalities); see \cite[Theorem\,8.2]{Brondsted83}.

A {\em vertex\/} (= an extreme point) of a convex set ${\sf P}\subseteq {\dv R}^{\calI}$ is a vector $x\in {\sf P}$ such that there is {\sf no open segment\/} $]y,z[$ with $y,z\in {\sf P}$ and $x\in\,]y,z[$, which is, in case of a polytope ${\sf P}$, another way of saying that $\{x\}$ is a face of ${\sf P}$ (of the dimension $0$). The set of vertices of (a polytope) ${\sf P}$ will be denoted by $\ext ({\sf P})$.
A well-known consequence of famous Krein-Milman theorem 
is that every polytope ${\sf P}$ has finitely many vertices and equals to the convex  hull of
the vertex set:
${\sf P}= \conv (\ext ({\sf P}))$; see \cite[Theorem\,7.2(c)]{Brondsted83} or \cite[Proposition\,2.2(i)]{Ziegler95}.

One of easy observations is that if ${\sf S}=\conv ({\sf T})$ for ${\sf S},{\sf T}\subseteq {\dv R}^{\calI}$ then $\ext ({\sf S})\subseteq {\sf T}$.
Another immediate observation is that if ${\sf Q}\subseteq {\sf P}$ are polytopes
in ${\dv R}^{\calI}$ then $x\in {\sf Q}\cap \ext ({\sf P})$ implies $x\in \ext ({\sf Q})$.
Every non-empty face ${\sf F}$ of a polytope ${\sf P}$ is again a polytope and
$\ext ({\sf F})={\sf F}\cap \ext ({\sf P})$; see \cite[Theorem\,7.3]{Brondsted83} or
\cite[Proposition\,2.3(i)]{Ziegler95}.
Further basic fact, which follows from properties of the operator ${\sf S}\mapsto\conv({\sf S})$
and Krein-Milman theorem, is this:
${\sf Q}:=\conv(\bigcup_{j\in J} {\sf P}_{j})$, with finite $J$, is a bounded polyhedron
whenever ${\sf P}_{j}$, $j\in J$, are bounded polyhedrons in ${\dv R}^{\cal I}$.
%

An {\em edge\/} of a polytope ${\sf P}$ is a closed segment $[y,z]\subseteq {\sf P}$ which is a face of
${\sf P}$ (of the dimension $1$); then necessarily $y,z\in\ext ({\sf P})$. Further special fact is as follows.
Given a polytope ${\sf P}\subseteq {\dv R}^{\calI}$ and a hyperplane ${\sf H}$ in ${\dv R}^{\calI}$ with
${\sf H}\cap {\sf P}\neq\emptyset$ the intersection ${\sf Q}:={\sf H}\cap {\sf P}$ is a polytope and
$x\in \ext ({\sf Q})$ iff {\sf either\/} $x\in {\sf H}\cap\ext ({\sf P})$ {\sf or\/} there is an edge
$[y,z]$ of ${\sf P}$ such that $x\in\,]y,z[$ and $[y,z]\cap{\sf H}=\{x\}$; use \cite[Theorem\,11.1(d)]{Brondsted83}.
\smallskip

A {\em polyhedral cone\/} in\/ ${\dv R}^{\calI}$ is a subset ${\sf C}$ of\/ ${\dv R}^{\calI}$ defined as the
conic hull of a non-empty finite set ${\sf S}$ of vectors from ${\dv R}^{\calI}$.
An equivalent definition 
is that ${\sf C}$ is specified by finitely many inequalities
$\langle \theta,x\rangle\geq 0$ for $x\in {\dv R}^{\calI}$ \cite[Theorem\,1.3]{Ziegler95}; thus, it is a non-empty polyhedron.
A polyhedral cone ${\sf C}\subseteq {\dv R}^{\calI}$ is called {\em pointed}, if $-{\sf C}\cap{\sf C}=\{\bf0\}$,
where $-{\sf C}:=\{-y\,:\ y\in {\sf C}\,\}$ and $\bf0$ denotes the zero vector in ${\dv R}^{\calI}$.
An equivalent condition is that there exists (non-zero) $\theta\in {\dv R}^{\calI}$ such that $\langle \theta,x\rangle>0$
for any $x\in {\sf C}\setminus\{\bf0\}$; see \cite[Proposition\,2]{Studeny93}.
It makes no problem to observe that if, moreover, ${\sf C}\setminus\{\bf0\}\neq\emptyset$, then, for each $\beta>0$, the
intersection of ${\sf C}$ with the hyperplane ${\sf H}:=\{ x\in {\dv R}^{\calI}\,:\
\langle \theta,x\rangle =\beta\,\}$ is a polytope; this is because one can find finite $\emptyset\neq {\sf S}\subset {\sf H}$
with ${\sf C}=\cone ({\sf S})$. The reader can verify (using alternate definitions of a face) that then
the mapping ${\sf F}\subseteq {\sf C} ~\mapsto~{\sf F}\cap {\sf H}$ established a one-to-one correspondence between
{\sf non-empty\/} faces ${\sf F}$ of ${\sf C}$ and (all) faces of the polytope ${\sf P}:={\sf C}\cap{\sf H}$: the inverse mapping
is ${\sf F}^{\prime}\subseteq {\sf C}\cap {\sf H} ~\mapsto~ \cone ({\sf F}^{\prime}\cup\{{\bf0}\})
\equiv \{{\bf0}\}\cup \{\alpha\cdot x\,:\ \alpha\geq 0 ~\&~ x\in{\sf F}^{\prime}\,\}$.
This correspondence preserves the inclusion ordering; thus, it is an {\em isomorphism\/} between the lattice of non-empty faces of
${\sf C}$ and the face-lattice of ${\sf P}$.

Every subset ${\sf S}\subseteq {\dv R}^{\calI}$ can be assigned its {\em dual cone\/}
$$
{\sf S}^{\dual}:=\{\, \theta\in {\dv R}^{\calI}\,:\ \langle\theta,x\rangle\geq 0\quad
\mbox{for any $x\in {\sf S}$}\,\}\,, ~~\mbox{which is clearly a closed convex cone.}
$$
A well-known elementary fact is that ${\sf C}\subseteq {\dv R}^{\calI}$ is a non-empty closed convex cone
iff ${\sf C}={\sf C}^{\dual\dual}$, which happens iff ${\sf C}={\sf S}^{\dual}$ for
some ${\sf S}\subseteq {\dv R}^{\calI}$; see for example \cite[Consequence~1]{Studeny93}.
Moreover, dual cone to a polyhedral cone is also a polyhedral cone; use \cite[Corollary\,7.12]{BK92}.
Thus, if one shows, for a polyhedral cone ${\sf C}\subseteq {\dv R}^{\calI}$,
and for a set ${\sf D}\subseteq {\dv R}^{\calI}$ that ${\sf D}={\sf C}^{\dual}$ then this
fact already implies that ${\sf C}$ and ${\sf D}$ are mutually dual polyhedral cones,
which means both ${\sf C}={\sf D}^{\dual}$ and ${\sf D}={\sf C}^{\dual}$.
Another useful fact from polyhedral geometry
is that the lattices of non-empty faces of mutually dual polyhedral cones are {\em anti-isomorphic}. More specifically, the mapping is as follows:
$$
\mbox{${\sf F}\subseteq {\sf C}$ a non-empty face of ${\sf C}$} ~\mapsto~~
{\sf F}^{\perp} := \{\, \theta\in {\sf D}\,:\ \langle\theta,x\rangle =0~~
\mbox{for all $x\in {\sf F}$}\,\}
$$
and the inverse mapping is of the same form (exchange ${\sf C}$ for ${\sf D}$); use \cite[Theorem\,7.41]{BK92} to derive that. This particular one-to-one correspondence reverses the inclusion ordering and has the property that
$\dim({\sf F})+\dim ({\sf F}^{\perp})=|\calI|$; use \cite[Theorem\,7.42]{BK92}.

\subsection{Some concepts from game theory}\label{ssec.game-concepts}
A ({\em transferable utility coalitional\,}) {\em game} over (a set of players) $N$ is modeled by a real function
$\m\colon \mathcal{P}(N)\to {\dv R}$ such that $\m(\emptyset)=0$, called the ``characteristic function" of the game.
The class of all games over $N$ will be denoted by ${\cal G}(N)$.
Given $m\in {\cal G}(N)$ and $\emptyset\neq S\subseteq N$, the restriction $\m_{S}$ of $\m$ to $\mathcal{P}(S)$
is called a {\em subgame} of the game $m\in {\cal G}(N)$.

The {\em core $C(\m)$ of\/} a game $\m\in{\cal G}(N)$ is the polyhedron
$$
C(\m) := \{\,[x_{i}]_{i\in N}\in {\dv R}^{N} \,:\ \sum_{i\in N} x_{i}=\m(N) ~\&~
\sum_{i\in S}x_{i} \geq \m(S)\; \mbox{for all $S\subseteq N$}\}\,.
$$
We say that a game $\m\in{\cal G}(N)$ is {\em balanced} if it has a non-empty core,
\begin{itemize}
\item {\em totally balanced} if every subgame of $\m$ is balanced,
\item {\em exact} if, for each coalition $S\subseteq N$, there exists a vector
$[x_{i}]_{i\in N} \in C(\m)$ in the core that is {\em tight for $S$}, which means that $\sum_{i\in S}x_{i} = \m(S)$.
\end{itemize}
The set $\mathcal{T}(N)$ of all totally balanced games over $N$ is known to be a polyhedral cone in ${\dv R}^{\caP}$; the same is true for the set $\mathcal{E}(N)$ of all exact games over $N$; see \cite[\S\,2]{KS19}.

Recall from \cite[Theorem\,1]{KZ82} that $\m$ is totally balanced iff it has a {\em finite min-representation}, which means that there exists a non-empty finite ${\cal X}\subseteq {\dv R}^{N}$ such that
$$
\m(S)~=~\min_{x\in {\cal X}}\, \sum_{i\in S} x_{i}\qquad \mbox{for any $S\subseteq N$}.
$$
It is a well-known fact that $\m$ is exact iff it has a {\em min-representation\/}
$\emptyset\neq {\cal X}\subseteq C(\m)$ {\em consisting of the elements in the core};
see \cite[Proposition\,1]{SK18}, for example.
Hence, we know that every exact game is totally balanced.

Following \cite{Oishi16}, given $\m\in {\cal E}(N)$, by the {\em anti-dual\/} of $\m$ we call the game
$$
\am(S) := \m(N\setminus S)-\m(N) \quad \mbox{for all $S\subseteq N$.}
$$
Note that, by habitual terminology in cooperative game theory, the $(-1)$ multiple of $\am$
is named the dual game of $\m$; see \cite[Definition 6.6.3]{PelegSudholter07}.
Nonetheless, for our purpose the concept of an anti-dual is more suitable as
$\m\in \mathcal{E}(N)$ iff $\am\in \mathcal{E}(N)$, see \cite[\S\,3.2]{KS19}.
In particular, if $\m$ is exact then both $\m$ and $\am$ are totally balanced.
On the other hand, $\m\in\mathcal{E}(N) ~\not\Rightarrow~ -\m\in\mathcal{E}(N)$ in general; thus,
duals of exact games need not be exact.

\section{Balanced and semi-balanced set systems}\label{sec.semi-balanced}

We first recall a classic well-known concept in cooperative game theory \cite{Shapley67}.

\begin{defin}[non-trivial and balanced set system]~\label{def.balanced}\rm \\
Assume that $N$ is a finite set with $|N|\geq 2$. By a {\em non-trivial\/} set system on $N$ will be meant
a system $\calS$ of its subsets such that $\emptyset\neq\calS\subseteq {\caP}\setminus \{\emptyset ,N\}$.\\[0.1ex]
A non-trivial set system $\calB$ is called {\em balanced\/} (on $N$) if $\chi_{N}$ is a {\sf conic\/} combination of vectors
$\{ \chi_{S}\,:\, S\in\calB\}$ with {\sf all coefficients non-zero}.
Such a system $\calB$ is called {\em minimal\/} if there is no balanced system $\calC$ on $N$ with $\calC\subset\calB$; $\calB$ will then be called briefly {\em min-balanced\/} (on $N$).
\end{defin}

In other words, $\calB$ is balanced on $N$ if $\chi_{N}$ is a linear combination of vectors $\{ \chi_{S}\,:\, S\in\calB\}$ with strictly positive coefficients. In the sequel we are going to generalize this concept.

\subsection{Semi-balanced set systems}\label{ssec.semi-balanced}
Our generalization is based on the following elementary concept.

\begin{defin}[semi-conic combination]~\rm \\
We shall say that a linear combination $\sum_{i\in I} \lambda_{i}\cdot x_{i}$ in ${\dv R}^{N}$
is {\em semi-conic\/} if {\sf at most one\/} of its coefficients is strictly {\sf negative}, that is,
$|\{j\in I\,:\  \lambda_{j}<0\}|\leq 1$.
\end{defin}

While linear combination concepts recalled in Section~\ref{ssec.linear} are standard in mathematics, the concept of a semi-conic combination is
a specific concept relevant to the topic of our study. Our terminology is motivated by the fact that the
remaining coefficients $\lambda_{i}$, $i\neq j$, in a such a linear combination must be non-negative.
Thus, any conic combination is also semi-conic. Nonetheless, despite this fact, the concepts of semi-conic
and conic combination differ substantially from each other.

\begin{remark}\em
In this side note we explain the principal difference between conic and semi-conic combinations. Recall that the
{\em conic hull\/} of a set ${\sf S}$ in a real Euclidean space is the collection $\cone ({\sf S})$ of all conic combinations of vectors from ${\sf S}$ and it is always a convex cone. Another well-known fact is that the mapping ${\sf S}\mapsto \cone({\sf S})$ is a closure operator in sense of abstract algebra \cite[\S\,V.1]{Birkhoff95}; this means that it is {\em extensive} [${\sf S}\subseteq \cone({\sf S})$], {\em monotone} [${\sf S}\subseteq {\sf T}$ implies $\cone ({\sf S})\subseteq \cone({\sf T})$] and {\em idempotent} [$\cone (\cone ({\sf S}))=\cone({\sf S})$].
One can analogously introduce the {\em semi-conic hull\/} $\scon ({\sf S})$ as the set of all semi-conic combinations of vectors from ${\sf S}$, but this set need not be convex. The operator ${\sf S}\mapsto \scon({\sf S})$ is then extensive and monotone but it is not idempotent. Consider, for example, the case ${\sf S}=\{ (1,0),(0,1)\}\subseteq {\dv R}^{2}$.
Thus, the mapping ${\sf S}\mapsto \scon ({\sf S})$ is not a closure operator.
\end{remark}

Now, we are ready to introduce basic concepts in our treatise.

\begin{defin}[semi-balanced set system, exceptional set]~\label{def.semi-bal}\rm \\
Assume that $N$ is a finite set with $|N|\geq 2$. 
We shall say that a non-trivial set system $\calS$ on $N$ is
{\em semi-balanced\/} (on $N$) if there is a {\sf constant vector\/} in ${\dv R}^{N}$ which is a {\sf semi-conic\/} combination of vectors $\{ \chi_{S}\,:\ S\in\calS\}$ with {\sf all coefficients non-zero}.
A semi-balanced system $\calS$ on $N$ will be called {\em minimal\/} if there is no semi-balanced system $\calC$ on $N$ with $\calC\subset\calS$. We will then say briefly that such a set system is {\em min-semi-balanced\/} (on $N$).\\[0.2ex]
A semi-balanced set system $\calS$ (on $N$) which is {\sf not balanced\/} (on $N$) will be called {\em purely semi-balanced\/} (on $N$). Analogously, by a {\em purely min-semi-balanced\/} system will be meant a min-semi-balanced system which is not balanced.\\[0.2ex]
Given a non-trivial set system $\calS$ on $N$, we will say that a set $T\in\calS$ is {\em exceptional\/} within $\calS$ if there exists a linear combination $\sum_{S\in\calS} \lambda_{S}\cdot\chi_{S}$
yielding a constant vector in ${\dv R}^{N}$ with $\lambda_{T}<0$ and $\lambda_{S}\geq 0$ for $S\in\calS\setminus\{T\}$.
\end{defin}


It follows directly from the definition that any balanced system is also semi-balanced.
Let us emphasize that both these concepts are relative to $N$ despite one may have $\bigcup\calS\subset N$ for a semi-balanced system $\calS$ on $N$ (see Example~\ref{ex.non-union} below); this is because the constant vector is required to be in ${\dv R}^{N}$. Note that, for a semi-balanced system, all coefficients in any of the considered
linear combinations must be strictly positive with one possible exception of a strictly negative coefficient (with an exceptional set).
\smallskip

If a set $T\subset N$ is exceptional within a system $\calS$
then it is exceptional within any larger non-trivial system $\calT\supseteq\calS$: put $\lambda_{S}:=0$ for $S\in\calT\setminus\calS$.
Of course, every non-trivial set system admitting an exceptional set contains a semi-balanced system because the considered linear combination is semi-conic.
Conversely, every purely semi-balanced system has at least one exceptional set. 
Note that, in case of a (purely) min-semi-balanced system, this exceptional set is uniquely determined, which fact follows from later  Lemmas~\ref{lem.semi-comb} and \ref{lem.min-semi-bal}.
\smallskip

A set system containing a semi-balanced system need not be semi-balanced.
In fact, even the union of two semi-balanced systems need not be a semi-balanced system as the next
Example~\ref{ex.non-union} shows.
On the other hand, the union of a semi-balanced system and a balanced system has to be a semi-balanced system:
consider a non-trivial convex combination of the respective  semi-conic combinations yielding constant vectors.
The same argument implies that the union of two balanced systems on $N$ is a balanced system on $N$.

\begin{example}\label{ex.non-union}\rm
Given $N:=\{a,b,c,d\}$, take set systems $\calS:=\{\,a,b,ab\,\}$ and $\calT:=\{\,c,bc,acd\,\}$. The equalities
$1\cdot\chi_{a}+1\cdot\chi_{b}+(-1)\cdot\chi_{ab}=\0$ and $(-1)\cdot\chi_{c}+1\cdot\chi_{bc}+1\cdot\chi_{acd}=\chi_{N}$
imply that both $\calS$ and $\calT$ is semi-balanced on $N$. To show that their union $\calD:=\calS\cup\calT$ is not semi-balanced on $N$
consider a linear combination $\sum_{S\in\calD}\lambda_{S}\cdot\chi_{S}=[r,r,r,r]\in {\dv R}^{N}$
having all coefficients non-zero. Realize that $\lambda_{acd}=\sum_{S\in\calD:\, d\in S}\lambda_{S}=r$. The equality
$\lambda_{acd}=r=\sum_{S\in\calD:\, c\in S}\lambda_{S}=\lambda_{c}+\lambda_{bc}+\lambda_{acd}$ then gives $\lambda_{c}+\lambda_{bc}=0$
and the assumption $\lambda_{c}\neq 0\neq\lambda_{bc}$ implies that [\,$\lambda_{c}<0$ or $\lambda_{bc}<0$\,].
Analogously, $\lambda_{acd}=r=\sum_{S\in\calD:\, a\in S}\lambda_{S}=\lambda_{a}+\lambda_{ab}+\lambda_{acd}$ means $\lambda_{a}+\lambda_{ab}=0$,
which implies that [\,$\lambda_{a}<0$ or $\lambda_{ab}<0$\,]. Therefore, one has $\lambda_{T}<0$ for at least two sets $T\in\calD$
and the considered linear combination is not semi-conic.
\end{example}

The following lemma contains a few elementary observations valid for semi-balanced set systems; in fact, they hold for a wider class of systems containing semi-balanced ones.

\begin{lem}\label{lem.semi-comb}\em
Given $|N|\geq 2$ and $\emptyset\neq \calS\subseteq {\caP}\setminus \{\emptyset,N\}$, let $\sum_{S\in\calS} \lambda_{S}\cdot\chi_{S}=\rho$ be a {\sf non-zero\/} semi-conic combination yielding a constant vector
$\rho=[r,\ldots,r]\in {\dv R}^{N}$ (with zero coefficients allowed).
Then one has $\sum_{S\in\calS} \lambda_{S}\geq r\geq 0$; moreover, $r>0$ in case of a conic combination.\\
In any case $\sum_{S\in\calS} \lambda_{S}>0$ and by a positive factor multiplication one gets
an {\sf affine\/} semi-conic combination $\sum_{S\in\calS} \lambda_{S}\cdot\chi_{S}$ yielding a constant
vector $\rho=[r,\ldots,r]\in {\dv R}^{N}$ with $r\in [0,1]$.\\
Finally, if the considered linear combination is not conic then $\calS$  has to contain at least three
different sets and the existence of such a set system forces $|N|\geq 3$.
\end{lem}

It follows from (the last claim in) Lemma~\ref{lem.semi-comb} that every purely semi-balanced system contains at least three sets and, thus, there is no such a set system on a two-element set.

\begin{proof}
The case of a conic combination is easy: choose $L\in\calS$ with $\lambda_{L}>0$, $i\in L$ and write $0<\lambda_{L}\leq \sum_{S\in\calS} \lambda_{S}\cdot\chi_{S}(i)=r=\sum_{S\in\calS:\,i\in S} \lambda_{S}\leq \sum_{S\in\calS} \lambda_{S}$. Thus, in the rest of the proof we assume that the combination is not conic.
\smallskip

Let $T\in\calS$ be the set with $\lambda_{T}<0$; note that $\emptyset\subset T\subset N$.
The choice of $j\in N\setminus T$ gives
$r=\sum_{S\in\calS} \lambda_{S}\cdot\chi_{S}(j)= \sum_{S\in\calS:\,j\in S} \lambda_{S}\geq 0$
because $\lambda_{S}\geq 0$ whenever $S\in \calS\setminus\{T\}$. In case of the choice
$i\in T$ one has $i\not\in S\in\calS ~\Rightarrow~ S\in \calS\setminus\{T\}~\Rightarrow~ \lambda_{S}\geq 0$,
which gives $\sum_{S\in\calS} \lambda_{S}\geq \sum_{S\in\calS:\,i\in S} \lambda_{S}=\sum_{S\in\calS} \lambda_{S}\cdot\chi_{S}(i)=r$. For the verification of the claim $\sum_{S\in\calS} \lambda_{S}>0$
assume without loss of generality that $\lambda_{S}\neq 0$ for any $S\in\calS$, for otherwise one can replace $\calS$ by $\calS^{\prime}:=\{S\in\calS\,:\
\lambda_{S}\neq 0\}$. We distinguish two cases:
\begin{itemize}
\item In case [\,$\exists\, L\in\calS\setminus\{T\}$ with $T\setminus L\neq\emptyset$\,] we choose
$k\in T\setminus L$. Then $\lambda_{S}\geq 0$ whenever $k\not\in S\in\calS$ and we write $\sum_{S\in\calS} \lambda_{S}\geq \lambda_{L}+ \sum_{S\in\calS:\,k\in S} \lambda_{S}
> \sum_{S\in\calS:\,k\in S} \lambda_{S}= r\geq 0$.
\item In case [\,$\forall\, L\in\calS\setminus\{T\}$ one has $T\subseteq L$\,] we first observe
$\bigcup\calS\setminus T\neq\emptyset$. Indeed, by contradiction: if $\bigcup\calS\subseteq T$
then $\forall\, L\in\calS\setminus \{T\}$ one has $L\subseteq T\subseteq L$, which means
$\calS=\{T\}$ contradicting $\sum_{S\in\calS} \lambda_{S}\geq 0$.
Thus, one can choose $k\in\bigcup\calS\setminus T$, fix $K\in \calS$ with $k\in K$ and write
$\sum_{S\in\calS} \lambda_{S}\geq r=\sum_{S\in\calS:\,k\in S} \lambda_{S}\geq \lambda_{K}>0$.
\end{itemize}
Thus, in both cases we have $\sum_{S\in\calS} \lambda_{S}> 0$.
\smallskip

In particular, given a non-zero semi-conic combination $\sum_{S\in\calS} \lambda_{S}\cdot\chi_{S}$ yielding a constant vector in ${\dv R}^{N}$
we put $\ell:=\sum_{S\in\calS} \lambda_{S}>0$ and observe that $\sum_{S\in\calS} (\ell^{-1}\cdot\lambda_{S})\cdot\chi_{S}$ is the required affine semi-conic combination.\smallskip

To verify the last claim take $T\in\calS$ with $\lambda_{T}<0$ and
assume for a contradiction $|\calS|\leq 2$. If $|\calS|=1$
take $j\in T$ and get a contradictory conclusion $0\leq r=\sum_{S\in\calS:\,j\in S} \lambda_{S}=\lambda_{T}<0$.
In case $|\calS|=2$ one has $\calS =\{T,L\}$ and contingent choice of $j\in T\setminus L$ leads to an analogous contradiction.
Hence, one must have $T\subset L$ and the choice of $j\in T$ and $k\in L\setminus T$
leads to another contradiction $r=\sum_{S\in\calS:\,k\in S} \lambda_{S}=\lambda_{L}>\lambda_{L}+\lambda_{T}=
\sum_{S\in\calS:\,j\in S} \lambda_{S}=r$.
Thus, $|\calS|\geq 3$; since $|N|=2 ~\Rightarrow~ |{\caP}\setminus \{\emptyset,N\}|=2$,
the existence of $\calS$ forces $|N|\geq 3$.
\end{proof}

We now provide equivalent definitions of min-semi-balanced systems.
Slightly longer proof of the next lemma is shifted to Appendix~\ref{app.lem-2}.

\begin{lem}\label{lem.min-semi-bal}\em
Given $|N|\geq 2$, let $\emptyset\neq\calS\subseteq {\caP}\setminus \{\emptyset,N\}$ be a non-trivial set system on $N$. Then the following conditions on $\calS$ are equivalent:
\begin{description}
\item[(a)] $\calS$ is a minimal set system such that there is a constant vector in ${\dv R}^{N}$ which can be
written as a non-zero semi-conic combination of vectors $\{ \chi_{S}\,:\ S\in\calS\}$,
\item[(b)] $\calS$ is a minimal semi-balanced set system on $N$,
\item[(c)] $\calS$ is semi-balanced on $N$, the vectors $\{ \chi_{S}\,:\ S\in\calS\}$ are affinely independent
and in case $\bigcup\calS=N$ even linearly independent,
\item[(d)] there is only one affine combination of vectors $\{ \chi_{S}\,:\ S\in\calS\}$ yielding
a constant vector in ${\dv R}^{N}$ and this unique combination is semi-conic and has all coefficients non-zero,
\item[(e)] there exists unique affine semi-conic combination of vectors $\{ \chi_{S}\,:\ S\in\calS\}$
which is a constant vector in ${\dv R}^{N}$ and this unique combination has all coefficients non-zero.
\end{description}
\end{lem}

Given (purely) min-semi-balanced system $\calS$ on $N$ with an exceptional set $T$,
by Lemma~\ref{lem.semi-comb} one can consider an affine combination $\sum_{S\in\calS} \lambda_{S}\cdot\chi_{S}$ yielding a constant vector where $\lambda_{T}<0$ and $\lambda_{S}\geq 0$ for $S\in\calS\setminus\{T\}$.
Then, by Lemma~\ref{lem.min-semi-bal}(d), such a combination is unique. This implies
that the exceptional set $T$ within $\calS$ is uniquely determined by $\calS$.

The following Example~\ref{ex.non-c} illustrates that the requirement concerning the case $\bigcup\calS=N$ in the condition {\bf (c)} of Lemma~\ref{lem.min-semi-bal} cannot be removed.

\begin{example}\label{ex.non-c}\rm
There exists a semi-balanced set system $\calS$ on $N$ with $\bigcup\calS=N$ where
vectors $\{ \chi_{S}\,:\ S\in\calS\}$ are affinely independent but not linearly independent.
Put $N:=\{a,b,c,d\}$ and $\calS:=\{\, a,b,ab,abc, abd\,\}$. Then one has
$$
1\cdot\chi_{a} + 1\cdot\chi_{b} + (-2)\cdot\chi_{ab} + 1\cdot\chi_{abc} + 1\cdot\chi_{abd} =\chi_{N}\,,
$$
which implies that $\calS$ is semi-balanced. The equality $\0=1\cdot\chi_{a} + 1\cdot\chi_{b}+ (-1)\cdot\chi_{ab}$
implies that $\{ \chi_{S}\,:\ S\in\calS\}$ are linearly dependent. On the other hand, the condition
$$
\alpha\cdot\chi_{a} + \beta\cdot\chi_{b} + \gamma\cdot\chi_{ab} + \delta\cdot\chi_{abc} + \varepsilon\cdot\chi_{abd} =\0
$$
together with $\alpha+ \beta+ \gamma+ \delta+ \varepsilon=0$ implies $\alpha=\ldots =\varepsilon=0$, which means
that $\{ \chi_{S}\,:\ S\in\calS\}$ are affine independent. Indeed, take $d$ first to derive
$\varepsilon=0$, then $c$ to get $\delta=0$, $a$ to obtain $\alpha=-\gamma$ and $b$ to obtain $\beta=-\gamma$;
thus, $0=\alpha+ \beta+ \gamma+ \delta+ \varepsilon=-\gamma$ gives the conclusion.
Hence, by Lemma~\ref{lem.min-semi-bal}, the system $\calS$ is {\sf not\/} min-semi-balanced on $N$; two of its proper subsystems that are semi-balanced are $\{\,a,b,ab\,\}$ and  $\{\,ab,abc, abd\,\}$.
\end{example}

The basic requirement in the condition {\bf (d)} of Lemma~\ref{lem.min-semi-bal} is geometric: 
the affine subspace $\{\,\sum_{S\in\calS} \lambda_{S}\cdot\chi_{S}\in {\dv R}^{N}\,:\ \sum_{S\in\calS} \lambda_{S}=1\}$
intersects the line $\{[r,\ldots,r]\in {\dv R}^{N}\,:\ r\in {\dv R}\}$ of constant vectors in exactly one vector.
The following Example~\ref{ex.non-d} shows that both additional requirements on this unique affine combination are necessary.

\begin{example}\label{ex.non-d}\rm
There exists a non-trivial system $\calD$ on $N$ such that only one affine combination of $\{ \chi_{S}\,:\ S\in\calD\}$ yields a constant vector in ${\dv R}^{N}$ and this unique combination is conic but has not all coefficients non-zero.
Put $N:=\{a,b,c\}$ and $\calD:=\{\,a,b,bc\,\}$. Then
$$
\frac{1}{2}\cdot\chi_{a} + 0\cdot\chi_{b} + \frac{1}{2}\cdot\chi_{bc}=\frac{1}{2}\cdot\chi_{N}\,,
$$
is the above-mentioned unique affine combination. Of course, this particular set system $\calD$ is not semi-balanced but its subsystem $\calD^{\prime}:=\{\,a,bc\,\}$ is even balanced.

There is also a non-trivial set system $\calC$ on $N$ with a unique affine combination of $\{ \chi_{S}\,:\ S\in\calC\}$
yielding a constant vector in ${\dv R}^{N}$, where all the coefficients in this combination are non-zero
but the combination is not semi-conic. Put $N:=\{a,b,c,d,e\}$ and consider $\calC:=\{\,ab,ac,ad,abc,abce\,\}$. Then
$$
(-1)\cdot\chi_{ab} + (-1)\cdot\chi_{ac} + 1\cdot\chi_{ad} + 1\cdot\chi_{abc} + 1\cdot\chi_{abce}=\chi_{N}
$$
is that unique combination. By Lemma~\ref{lem.semi-comb}, the system $\calC$ is not semi-balanced on $N$.
\end{example}

A consequence of Lemma~\ref{lem.min-semi-bal} is the observation that the
concept of a min-semi-balanced system generalizes the one of a min-balanced system.

\begin{cor}\label{cor.min-balanced}\rm
Given $|N|\geq 2$, a {\sf balanced\/} system $\calB$ on $N$ is minimal within the class of balanced systems
on $N$ iff it is minimal within the class of semi-balanced systems on $N$. In other words,
$\calB$ is min-balanced (on $N$) iff it is balanced and min-semi-balanced (on $N$).
\end{cor}

\begin{proof}
The fact that $\calB$ is balanced on $N$ implies $\bigcup\calB=N$. By Lemma~\ref{lem.min-semi-bal}{\bf (c)},
its minimality within semi-balanced systems means that $\{ \chi_{S}\,:\ S\in\calB\}$ are linearly independent,
while, by \cite[Lemma\,2.1]{KS19}, this condition characterizes its minimality within balanced systems.
\end{proof}

Lemma~\ref{lem.min-semi-bal}{\bf (c)} also sets a limit on the number of sets in a min-semi-balanced system.

\begin{cor}\label{cor.cardinality-semi-balanced}\rm
Given $|N|\geq 2$ and a {\sf min-semi-balanced\/} system $\calS$ on $N$
one has $|\calS|\leq |N|$.
\end{cor}

\begin{proof}
If $\bigcup\calS\subset N$ then affine independence of $\{\chi_{S}\,:\, S\in\calS\}$ gives $|\calS|\leq |\bigcup\calS|+1\leq|N|$ because affinely independent set in ${\dv R}^{M}$ has at most $|M|+1$ elements. In case $\bigcup\calS =N$  the linear independence of the respective set of vectors implies directly $|\calS|\leq |N|$.
\end{proof}

\subsection{Inequalities assigned to semi-balanced systems}\label{ssec.inequalities}
Given a non-trivial set system $\calS$ on $N$, any non-zero semi-conic combination
$\sum_{S\in\calS} \lambda_{S}\cdot\chi_{S}$ yielding a constant vector $[r,\ldots,r]\in {\dv R}^{N}$
gives an inequality
$$
r\cdot\m (N) ~\geq~ \sum_{S\in\calS} \lambda_{S}\cdot\m(S)\qquad \mbox{for $\m\in{\dv R}^{\caP}$},
$$
which appears to be valid for all exact games $\m$.
To ensure one-to-one correspondence between the inequalities and semi-conic combinations we limit ourselves
to affine combinations, which is possible owing to Lemma~\ref{lem.semi-comb}.
The formal definition is as follows.

\begin{defin}[vectors of coefficients, inequalities induced by set systems]~\label{def.theta}\rm \\
Let $\calS$ be a non-trivial set system on $N$, $|N|\geq 2$.
Any {\sf affine\/} semi-conic combination $\sum_{S\in\calS} \lambda_{S}\cdot\chi_{S}$ yielding a constant vector
$[r,\ldots,r]\in {\dv R}^{N}$ is assigned a vector $\theta\in {\dv R}^{\caP}$:
\begin{equation}
\theta (S)= -\lambda_{S}~~\mbox{for $S\in\calS$},\quad
\theta(N)=r,\quad \theta (\emptyset)= 1-r, \quad
\theta (L)=0~~\mbox{for other $L\subseteq N$.}
\label{eq.theta}
\end{equation}
The vector $\theta$ is then interpreted as the coefficient vector in an inequality
\begin{equation}
0\leq \langle \theta,\m\rangle ~:=~ \sum_{S\subseteq N}\, \theta(S)\cdot\m(S)\qquad
\mbox{for vectors $\m\in {\dv R}^{\caP}$ with $\m(\emptyset)=0$.}
\label{eq.theta-ineq}
\end{equation}
The symbol $\Theta_{\calS}$ will denote the set of vectors $\theta$ for all such affine semi-conic combinations.
Provided $|\Theta_{\calS}|=1$, the only vector in $\Theta_{\calS}$ will be denoted by $\theta_{\calS}$.
\end{defin}

Note that every coefficient vector $\theta$ is given by a suitable affine combination
and that the values $\theta(S)$ for $S\in\calS$ correspond to the coefficients in the combination;
however, the remaining contingent non-zero values $\theta(N)$ and $\theta(\emptyset)$ are determined by them.

Here are some elementary observations on the set of coefficient vectors.

\begin{cor}\label{cor.theta}~\rm
Given $|N|\geq 2$ and a non-trivial set system $\calS$ on $N$, one has $\Theta_{\calS}\neq\emptyset$
iff $\calS$ contains a semi-balanced system. The inclusion $\calT\subseteq\calS$ of two non-trivial systems on $N$ implies $\Theta_{\calT}\subseteq\Theta_{\calS}$. The set $\Theta_{\calS}$ is the union of sets
$\Theta_{\calT}$ for semi-balanced systems $\calT$ with $\calT\subseteq\calS$.
One has $|\Theta_{\calS}|=1$ iff $\calS$ contains just one semi-balanced system on $N$.
In particular, every min-semi-balanced system $\calS$ on $N$ satisfies $|\Theta_{\calS}|=1$.
\end{cor}

These facts mean that the inequalities \eqref{eq.theta-ineq} from Definition~\ref{def.theta} are just those
that are assigned to semi-balanced systems. The last claim says that only a non-min-semi-balanced system $\calS$ may have non-singleton $\Theta_{\calS}$ and, thus, be assigned several inequalities.

On the other hand, the substantial inequalities appear to be those assigned
to min-semi-balanced systems; see later Corollary~\ref{cor.exa-char}(iii). Thus, the inequalities assigned to other non-trivial systems are superfluous for the description of the exact cone.
Nonetheless, in order to follow the analogy with former results by Lohmann {\em et al.\/} \cite{LBH12},
we have also assigned the inequalities to non-minimal semi-balanced set systems;
see later Corollary~\ref{cor.exa-char}(ii).

\begin{proof}
The first claim follows from Lemma~\ref{lem.semi-comb} and Definition~\ref{def.semi-bal}, further two ones are direct
consequences of Definition~\ref{def.theta}. As concerns the fourth claim, the necessity of the uniqueness of a
semi-balanced subsystem can be shown by contradiction with the help of Lemma~\ref{lem.semi-comb}.
For its sufficiency realize that the unique semi-balanced system $\calT$ on $N$ with $\calT\subseteq\calS$
is necessarily minimal. Given $\theta\in\Theta_{\calS}$, the set system
$\{S\subset N:\ S\neq\emptyset ~\&~ \theta(S)\neq 0\}\subseteq\calS$ is semi-balanced on $N$, and, thus, it has to coincide with $\calT$.
The condition {\bf (e)} in Lemma~\ref{lem.min-semi-bal} (for $\calT$) implies the uniqueness of $\theta\in\Theta_{\calS}$; hence, $|\Theta_{\calS}|\leq 1$. This implies the last claim.
\end{proof}

The next example shows that the set $\Theta_{\calS}$ need not be convex.

\begin{example}\label{ex.non-convex}\rm
Put $N:=\{a,b,c,d\}$ and $\calB:=\{\,a,b,c,d,ab,cd\,\}$.
The conic combination
$$
\frac{1}{2}\cdot\chi_{a}+\frac{1}{2}\cdot\chi_{b}+\frac{1}{2}\cdot\chi_{c}+\frac{1}{2}\cdot\chi_{d} + \frac{1}{2}\cdot\chi_{ab}+ \frac{1}{2}\cdot\chi_{cd}
= \chi_{N}
$$
means that $\calB$ is balanced on $N$. A semi-conic combination
$1\cdot\chi_{a}+1\cdot\chi_{b}+(-1)\cdot\chi_{ab}=\0$ leads to a coefficient vector $\theta^{1}\in\Theta_{\calB}$:
$$
\theta^{1}(a)=-1, ~~\theta^{1}(b)=-1,~~ \theta^{1}(ab)=+1, ~~\theta^{1}(\emptyset)=+1,\quad \theta^{1}(L)=0
~~\mbox{otherwise}.
$$
Analogously, $1\cdot\chi_{c}+1\cdot\chi_{d}+(-1)\cdot\chi_{cd}=\0$ gives rise to the vector $\theta^{2}\in\Theta_{\calB}$:
$$
\theta^{2}(c)=-1, ~~\theta^{2}(d)=-1,~~ \theta^{2}(cd)=+1, ~~\theta^{2}(\emptyset)=+1,\quad \theta^{1}(L)=0
~~\mbox{otherwise}.
$$
Their convex combination $\frac{1}{2}\cdot\theta^{1}+\frac{1}{2}\cdot\theta^{2}$, however, does not belong to $\Theta_{\calB}$
because the corresponding affine combination $\frac{1}{2}(\chi_{a}+\chi_{b}+\chi_{c}+\chi_{d}-\chi_{ab}-\chi_{cd})=\0$ is
not semi-conic.
\end{example}

On the other hand, $\Theta_{\calS}$ is always the union of finitely many closed convex sets.

\begin{cor}\label{cor.rational}~\rm
Given a semi-balanced system $\calS$ on $N$, $|N|\geq 2$, the set $\Theta_{\calS}$ is the union of finitely many
rational polyhedrons. In particular, if $\calS$ is min-semi-balanced then the unique vector $\theta_{\calS}$ in $\Theta_{\calS}$ has rational components: $\theta_{\calS}\in {\dv Q}^{\caP}$.
\end{cor}

\begin{proof}
In fact, $\Theta_{\calS}$ is the union over $T\in\calS$ of sets $\Theta_{\calS:T}$ consisting of
$\theta\in {\dv R}^{\caP}$ such that
\begin{itemize}
\item $\theta(N)+\theta(\emptyset)=1=\sum_{S\in\calS}-\theta(S)$,\qquad  $\sum_{S\in\calS} -\theta(S)\cdot\chi_{S}=\theta(N)\cdot\chi_{N}$,
\item $\forall\, S\in\calS\setminus\{T\}\quad \theta(S)\leq 0$,
\qquad $\forall\, L\in \caP\setminus (\calS\cup\{\emptyset,N\})\quad \theta (L)=0$.
\end{itemize}
Indeed, the above conditions defining $\Theta_{\calS:T}$ are the rewriting of \eqref{eq.theta} and of the requirements on
the respective linear combination $\sum_{S\in\calS} \lambda_{S}\cdot\chi_{S}$, where $\lambda_{T}$ is allowed to be negative.
These constraints have clearly rational coefficients.
In case of a minimal $\calS$ one has $|\Theta_{\calS}|=1$ by Corollary~\ref{cor.theta}.
Thus, $\{\theta_{\calS}\}=\Theta_{\calS:T}$ for some $T\in\calS$ then.
Because $\Theta_{\calS:T}$ is specified by rational constraints it is a rational polyhedron.
Another well-known fact from polyhedral geometry that every vertex of a rational polyhedron
has rational components; see \cite[Statement\,3]{Studeny93} for example.
This gives $\theta_{\calS}\in {\dv Q}^{\caP}$.
\end{proof}

Note that later Lemma~\ref{lem.theta} implies that every set $\Theta_{\calS:T}$ from the above proof is bounded.

\section{Characterization of exact games}\label{sec.char-exact}
The inequalities of the form \eqref{eq.theta-ineq} from Definition~\ref{def.theta} allow one to delimit the cone of exact games. To this end we introduce the following convex sets.

\begin{defin}[auxiliary cones and polyhedrons]~\label{def.aux-cone} \rm \\
Given $|N|\geq 2$ and $\emptyset\neq D\subseteq N$ we put
\begin{eqnarray*}
\lefteqn{\hspace*{-0.6cm}\Theta^{N}_{D} ~~:=~\, \{\,\, \theta\in {\dv R}^{\caP}\, :\
 \theta (S)\leq 0 \qquad~~\mbox{for any $S\subseteq N$ such that $S\not\in\{\emptyset, D, N\}$,} }\\
 &&
~~\qquad\qquad ~~~~\sum_{L\subseteq N} \theta (L)=0 ~~~\mbox{and}~ \sum_{L\subseteq N:\, i\in L} \theta (L)=0 ~~\mbox{for any $i\in N$}\,\},\\
\tilde{\Theta}^{N}_{D} &:=& \Theta^{N}_{D}\,\cap\, \{\,\theta\in {\dv R}^{\caP}\,:\ \theta(N)+\theta(\emptyset)=1\,\},\\
\Delta~ &:=& \conv (\bigcup_{D:\,\emptyset\neq D\subseteq N} \tilde{\Theta}^{N}_{D})\,\,.
\end{eqnarray*}
\end{defin}
\medskip

Observe that $\tilde{\Theta}^{N}_{D} = \Theta^{N}_{D}\,\cap\, \{\theta\in {\dv R}^{\caP}\,:\
\sum_{L:\,\emptyset\neq L\subset N} \theta(L)=-1\,\}$,  which re-writing allows one to ignore
the component for $\emptyset$ and interpret these convex sets as subsets of ${\dv R}^{\caP\setminus\{\emptyset\}}$.
It follows directly from the definition
that $\tilde{\Theta}^{N}_{N}\subseteq \tilde{\Theta}^{N}_{D}$ for any $\emptyset\neq D\subset N$; thus,
one can, alternatively, consider the union over $\emptyset\neq D\subset N$ is the definition of $\Delta$.
The proof of the next lemma, based on some facts from \cite[\S\,5.1]{KS19},
is moved to Appendix~\ref{app.lem-theta}.

\begin{lem}\label{lem.theta}~\rm
Given $|N|\geq 2$, every set $\tilde{\Theta}^{N}_{D}$, where $\emptyset\neq D\subseteq N$, is a bounded polyhedron.
Every vector $\theta\in\Theta^{N}_{D}$ satisfies both $\theta(N)\geq 0$ and $\theta(\emptyset)\geq 0$ and
every {\sf non-zero\/} vector $\theta\in\Theta^{N}_{D}$ satisfies $\theta(N)+\theta(\emptyset)>0$.
Given $\m\in {\dv R}^{\caP}$ with $\m(\emptyset)=0$, one has
\begin{equation}
\m\in {\cal E}(N) ~~\Leftrightarrow~~ [\,\forall\, \theta\in \bigcup_{\emptyset\neq D\subseteq N}\tilde{\Theta}^{N}_{D}\qquad \langle \theta,\m\rangle\geq 0\,]\,.
\label{eq.exact-char}
\end{equation}
\end{lem}
\smallskip

Observe that the set $\Theta_{\calS:T}$ from the proof of Corollary~\ref{cor.rational}, where $\calS:=\caP\setminus\{\emptyset,N\}$ and $T\in\calS$, coincides with $\tilde{\Theta}^{N}_{T}$.
Indeed, realize that the conditions $\sum_{L\subseteq N:\, i\in L} \theta (L)=0$
for $i\in N$ mean $\sum_{S\in\calS} -\theta(S)\cdot\chi_{S}=\theta(N)\cdot\chi_{N}$.
Hence, $\Theta_{\calS:T}$ is bounded by the first claim in Lemma~\ref{lem.theta}.
The first claim in Lemma~\ref{lem.theta} also allows one to observe that $\Delta$ is a bounded polyhedron as well;
this follows from basic facts in polyhedral geometry recalled in Section~\ref{ssec.polyhedrals}.

The second claim  in Lemma~\ref{lem.theta} means that each set $\Theta^{N}_{D}$, $\emptyset\neq D\subseteq N$, is a pointed polyhedral cone. An equivalent formulation is that every non-zero vector $\theta\in\Theta^{N}_{D}$ satisfies
$\sum_{L:\,\emptyset\neq L\subset N} \theta(L)<0$, which is relevant if $\Theta^{N}_{D}$ is interpreted as a subset of ${\dv R}^{\caP\setminus\{\emptyset\}}$.
\bigskip

Further auxiliary observation says that the vertices of the (bounded) polyhedrons from
Definition~\ref{def.aux-cone} correspond to (certain) min-semi-balanced set systems on $N$. Thus,
together with Lemma~\ref{lem.theta}, it puts in relation exact games and semi-balanced systems.
Its proof is also shifted to Appendix~\ref{app.lem-extr-theta}.

\begin{lem}\label{lem.extr-theta}~\rm
Given $|N|\geq 2$ and $\emptyset\neq D\subseteq N$, every vertex of $\tilde{\Theta}^{N}_{D}$ has either the form
$\theta_{\calB}$, where $\calB$ is a min-balanced set system on $N$, or the form $\theta_{\calS}$, where $\calS$ is a min-semi-balanced system on $N$ having $D$ as the exceptional set.\\
Conversely, in case $\emptyset\neq D\subset N$, every vector $\theta_{\calS}$, where
$\calS$ is a min-semi-balanced system on $N$ having $D$ as the exceptional set, is a vertex of
$\tilde{\Theta}^{N}_{D}$: $\theta_{\calS}\in \ext (\tilde{\Theta}^{N}_{D})$.
\end{lem}
\smallskip

Note in this context that one can show, using the same arguments as in the proof of Lemma~\ref{lem.extr-theta},
that $\ext(\tilde{\Theta}^{N}_{N})$ consists just of the vectors $\theta_{\calB}$ where $\calB$ is a min-balanced set system on $N$;
nonetheless, this observation is not necessarily needed to derive our results. On the other hand,
the delimitation of $\ext (\tilde{\Theta}^{N}_{D})$ for $D\subset N$ in the first claim of Lemma~\ref{lem.extr-theta}
is not tight. Analogous arguments can be used to show that $\theta_{\calB}\in \ext (\tilde{\Theta}^{N}_{D})$
for every min-balanced system $\calB$ on $N$ with $D\in\calB$.
For a min-balanced system $\calB$ on $N$ such that $D\not\in\calB$, however,
the vector $\theta_{\calB}$ may or may not be a vertex of $\tilde{\Theta}^{N}_{D}$ as the next example shows.

\begin{example}\label{ex.not-all}\rm
Put $N:=\{a,b,c\}$ and $D:=\{ a,b\}$; the set system $\calB:=\{\,a,b,c\,\}$ is then min-balanced on $N$.
The corresponding vector $\theta_{\calB}\in {\dv R}^{\caP}$ is given by
$$
\theta_{\calB}(N)=+\frac{1}{3}, ~~\theta_{\calB}(a)=\theta_{\calB}(b)=\theta_{\calB}(c)=-\frac{1}{3}, ~~\theta_{\calB}(\emptyset)=+\frac{2}{3},
\quad \theta_{\calB}(L)=0
~~\mbox{otherwise},
$$
evidently belongs to $\tilde{\Theta}^{N}_{D}$. Consider another min-balanced system $\calC:=\{\,c,ab\,\}$ on $N$ with
$$
\theta_{\calC}(N)=\theta_{\calC}(\emptyset)=+\frac{1}{2}, ~~\theta_{\calC}(c)=\theta_{\calC}(ab)=-\frac{1}{2},
\quad \theta_{\calC}(L)=0 ~~\mbox{otherwise},
$$
and a min-semi-balanced system $\calD:=\{\,a,b,ab\,\}$ on $N$ with
$$
\theta_{\calD}(ab)=\theta_{\calD}(\emptyset)=+1, ~~\theta_{\calD}(a)=\theta_{\calD}(b)=-1,
\quad \theta_{\calC}(L)=0 ~~\mbox{otherwise}\,.
$$
These two vectors both belong to $\tilde{\Theta}^{N}_{D}$ and
one has $\theta_{\calB}=\frac{2}{3}\cdot\theta_{\calC}+\frac{1}{3}\cdot\theta_{\calD}$.
In particular, $\theta_{\calB}$ is not a vertex of $\tilde{\Theta}^{N}_{D}$. On the other hand,
the min-balanced system $\calB^{\prime}:=\{\,a,bc\,\}$ on $N$ with
$$
\theta_{\calB^{\prime}}(N)=\theta_{\calB^{\prime}}(\emptyset)=+\frac{1}{2}, ~~\theta_{\calB^{\prime}}(a)=\theta_{\calB^{\prime}}(bc)=-\frac{1}{2},
\quad \theta_{\calB^{\prime}}(L)=0 ~~\mbox{otherwise},
$$
also complies with $D\not\in\calB^{\prime}$ and it makes no problem to show that $\theta_{\calB^{\prime}}\in \ext (\tilde{\Theta}^{N}_{D})$.
\end{example}

Lemmas~\ref{lem.theta} and \ref{lem.extr-theta} allow one to characterize exact games in terms
of semi-balanced systems.

\begin{cor}\label{cor.exa-char}~\rm
Given $|N|\geq 2$, consider a function $\m\in {\dv R}^{\caP}$ such that $\m(\emptyset)=0$. Then the following
conditions on the game $\m$ are equivalent:
\begin{description}
\item[(i)] ~\,\,$\m$ is an exact game over $N$, that is, $\m\in \calE(N)$,
\item[(ii)] ~for every semi-balanced set $\calS$ system on $N$ one has $\langle \theta,\m\rangle\geq 0$
for all $\theta\in\Theta_{\calS}$,
\item[(iii)] for every min-semi-balanced set $\calS$ system on $N$ one has
$\langle \theta_{\calS},\m\rangle\geq 0$.
\end{description}
\end{cor}

Note that one is entitled to write $\theta_{\calS}$ in {\bf (iii)} since, by Corollary~\ref{cor.theta}, $|\Theta_{\calS}|=1$ then.

\begin{proof}
The implication {\bf (i)}\,$\Rightarrow${\bf (ii)} follows from \eqref{eq.exact-char}
in Lemma~\ref{lem.theta}, one only needs to realize that, given a semi-balanced system $\calS$ on $N$, one has $\Theta_{\calS}\subseteq\tilde{\Theta}^{N}_{D}$ for some $\emptyset\neq D\subseteq N$. Indeed,
any $\theta\in\Theta_{\calS}$ is defined in \eqref{eq.theta} from an affine (semi-conic) combination $\sum_{S\in\calS} \lambda_{S}\cdot\chi_{S}$ yielding a constant vector in ${\dv R}^{N}$, which
gives both $\theta(N)+\theta(\emptyset)=-\sum_{L:\emptyset\neq L\subset N} \theta(L)=+1$ and
$\sum_{L\subseteq N:i\in L} \theta(L)=0$ for any $i\in N$.
Thus, if it is a conic combination then $\theta\in\tilde{\Theta}^{N}_{N}$, otherwise one has
$\theta\in\tilde{\Theta}^{N}_{T}$ for the only (exceptional) $T\in\calS$ with $\lambda_{T}<0$.
\smallskip

The implication {\bf (ii)}\,$\Rightarrow${\bf (iii)} is immediate.
\smallskip

To verify {\bf (iii)}\,$\Rightarrow${\bf (i)} we apply the first claim in Lemma~\ref{lem.extr-theta}
to observe that $\langle \theta,\m\rangle\geq 0$ for any $\theta\in \ext (\tilde{\Theta}^{N}_{D})$
and $\emptyset\neq D\subseteq N$. Hence, the same holds for any $\theta\in \tilde{\Theta}^{N}_{D}$
and arbitrary $\emptyset\neq D\subseteq N$. In particular, one can use \eqref{eq.exact-char} in Lemma~\ref{lem.theta}
to derive {\bf (i)}.
\end{proof}
\smallskip

Let us remark that the observations from Corollary~\ref{cor.exa-char} are analogous to former results
by Lohmann {\em et al.} \cite{LBH12}, specifically to Theorem~3.4 in \cite[\S\,3]{LBH12}
and Theorem~5.1 in \cite[\S\,5]{LBH12}. The proviso is that Lohmann {\em et al.\/} used
a different formal way to associate set systems with inequalities -- see later Remark~\ref{rem.lohmnamm}
for the explanation. We believe that our approach and presentation offers an elegant geometric interpretation
and simpler arguments. Moreover, it can be extended to get the following characterization of facet-defining inequalities.

\begin{cor}\label{cor.facet-char}~\rm
Given $|N|\geq 2$, the inequality \eqref{eq.theta-ineq}, that is, $0\leq \langle \theta,\m\rangle$ for $\m\in {\dv R}^{\caP}$, 
with a coefficient vector $\theta\in {\dv R}^{\caP}$, where $\sum_{S\subseteq N} \theta (S)=0$ and
$\sum_{L:\,\emptyset\neq L\subset N} \theta(L)=-1$, is facet-defining for the
cone of exact games $\calE (N)$ iff $\theta\in \ext (\Delta)$.
\end{cor}

Note that the above requirements $\sum_{S\subseteq N} \theta (S)=0$ and
$\sum_{L:\,\emptyset\neq L\subset N} \theta(L)=-1$ on a coefficient vector $\theta$
are solely technical constraints which can, without loss of generality,
be assumed to hold for any facet-defining inequality for the exact cone.
The former requirement is related to the fact that $\m(\emptyset)=0$ for any game $\m\in {\cal G}(N)$
and ${\cal E}(N)$ is, in fact, a full-dimensional cone in ${\dv R}^{\caP\setminus\{\emptyset\}}$: it is a
convention on the value $\theta(\emptyset)$. The latter requirement is related to the fact that
facet-defining inequalities are determined uniquely up to a positive multiple: it is
a particular convention about the choice of the multiplicative factor.

The proof is based on a geometric consideration concerning the duality of polyhedral cones.
It is more convenient technically to imagine both the cone of exact games and its dual cone
within the space ${\dv R}^{\caP\setminus\{\emptyset\}}$
because then the dual cone becomes pointed.

\begin{proof}
Consider the space ${\dv R}^{\caP\setminus\{\emptyset\}}$ and interpret the polyhedrons
$\tilde{\Theta}^{N}_{D}$ from Definition~\ref{def.aux-cone} as its subsets.
Introduce a polyhedral cone $\bar{\Delta}\subseteq {\dv R}^{\caP\setminus\{\emptyset\}}$ as the conic hull
of $\Delta$. Since every $\theta\in\Delta$ satisfies $\sum_{\emptyset\neq L\subset N} \theta(L)=-1$
every non-zero $\theta\in \bar{\Delta}$ satisfies $\sum_{\emptyset\neq L\subset N} \theta(L)<0$
meaning that $\bar{\Delta}$ is pointed and its extreme rays are generated by vertices
$\theta\in\ext (\Delta)$.
The first claim in Lemma~\ref{lem.theta} allows one to observe that $\Delta$ is a bounded polyhedron,
while the third claim in Lemma~\ref{lem.theta} means that the cone $\calE(N)$ is dual to $\bar{\Delta}$:
$\calE(N)=(\bar{\Delta})^{\dual}$.

Hence, by the basic facts from polyhedral geometry recalled in Section~\ref{ssec.polyhedrals},
the cones $\bar{\Delta}$ and $\calE(N)$ are mutually dual polyhedral cones and
the lattices of their non-empty faces are anti-isomorphic. In particular, extreme rays of $\bar{\Delta}$
correspond to facets of $\calE(N)$. Moreover, the lattice of non-empty faces of $\bar{\Delta}$ is isomorphic to the face-lattice of $\Delta$, which implies that the extreme rays of $\bar{\Delta}$ are just those
rays that are generated by vertices of $\Delta$. In other words, a coefficient vector $\theta$ normalized by $\sum_{\emptyset\neq L\subset N} \theta(L)=-1$ yields a facet-defining inequality for $\calE(N)$
iff $\theta$ is a vertex of $\Delta$.
\end{proof}

\begin{remark}\label{rem.lohmnamm}\rm
Lohmann {\em et al.\/} introduced in \cite{LBH12} the concept of
an {\em exact balanced collection\/} of sets with the intention to use
such set systems to generate linear inequalities specifying the cone of exact games.
These collections of sets often coincide with our semi-balanced set systems, but there is one (substantial)
technical difference in their approach. It concerns of the inequalities $\langle \theta,\m\rangle\geq 0$ for
$\m\in \calE(N)$ with $\theta (N)=0$. They ascribe such inequalities to the so-called ``minimal sub-balanced"
collections which are certain set systems always involving the grand coalition $N$.
To give an example of the difference consider the vector $\theta\in {\dv R}^{\caP}$ where
$N=\{a,b,c\}$ given by
$$
\theta(\emptyset) =+1,\quad \theta (ab)=+1,\qquad \theta(a)=-1,\quad \theta(b)=-1, \qquad \theta(L)=0~~\mbox{for other $L\subseteq N$}\,.
$$
One has $\theta\in\Theta^{N}_{D}$ for $D=ab$ and $\langle \theta,\m\rangle\geq 0$ is facet-defining
for $\m\in \calE(N)$. Our approach is to associate this inequality with a set system
$\{\, a,b,ab\,\}$ while Lohmann {\em et al.\/} \cite{LBH12} ascribe that inequality to the set system
$\{\, a,b,ab, N\,\}$. We have two arguments why their approach is not appropriate for our purpose:
\begin{itemize}
\item there is no one-to-one correspondence between set systems and inequalities in their approach
although technically all inequalities assigned to a minimal sub-balanced collection are equivalent
(see \cite[Theorem\,3.9]{LBH12}),
\item their approach does not allow one to reveal one important relation of
{\em complementarity} among set systems which corresponds to the respective relation
of {\em conjugacy\/} between facet-defining inequalities for $\mathcal{E}(N)$ (see \cite[Lemma\,3.4]{KS19}).
The reader can find further details in Section~\ref{ssec.complementarity}.
\end{itemize}
That is why we believe our way of inequality description is more appropriate.
\end{remark}

\section{Properties of semi-balanced systems}\label{sec.properties}
In this section we discuss some structural relations among semi-balanced set systems.

\subsection{Complementarity of set systems}\label{ssec.complementarity}
A substantial fact about the cone of exact games is that its facet-defining inequalities come in pairs
of mutually conjugate inequalities, as shown already in \cite[Lemma~3.4]{KS19}. Therefore, the corresponding set systems
also come in pairs of mutually complementary systems.

\begin{defin}[conjugate inequality, complementary set system]~\label{def.complement}\rm \\
Assume $|N|\geq 2$. The {\em conjugate inequality\/} to the inequality \eqref{eq.theta-ineq}, that is, to
the inequality $0\leq \langle \theta,\m\rangle$ for $\m\in {\dv R}^{\caP}$
with a coefficient vector $\theta\in {\dv R}^{\caP}$, is the inequality
\begin{equation}
0\leq \langle \theta^{\refl},\m\rangle~~\mbox{for $\m\in {\dv R}^{\caP}$,}\qquad \mbox{where~~ $\theta^{\refl}(L):=\theta(N\setminus L)$ for any $L\subseteq N$.}
\label{eq.conjug}
\end{equation}
Given a non-trivial set system $\calS$ on $N$, its {\em complementary system\/} is the set system
$$
\calS^{\refl} := \{ N\setminus S\,:\ S\in\calS\}\,.
$$
\end{defin}
\smallskip

Of course, both concepts are relative to $N$. Here are the relevant observations.

\begin{lem}\label{lem.conjug}~\rm
Given $|N|\geq 2$, let $\calS$ be a non-trivial set system on $N$. Then $\calS$ is semi-balanced iff
$\calS^{\refl}$ is semi-balanced and $\Theta_{\calS^{\refl}}=\{ \theta^{\refl}\,:\ \theta\in\Theta_{\calS}\}$.
An analogous statement holds for balanced systems.
In particular, $\calS$ is min-semi-balanced iff $\calS^{\refl}$ is min-semi-balanced and the same
holds for min-balanced systems on $N$. Clearly, $\theta^{\refl}_{\calS}:=(\theta_{\calS})^{\refl}=\theta_{\calS^{\refl}}$ then.\\
Given $\theta\in {\dv R}^{\caP}$ with $\sum_{S\subseteq N} \theta (S)=0$ and
$\sum_{\emptyset\neq L\subset N} \theta(L)=-1$, the inequality \eqref{eq.theta-ineq}
is facet-defining for $\calE (N)$ iff its conjugate inequality \eqref{eq.conjug}
is facet-defining for $\calE (N)$.
\end{lem}

\begin{proof}
By Lemma~\ref{lem.semi-comb}, $\calS$ is semi-balanced if $r\cdot\chi_{N}=\sum_{S\in\calS} \lambda_{S}\cdot\chi_{S}$
with $r\in [0,1]$ and an affine semi-conic combination on the right-hand side (which has all its coefficients non-zero).
One can multiply that by $(-1)$ and add to that the equality $\chi_{N}=\sum_{S\in\calS} \lambda_{S}\cdot\chi_{N}$ to get
$(1-r)\cdot\chi_{N}=\sum_{S\in\calS} \lambda_{S}\cdot\chi_{N\setminus S}=\sum_{L\in\calS^{\refl}} \lambda_{N\setminus L}\cdot\chi_{L}$,
which means that $\calS^{\refl}$ is semi-balanced; one can then put $r^{\refl}=1-r$ and $\lambda^{\refl}_{L}=\lambda_{N\setminus L}$ for
$L\in\calS^{\refl}$. Thus, the relation $\Theta_{\calS^{\refl}}=\{ \theta^{\refl}\,:\ \theta\in\Theta_{\calS}\}$ follows from 
Definition~\ref{def.theta}. The same argument works for balanced systems: the linear combination is even conic then.
The relation of non-trivial systems $\calT\subseteq\calS$ iff $\calT^{\refl}\subseteq\calS^{\refl}$ then implies the consequences concerning
minimal such systems.

The last claim follows can be derived from Corollary~\ref{cor.facet-char} using Definition~\ref{def.aux-cone}.
Note that $\theta\in\tilde{\Theta}^{N}_{N} ~\Leftrightarrow~ \theta^{\refl}\in\tilde{\Theta}^{N}_{N}$, and, for every
$\emptyset\neq D\subset N$, $\theta\in\tilde{\Theta}^{N}_{D} ~\Leftrightarrow~ \theta^{\refl}\in\tilde{\Theta}^{N}_{N\setminus D}$,
which allows one to deduce $\theta\in\Delta ~\Leftrightarrow~ \theta^{\refl}\in\Delta$.
Since $\theta\mapsto \theta^{\refl}$ is a linear mapping one has
$\theta\in\ext(\Delta) ~\Leftrightarrow~ \theta^{\refl}\in\ext(\Delta)$ and the rest follows from Corollary~\ref{cor.facet-char}.
\end{proof}
\smallskip

Note that $T\in\calS$ is exceptional within a non-trivial set system $\calS$ on $N$ (see Definition~\ref{def.semi-bal}) iff $N\setminus T\in\calS^{\refl}$ is exceptional within $\calS^{\refl}$: given a combination $r\cdot\chi_{N}=\sum_{S\in\calS} \lambda_{S}\cdot\chi_{S}$
multiply it by $(-1)$ and add $(\sum_{S\in\calS} \lambda_{S})\cdot\chi_{N}$ to get
$(\sum_{S\in\calS} \lambda_{S}-r)\cdot\chi_{N}=\sum_{S\in\calS} \lambda_{S}\cdot\chi_{N\setminus S}=\sum_{R\in\calS^{\refl}} \lambda_{N\setminus R}\cdot\chi_{R}$. In particular, a non-trivial system $\calS$ on $N$ is purely min-semi-balanced iff
the same holds for its complementary system $\calS^{\refl}$.

\subsection{Basic classification of min-semi-balanced systems}\label{ssec.classification}
The uniqueness condition {\bf (d)} in Lemma~\ref{lem.min-semi-bal} on an affine combination
yielding a constant vector allows one to classify min-semi-balanced systems by the values of the constant.

\begin{lem}\label{lem.clasify}~\rm
Given $|N|\geq 2$ and a min-semi-balanced set system $\calS$ on $N$, let $\sum_{S\in\calS} \lambda_{S}\cdot\chi_{S}$
be the unique affine (semi-conic) combination yielding a constant vector $r\cdot\chi_{N}$, $r\in [0,1]$.\\
One has then $r=0$ iff there exists a min-balanced system $\calB$ on $M\subset N$, $|M|\geq 2$, such that $\calS=\calB\cup\{M\}$;
another equivalent condition is $\bigcup\calS\subset N$.\\
One has $r=1$ iff $\calS$ is a complementary system to a system $\calS^{\refl}$ with $\bigcup\calS^{\refl}\subset N$;
another equivalent condition is $\bigcap\calS\neq\emptyset$.\\
On the other hand, every min-balanced system $\calS$ on $N$ satisfies $0<r<1$.
\end{lem}

\begin{proof}
Recall that all coefficients $\lambda_{S}$, $S\in\calS$, are non-zero.
If $r=0$ then there is $T\in\calS$ with $\lambda_{T}<0$,
for otherwise $\lambda_{S}>0$ for $S\in\calS$ and $\calS\neq\emptyset$ contradict $\sum_{S\in\calS} \lambda_{S}\cdot\chi_{S}=\0$.
Given $i\in N\setminus T$ one has $\lambda_{S}>0$ whenever $i\in S\in\calS$ and
$0=\sum_{S\in\calS} \lambda_{S}\cdot\chi_{S}(i)=\sum_{S\in\calS: i\in S} \lambda_{S}$ implies that there is no
$S\in\calS$ with $i\in S$. Hence, $\bigcup\calS=T$ and $\sum_{S\in\calS\setminus\{T\}} \lambda_{S}\cdot\chi_{S}=(-\lambda_{T})\cdot\chi_{T}$.
This implies $\calS\setminus\{T\}\neq 0$ and, thus, forces $|T|\geq 2$. The latter equality also means that
$\calS\setminus\{T\}$ is balanced on $T$ and one can put $M:=T$ and $\calB:=\calS\setminus\{T\}$.
To show that $\calB$ is minimal assume for a contradiction that $\calD\subset\calB$ exists which is balanced on $T$, that is,
there are $\sigma_{S}>0$, $S\in\calD$, with $\sum_{S\in\calD} \sigma_{S}\cdot\chi_{S}=\chi_{T}$.
Hence, $\sum_{S\in\calD} \sigma_{S}\cdot\chi_{S}+ (-1)\cdot\chi_{T}=\0$
and, by Lemma~\ref{lem.semi-comb}, one has $-1+\sum_{S\in\calD} \sigma_{S}>0$ and
one can multiply it to get an affine combination (within ${\dv R}^{N}$) different from $\sum_{S\in\calS} \lambda_{S}\cdot\chi_{S}=\0$,
which contradicts the minimality of $\calS$. Thus, $\calB$ has to be min-balanced on $T$.
The existence of such $\calB$ then implies $\bigcup\calS\subset N$. Conversely, if $\bigcup\calS\subset N$ then the choice
of $i\in N\setminus\bigcup\calS$ gives $r=\sum_{S\in\calS:i\in S} \lambda_{S}=0$.

By Lemma~\ref{lem.conjug}, $\calS^{\refl}$ is also min-semi-balanced and the respective unique affine combination
for $\calS^{\refl}$ is $\sum_{L\in\calS^{\refl}} \lambda_{N\setminus L}\cdot\chi_{L}=(1-r)\cdot\chi_{N}$.
Thus, $r=1$ iff the previous case occurs for $\calS^{\refl}$. The formula $\bigcap\calS= N\setminus \bigcup\calS^{\refl}$ then gives the other equivalent condition.

The definition of a min-balanced system $\calS$ on $N$ implies $\bigcup\calS=N$. By Lemma~\ref{lem.conjug}, the same is true for the complementary system: $\bigcup\calS^{\refl}=N$, which means $\bigcap\calS=\emptyset$. Then use the previous two claims.
\end{proof}

Thus, one can distinguish at least three classes of min-semi-balanced systems $\calS$:
\begin{itemize}
\item those with $\bigcup\calS\subset N$, which are extensions of min-balanced systems on strict subsets,
\item those with $\bigcap\calS\neq\emptyset$, which can be viewed as their complementary systems, and
\item min-balanced systems on $N$, which satisfy both $\bigcup\calS=N$ and $\bigcap\calS=\emptyset$.
\end{itemize}
Nevertheless, as the next example shows, there is the fourth class of min-semi-balanced systems:
these satisfy both $\bigcup\calS=N$ and $\bigcap\calS=\emptyset$ but they are {\sf not\/} balanced on $N$.

\begin{example}\label{ex.4th-type}\rm
Put $N:=\{a,b,c,d\}$ and $\calS:=\{\,a,ab,bc,abd\,\}$.
The semi-conic combination
$$
1\cdot\chi_{a}+ (-1)\cdot\chi_{ab}+ 1\cdot\chi_{bc}+ 1\cdot\chi_{abd} = \chi_{N}
$$
implies that $\calS$ is semi-balanced on $N$. Since $\bigcup\calS=N$ and
$\{\chi_{S}\,:\ S\in\calS\}$ are linearly independent, by Lemma~\ref{lem.min-semi-bal}{\bf (c)}, $\calS$
is min-semi-balanced. Clearly, the unique linear combination yielding $\chi_{N}$ is not conic;
hence, $\calS$ is not balanced.
\end{example}

\begin{remark}\rm
By Lemma~\ref{lem.clasify}, every min-balanced set system $\calB$ on a strict subset $M\subset N$ leads to a semi-balanced system $\calB\cup \{M\}$ on $N$. Note in this context that $\calB$ itself is never semi-balanced on $N$. This is because the vectors $\{ \chi_{S}\,:\ S\in\calB\}$ are then linearly independent
\cite[Lemma~2.1]{KS19}. On the other hand, $\bigcup\calB=M\subset N$ implies that a contingent semi-conic combination of $\{ \chi_{S}\,:\ S\in\calB\}$ can only yield the zero constant vector $\0$, while the linear independence implies that only the zero linear combination yields $\0$.
\end{remark}

\subsection{Pictorial representation of min-semi-balanced systems}\label{ssec.picto-grams}
We propose to use certain special pictures to represent (permutational types of) minimal semi-balanced set systems.
In fact, our diagrams additionally encode the corresponding linear inequalities \eqref{eq.theta-ineq}.
Corollary~\ref{cor.rational} says that, given a min-semi-balanced system $\calS$ on $N$, the vector $\theta_{\calS}$ has rational components, that is, it is
given by an affine rational semi-conic combination $\sum_{S\in\calS} \lambda_{S}\cdot\chi_{S}=r\cdot\chi_{N}$.
One can multiply this by a natural number $\ell$ so that $\alpha_{S}:=\ell\cdot\lambda_{S}$, $S\in\calS$, become integers with no common prime
divisor. Then $\alpha_{N}:=\ell\cdot r\in {\dv Z}^{+}$ and one can introduce $\alpha_{\emptyset}:=-\alpha_{N}+\sum_{S\in\calS} \alpha_{S}\in {\dv Z}^{+}$ (use Lemma~\ref{lem.semi-comb}).
Thus, one gets
$$
\sum_{S\in\calS} \alpha_{S}\cdot\chi_{S}= \alpha_{\emptyset}\cdot\chi_{\emptyset} +\alpha_{N}\cdot\chi_{N}\,,\quad \mbox{where all the coefficients $\alpha_{S}$ are integers.}
$$
One can have $\alpha_{\emptyset}=0$ or $\alpha_{N}=0$, while the remaining coefficients are non-zero.
Provided there is $T\in\calS$ with $\lambda_{T}<0$ one can re-write that in the form
$$
\sum_{S\in\calS\setminus \{T\}} \alpha_{S}\cdot\chi_{S}= \alpha_{\emptyset}\cdot\chi_{\emptyset} + (-\alpha_{T})\cdot\chi_{T} +\alpha_{N}\cdot\chi_{N}\,\quad \mbox{with non-negative integers as coefficients.}
$$
A diagram representing a set system $\calS$ has the form of a pair of
two-dimensional arrays whose entries are colorful boxes; the arrays encode the sides of the above vector equality. The rows of these arrays correspond to the elements of the base set $N$ (= players); they are labeled if the diagram represents a particular set system $\calS$ and they are unlabeled if it represents a permutational type of such systems.

The columns of the arrays encode sets $S$ from the enlarged system $\calS\cup\{\emptyset, N\}$. Each of the sets has its own color; however, the black color is reserved for the grand coalition $N$, a fully blank (= white) column implicitly encodes the empty set and the grey color is reserved for a contingent set $T$ with a negative coefficient $\lambda_{T}<0$ (= an exceptional set in $\calS$).
The other sets from $\calS$ have bright colors then.
The column representing a set $S$ has boxes of the respective color just in rows corresponding to elements of $S$.
To express the value of the respective coefficient $\alpha_{S}\in {\dv Z}$ in the inequality the respective column
is repeated $|\alpha_{S}|$-times.

The left array is composed of columns which correspond to sets with positive coefficient $\lambda_{S}>0$, $S\in\calS$, while the array on the right-hand side has either fully black columns, fully blank (= white) columns and possibly
columns containing grey boxes.

\begin{figure}[h]
\setlength{\unitlength}{1mm}
\begin{center}
\scalebox{0.7}{
\begin{picture}(150,58)
\thinlines %
\put(3,5){\makebox(0,0){\Large $e$}}%
\put(3,15){\makebox(0,0){\Large $d$}}%
\put(3,25){\makebox(0,0){\Large $c$}}%
\put(3,35){\makebox(0,0){\Large $b$}}%
\put(3,45){\makebox(0,0){\Large $a$}}%
\put(75,25){\makebox(0,0){\Large $=$}}%
\put(10,0){\line(1,0){60}} %
\put(10,10){\line(1,0){60}} %
\put(10,20){\line(1,0){60}} %
\put(10,30){\line(1,0){60}} %
\put(10,40){\line(1,0){60}} %
\put(10,50){\line(1,0){60}} %
\put(10,0){\line(0,1){50}} %
\put(20,0){\line(0,1){50}} %
\put(30,0){\line(0,1){50}} %
\put(40,0){\line(0,1){50}} %
\put(50,0){\line(0,1){50}} %
\put(60,00){\line(0,1){50}} %
\put(70,0){\line(0,1){50}} %
%
%
\put(80,0){\line(1,0){60}} %
\put(80,10){\line(1,0){60}} %
\put(80,20){\line(1,0){60}} %
\put(80,30){\line(1,0){60}} %
\put(80,40){\line(1,0){60}} %
\put(80,50){\line(1,0){60}} %
\put(80,0){\line(0,1){50}} %
\put(90,0){\line(0,1){50}} %
\put(100,0){\line(0,1){50}} %
\put(110,0){\line(0,1){50}} %
\put(120,0){\line(0,1){50}} %
\put(130,00){\line(0,1){50}} %
\put(140,0){\line(0,1){50}} %
%
%
\put(10.3,40.3){\textcolor{green}{\rule{9.4mm}{9.4mm}}} %
\put(10.3,20.3){\textcolor{green}{\rule{9.4mm}{9.4mm}}} %
\put(20.3,30.3){\textcolor{myorange}{\rule{9.4mm}{9.4mm}}} %
\put(20.3,20.3){\textcolor{myorange}{\rule{9.4mm}{9.4mm}}} %
\put(30.3,40.3){\textcolor{blue}{\rule{9.4mm}{9.4mm}}} %
\put(30.3,30.3){\textcolor{blue}{\rule{9.4mm}{9.4mm}}} %
\put(30.3,10.3){\textcolor{blue}{\rule{9.4mm}{9.4mm}}} %
\put(40.3,40.3){\textcolor{blue}{\rule{9.4mm}{9.4mm}}} %
\put(40.3,30.3){\textcolor{blue}{\rule{9.4mm}{9.4mm}}} %
\put(40.3,10.3){\textcolor{blue}{\rule{9.4mm}{9.4mm}}} %
%
\put(50.3,40.3){\textcolor{red}{\rule{9.4mm}{9.4mm}}} %
\put(50.3,30.3){\textcolor{red}{\rule{9.4mm}{9.4mm}}} %
\put(50.3,0.3){\textcolor{red}{\rule{9.4mm}{9.4mm}}} %
\put(60.3,40.3){\textcolor{red}{\rule{9.4mm}{9.4mm}}} %
\put(60.3,30.3){\textcolor{red}{\rule{9.4mm}{9.4mm}}} %
\put(60.3,0.3){\textcolor{red}{\rule{9.4mm}{9.4mm}}} %
\put(90.3,40.3){\textcolor{mygrey}{\rule{9.4mm}{9.4mm}}} %
\put(90.3,30.3){\textcolor{mygrey}{\rule{9.4mm}{9.4mm}}} %
\put(100.3,40.3){\textcolor{mygrey}{\rule{9.4mm}{9.4mm}}} %
\put(100.3,30.3){\textcolor{mygrey}{\rule{9.4mm}{9.4mm}}} %
\put(110.3,40.3){\textcolor{mygrey}{\rule{9.4mm}{9.4mm}}} %
\put(110.3,30.3){\textcolor{mygrey}{\rule{9.4mm}{9.4mm}}} %
\put(120.3,40.3){\textcolor{black}{\rule{9.4mm}{9.4mm}}} %
\put(120.3,30.3){\textcolor{black}{\rule{9.4mm}{9.4mm}}} %
\put(120.3,20.3){\textcolor{black}{\rule{9.4mm}{9.4mm}}} %
\put(120.3,10.3){\textcolor{black}{\rule{9.4mm}{9.4mm}}} %
\put(120.3,0.3){\textcolor{black}{\rule{9.4mm}{9.4mm}}} %
\put(130.3,40.3){\textcolor{black}{\rule{9.4mm}{9.4mm}}} %
\put(130.3,30.3){\textcolor{black}{\rule{9.4mm}{9.4mm}}} %
\put(130.3,20.3){\textcolor{black}{\rule{9.4mm}{9.4mm}}} %
\put(130.3,10.3){\textcolor{black}{\rule{9.4mm}{9.4mm}}} %
\put(130.3,0.3){\textcolor{black}{\rule{9.4mm}{9.4mm}}} %
\end{picture}
}
\end{center}
\caption{A picture representing the set system from Example~\ref{ex.pictogram}.\label{fig.pictogram}}
\end{figure}
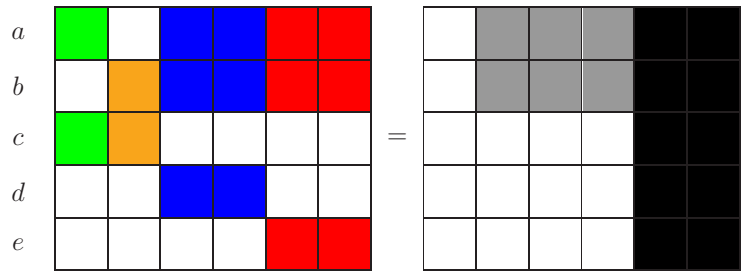

\begin{example}\label{ex.pictogram}\rm
Put $N:=\{a,b,c,d,e\}$ and $\calS:=\{\,ab,ac,bc,abd,abe\,\}$.
The equality relation
$$
1\cdot\chi_{ac}+ 1\cdot\chi_{bc}+ 2\cdot\chi_{abd}+ 2\cdot\chi_{abe} = 1\cdot\chi_{\emptyset} +3\cdot\chi_{ab} +2\cdot\chi_{N}
$$
allows one to observe that $\calS$ is semi-balanced on $N$.
As the vectors $\{\chi_{S}\,:\ S\in\calS\}$ are linearly independent $\calS$ is minimal. A picture representing this set system is in Figure~\ref{fig.pictogram}.
\end{example}

Note that, in any row, the numbers of boxes in the left and right array coincide: this is because of
the equality $\sum_{S\in\calS} \alpha_{S}\cdot\chi_{S}(i)=\alpha_{N}\cdot\chi_{N}(i)$
for any $i\in N$.
Another interesting observation is that the diagram for the complementary system
$\calS^{\refl}$ can easily by obtained by ``reflection" from the diagram for $\calS$:
the boxes are interchanged with non-boxes and the colors for columns are kept, under a convention
that the color for blank columns is black.

One can also easily recognize on basis of the diagram for $\calS$ to which of the four basic classes it belongs.
Systems $\calS$ with $\bigcup\calS\subset N$, that is, those with the constant $r=0$, have no black column.
Their complementary systems $\calS$ with $\bigcap\calS\neq\emptyset$, that is, those with the constant $r=1$, have no
blank (= white) column. The min-balanced systems on $N$ have no grey column and other min-semi-balanced
systems have blank, grey and black columns.

\section{Purely min-semi-balanced systems}\label{sec.purely-min}
Let us first discuss the situation when only two players exist, that is, $|N|=2$.
Then, as explained below Lemma~\ref{lem.semi-comb}, there is no purely semi-balanced system on $N$.
In fact, the only semi-balanced system over $N=\{a,b\}$ is the min-balanced system
$\calB:=\{\, a,b\,\}$. Thus, by Corollary~\ref{cor.exa-char}, a game $\m\in {\dv R}^{\caP}$, $\m(\emptyset)=0$,
is exact iff $\langle \theta_{\calB},\m\rangle\geq 0$, and the cone of exact games
is specified by a single inequality $\m (ab)-\m(a)-\m(b)\geq 0$.

\subsection{How to get purely min-semi-balanced systems}\label{ssec.get-purely}
Therefore, in the sequel we limit our attention to a non-trivial case $|N|\geq 3$,
when there exist purely semi-balanced systems on $N$. We first establish
some relation between min-balanced and purely min-semi-balanced systems.

\begin{lem}\label{lem.balanc-corresp}~\rm
Assume $|N|\geq 3$.
If $\calB$ is a min-balanced set system on $N$ such that $|\calB|\geq 3$ and $Z\in\calB$ then
$Y:=N\setminus Z$ is {\sf not} in $\calB$, the set system
$\calS:=(\calB\setminus\{Z\})\cup\{Y\}$ is purely min-semi-balanced on $N$ and $Y$ is the exceptional set within $\calS$.\\[0.2ex]
Conversely, if $\calS$ is a (purely) min-semi-balanced system on $N$ and $Y\in\calS$ the exceptional
set within $\calS$ then $Z:=N\setminus Y$ is {\sf not} in $\calS$ and $\calB:=(\calS\setminus\{Y\})\cup\{Z\}$
is a min-balanced system on $N$ such that $|\calB|\geq 3$.
\end{lem}

\begin{proof}
{\bf I.}\, Let $\calB$ be a balanced system on $N$ and $Z\in\calB$ such that $Y:=N\setminus Z\not\in\calB$.
Then $\calS:=(\calB\setminus\{Z\})\cup\{Y\}$ is semi-balanced on $N$. Indeed, there exists a conic combination $\sum_{S\in\calB} \lambda_{S}\cdot\chi_{S}=\chi_{N}$, $\lambda_{S}>0$ for $S\in\calB$.
We add $-\lambda_{Z}\cdot\chi_{N}$ to that and obtain a semi-conic combination $\sum_{S\in\calB\setminus\{Z\}} \lambda_{S}\cdot\chi_{S} + (-\lambda_{Z})\cdot\chi_{Y}=(1-\lambda_{Z})\cdot\chi_{N}$
yielding a constant vector.

{\bf II.}\, Analogously, given a semi-balanced system $\calS$ over $N$ with an exceptional set $Y\in\calS$ and $Z:=N\setminus Y\not\in\calS$
the set system $\calB:=(\calS\setminus\{Y\})\cup\{Z\}$ is balanced on $N$. Indeed, given a semi-conic combination $\sum_{S\in\calS} \lambda_{S}\cdot\chi_{S}=r\cdot\chi_{N}$ with $r\in {\dv R}$, where
$\lambda_{Y}<0$ and $\lambda_{S}>0$ for $S\in\calS\setminus\{Y\}$
we add $-\lambda_{Y}\cdot\chi_{N}$ to that and get
$\sum_{S\in\calS\setminus\{Y\}} \lambda_{S}\cdot\chi_{S} + (-\lambda_{Y})\cdot\chi_{Z}=(r-\lambda_{Y})\cdot\chi_{N}$, which is a conic combination yielding a constant vector. Then use Lemma~\ref{lem.semi-comb}.

{\bf III.}\, Let $\calB$ be a min-balanced system on $N$ with $|\calB|\geq 3$ and $Z\in\calB$. Then $Y:=N\setminus Z\not\in\calB$
as otherwise $\calD:=\{Y,Z\}\subset\calB$ is a balanced system
contradicting the minimality of $\calB$. Step {\bf I.\/} implies that
$\calS:=(\calB\setminus\{Z\})\cup\{Y\}$ is semi-balanced on $N$, specifically,
that there exists a semi-conic combination $\sum_{S\in\calS} \lambda_{S}\cdot\chi_{S}$
yielding a constant vector in ${\dv R}^{N}$ with $\lambda_{Y}<0$.
To show that $\calS$ is minimal assume for a contradiction that a min-semi-balanced system $\calS^{\prime}\subset\calS$ exists.
Then necessarily $Y\in\calS^{\prime}$ as otherwise $\calS^{\prime}\subset\calB$ contradicts the minimality of $\calB$ (use Corollary~\ref{cor.min-balanced}).
Let $\sum_{S\in\calS^{\prime}} \sigma_{S}\cdot\chi_{S}$ be the (unique) affine semi-conic combination yielding a constant vector in ${\dv R}^{N}$.
Observe that $\sigma_{Y}<0$ as otherwise one can put $\sigma_{S}:=0$ for $S\in\calS\setminus\calS^{\prime}$,
$\tau_{S}:=\alpha\cdot\sigma_{S}+(1-\alpha)\cdot\lambda_{S}$ for
$S\in\calS$ with $\alpha:=-\lambda_{Y}\cdot (\sigma_{Y}-\lambda_{Y})^{-1}\in (0,1)$ to
get a semi-conic combination  $\sum_{S\in\calS} \tau_{S}\cdot\chi_{S}$
yielding a constant vector where $\tau_{Y}=0$; this means that
$\calE:=\{S\in\calS\,:\ \tau_{S}\neq 0\}\subset \calB$ is a semi-balanced system
which fact contradicts the minimality of $\calB$, by Corollary~\ref{cor.min-balanced}.
Thus, $Y$ has to be an exceptional set within $\calS^{\prime}$ and, by step {\bf II.} applied to $\calS^{\prime}$,
the system $\calB^{\prime}:=(\calS^{\prime}\setminus\{Y\})\cup\{Z\}\subset\calB$ is a balanced system
on $N$ contradicting the minimality of $\calB$. Thus, the first claim in Lemma~\ref{lem.balanc-corresp}
has been verified.

{\bf IV.}\, Let $\calS$ be a min-semi-balanced system on $N$ with the exceptional set $Y\in\calS$ within it.
Then $Z:=N\setminus Y\not\in\calS$
as otherwise $\calD:=\{Y,Z\}\subset\calS$ is a balanced system
contradicting the minimality of $\calS$ (note that $|\calS|\geq 3$ by Lemma~\ref{lem.semi-comb}). Step {\bf II.\/} implies that $\calB:=(\calS\setminus\{Y\})\cup\{Z\}$ is balanced on $N$. To show that $\calB$ is minimal assume for a contradiction that a balanced system $\calB^{\prime}\subset\calB$ on $N$ exists. Then necessarily $Z\in\calB^{\prime}$ as otherwise $\calB^{\prime}\subset\calS$ contradicts the minimality of $\calS$. By step {\bf I.} applied to $\calB^{\prime}$, the set system $\calS^{\prime}:=(\calB^{\prime}\setminus\{Z\})\cup\{Y\}\subset\calS$
is semi-balanced on $N$, which contradicts the minimality of $\calS$. Of course, $|\calB|=|\calS|\geq 3$,
which concludes the proof of the second claim in Lemma~\ref{lem.balanc-corresp}.
\end{proof}
\medskip

Thus, by Lemma~\ref{lem.balanc-corresp}, there is one-to-many correspondence $\calB \leftrightarrow \calS$ between
min-balanced systems $\calB$ on $N$ satisfying $|\calB|\geq 3$ and purely min-semi-balanced systems $\calS$ on $N$
which is realized by the mutual exchange of a set $Z\in\calB$ and of its complement $Y\in\calS$. 
It allows one to generate a complete list of min-semi-balanced systems on $N$ on basis of the list of all balanced systems on $N$.
Note that the fact that balanced systems on $N$ induce in this way other semi-balanced systems on $N$  has already been recognized in
\cite[\S\,4,\,Theorem~4.4]{LBH12}.  Nevertheless, the above correspondence was not revealed there in its full scope for the reason mentioned in Remark~\ref{rem.lohmnamm}.

Let us remark in this context that the discussed transition from a min-balanced system $\calB$
to a purely min-semi-balanced system $\calS$ (and back) can be recognized easily on basis of 
their diagrams/pictures from Section~\ref{ssec.picto-grams}. Indeed, if a diagram represents
a min-balanced system $\calB$ and has $\alpha_{Z}$ columns representing a set $Z\in\calB$ then
the diagram representing $\calS:=(\calB\setminus\{Z\})\cup\{Y\}$, where $Y:=N\setminus Z$, can be
obtained from it by the removal those $\alpha_{Z}$ columns of bright color from the left array,
$\alpha_{Z}$ blank columns and $\alpha_{Z}$ black columns from the right array and by adding
$\alpha_{Z}$ grey columns representing the set $Y$ to the right array. Of course, the transition back
can be done by an inverse operation with the diagrams. The following example illustrates the procedure.

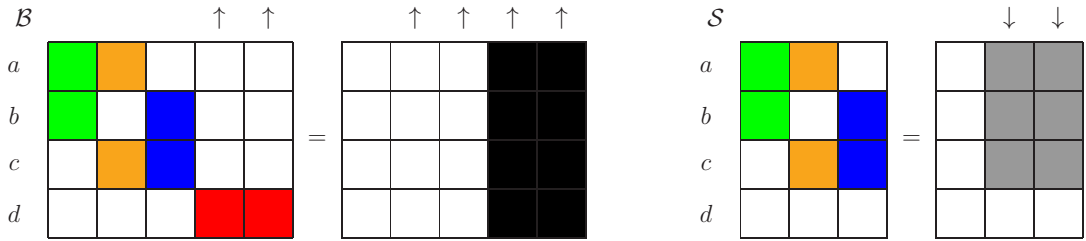
\begin{figure}[h]
\setlength{\unitlength}{1mm}
\begin{center}
\scalebox{0.65}{
\hspace*{2mm}
\begin{picture}(130,50)
\thinlines %
\put(3,5){\makebox(0,0){\Large $d$}}%
\put(3,15){\makebox(0,0){\Large $c$}}%
\put(3,25){\makebox(0,0){\Large $b$}}%
\put(3,35){\makebox(0,0){\Large $a$}}%
\put(65,20){\makebox(0,0){\Large $=$}}%
\put(5,45){\makebox(0,0){\Large $\calB$}}%
\put(45,45){\makebox(0,0){\Large $\uparrow$}}%
\put(55,45){\makebox(0,0){\Large $\uparrow$}}%
\put(85,45){\makebox(0,0){\Large $\uparrow$}}%
\put(95,45){\makebox(0,0){\Large $\uparrow$}}%
\put(105,45){\makebox(0,0){\Large $\uparrow$}}%
\put(115,45){\makebox(0,0){\Large $\uparrow$}}%
\put(10,0){\line(1,0){50}} %
\put(10,10){\line(1,0){50}} %
\put(10,20){\line(1,0){50}} %
\put(10,30){\line(1,0){50}} %
\put(10,40){\line(1,0){50}} %
\put(10,0){\line(0,1){40}} %
\put(20,0){\line(0,1){40}} %
\put(30,0){\line(0,1){40}} %
\put(40,0){\line(0,1){40}} %
\put(50,0){\line(0,1){40}} %
\put(60,0){\line(0,1){40}} %
%
%
\put(70,0){\line(1,0){50}} %
\put(70,10){\line(1,0){50}} %
\put(70,20){\line(1,0){50}} %
\put(70,30){\line(1,0){50}} %
\put(70,40){\line(1,0){50}} %
\put(70,0){\line(0,1){40}} %
\put(80,0){\line(0,1){40}} %
\put(90,0){\line(0,1){40}} %
\put(100,0){\line(0,1){40}} %
\put(110,0){\line(0,1){40}} %
\put(120,00){\line(0,1){40}} %
%
%
\put(10.3,30.3){\textcolor{green}{\rule{9.4mm}{9.4mm}}} %
\put(10.3,20.3){\textcolor{green}{\rule{9.4mm}{9.4mm}}} %
\put(20.3,30.3){\textcolor{myorange}{\rule{9.4mm}{9.4mm}}} %
\put(20.3,10.3){\textcolor{myorange}{\rule{9.4mm}{9.4mm}}} %
\put(30.3,20.3){\textcolor{blue}{\rule{9.4mm}{9.4mm}}} %
\put(30.3,10.3){\textcolor{blue}{\rule{9.4mm}{9.4mm}}} %
\put(40.3,0.3){\textcolor{red}{\rule{9.4mm}{9.4mm}}} %
\put(50.3,0.3){\textcolor{red}{\rule{9.4mm}{9.4mm}}} %
\put(100.3,30.3){\textcolor{black}{\rule{9.4mm}{9.4mm}}} %
\put(100.3,20.3){\textcolor{black}{\rule{9.4mm}{9.4mm}}} %
\put(100.3,10.3){\textcolor{black}{\rule{9.4mm}{9.4mm}}} %
\put(100.3,0.3){\textcolor{black}{\rule{9.4mm}{9.4mm}}} %
\put(110.3,30.3){\textcolor{black}{\rule{9.4mm}{9.4mm}}} %
\put(110.3,20.3){\textcolor{black}{\rule{9.4mm}{9.4mm}}} %
\put(110.3,10.3){\textcolor{black}{\rule{9.4mm}{9.4mm}}} %
\put(110.3,0.3){\textcolor{black}{\rule{9.4mm}{9.4mm}}} %
\end{picture}
\hspace*{9mm}
\begin{picture}(90,50)
\thinlines %
\put(3,5){\makebox(0,0){\Large $d$}}%
\put(3,15){\makebox(0,0){\Large $c$}}%
\put(3,25){\makebox(0,0){\Large $b$}}%
\put(3,35){\makebox(0,0){\Large $a$}}%
\put(45,20){\makebox(0,0){\Large $=$}}%
\put(5,45){\makebox(0,0){\Large $\calS$}}%
\put(65,45){\makebox(0,0){\Large $\downarrow$}}%
\put(75,45){\makebox(0,0){\Large $\downarrow$}}%
\put(10,0){\line(1,0){30}} %
\put(10,10){\line(1,0){30}} %
\put(10,20){\line(1,0){30}} %
\put(10,30){\line(1,0){30}} %
\put(10,40){\line(1,0){30}} %
\put(10,0){\line(0,1){40}} %
\put(20,0){\line(0,1){40}} %
\put(30,0){\line(0,1){40}} %
\put(40,0){\line(0,1){40}} %
%
%
\put(50,0){\line(1,0){30}} %
\put(50,10){\line(1,0){30}} %
\put(50,20){\line(1,0){30}} %
\put(50,30){\line(1,0){30}} %
\put(50,40){\line(1,0){30}} %
\put(50,0){\line(0,1){40}} %
\put(60,0){\line(0,1){40}} %
\put(70,0){\line(0,1){40}} %
\put(80,0){\line(0,1){40}} %
%
%
\put(10.3,30.3){\textcolor{green}{\rule{9.4mm}{9.4mm}}} %
\put(10.3,20.3){\textcolor{green}{\rule{9.4mm}{9.4mm}}} %
\put(20.3,30.3){\textcolor{myorange}{\rule{9.4mm}{9.4mm}}} %
\put(20.3,10.3){\textcolor{myorange}{\rule{9.4mm}{9.4mm}}} %
\put(30.3,20.3){\textcolor{blue}{\rule{9.4mm}{9.4mm}}} %
\put(30.3,10.3){\textcolor{blue}{\rule{9.4mm}{9.4mm}}} %
\put(60.3,30.3){\textcolor{mygrey}{\rule{9.4mm}{9.4mm}}} %
\put(60.3,20.3){\textcolor{mygrey}{\rule{9.4mm}{9.4mm}}} %
\put(60.3,10.3){\textcolor{mygrey}{\rule{9.4mm}{9.4mm}}} %
\put(70.3,30.3){\textcolor{mygrey}{\rule{9.4mm}{9.4mm}}} %
\put(70.3,20.3){\textcolor{mygrey}{\rule{9.4mm}{9.4mm}}} %
\put(70.3,10.3){\textcolor{mygrey}{\rule{9.4mm}{9.4mm}}} %
\end{picture}
}
\end{center}
\caption{Pictures representing the set systems $\calB$ and $\calS$ from Example~\ref{ex.transition}.\label{fig.transition}}
\end{figure}

\begin{example}\label{ex.transition}\rm
Put $N:=\{a,b,c,d\}$ and $\calB:=\{\,ab,ac,bc,d\,\}$. One can observe that $\calB$ is min-balanced on $N$
using the equality relation
$$
\frac{1}{2}\cdot\chi_{ab}+ \frac{1}{2}\cdot\chi_{ac}+ \frac{1}{2}\cdot\chi_{bc}+ 1\cdot\chi_{d}= \chi_{N}\,.
$$
The choice of a set $Z:=d$ from $\calB$ leads, by Lemma~\ref{lem.balanc-corresp}, to a min-semi-balanced
system $\calS:=(\calB\setminus\{Z\})\cup\{Y\}=\{\,ab,ac,bc,abc\,\}$. The diagrams for both systems are
in Figure~\ref{fig.transition}.
We observe that two columns representing $Z$, $\emptyset$ and $N$ from the diagram for $\calB$ are missing
in the diagram for $\calS$ and replaced there by two columns representing $Y:=N\setminus Z= abc$.
\end{example}

\subsection{The case of balanced systems}\label{ssec.rid-balanced}
We now show that the inequalities corresponding to balanced systems on $N$, $|N|\geq 3$, are superfluous.
The next observation follows from Lemma~\ref{lem.balanc-corresp}.

\begin{cor}\label{cor.balanc-decomp}~\rm
Given $|N|\geq 3$, let $\calB$ be a min-balanced set system on $N$ with $|\calB|\geq 3$. Then
$\theta_{\calB}$ is a (non-trivial) convex combination of $\theta_{\calD}$ for a min-balanced
system $\calD$ on $N$ with $|\calD|=2$ and of $\theta_{\calS}$ for a purely min-semi-balanced system
$\calS$ on $N$ with $|\calS|=|\calB|$.
\end{cor}

\begin{proof}
We take $Z\in\calB$, put $Y:=N\setminus Z$ and observe, by
Lemma~\ref{lem.balanc-corresp}, that $\calS:=(\calB\setminus\{Z\})\cup\{Y\}$ is min-semi-balanced on $N$.
Consider the unique affine conic combination $\sum_{S\in\calB} \lambda_{S}\cdot\chi_{S}$ yielding a constant vector in ${\dv R}^{N}$ and
the semi-conic combination $\sum_{S\in\calS} \lambda_{S}\cdot\chi_{S}$ such that $\lambda_{Y}=-\lambda_{Z}$; it also
yields a constant vector in ${\dv R}^{N}$ by the argument from step {\bf I.} in the proof of Lemma~\ref{lem.balanc-corresp}.
Note that one has $1-2\lambda_{Z}>0$ by Lemma~\ref{lem.semi-comb}
applied to the latter combination. Thus, $(1-2\lambda_{Z})^{-1}$-multiple of it is an affine combination.
Put $\calD:=\{ Y,Z\}$, which is a min-balanced system on $N$, and the respective affine conic combination
$\frac{1}{2}\cdot\chi_{Y}+\frac{1}{2}\cdot\chi_{Z}$ yields a constant vector in ${\dv R}^{N}$. Then
\begin{eqnarray*}
\lefteqn{\hspace*{-4mm}\sum_{S\in\calB} \lambda_{S}\cdot\chi_{S} ~=~ 2\lambda_{Z}\cdot [\,\frac{1}{2}\cdot\chi_{Y}+\frac{1}{2}\cdot\chi_{Z} \,]}\\
~~ &+& (1-2\lambda_{Z})\cdot [\,\sum_{S\in\calB\setminus\{Z\}} (1-2\lambda_{Z})^{-1}\cdot \lambda_{S}\cdot\chi_{S} +
(1-2\lambda_{Z})^{-1}\cdot (-\lambda_{Z})\cdot \chi_{Y} \,]\,,
\end{eqnarray*}
and using \eqref{eq.theta} derive that $\theta_{\calB}=2\lambda_{Z}\cdot\theta_{\calD} +(1-2\lambda_{Z})\cdot \theta_{\calS}$.
\end{proof}

Corollary~\ref{cor.balanc-decomp} says that the vector $\theta_{\calB}$ for a min-balanced
system $\calB$ with $|\calB|\geq 3$ is a convex combination of vectors for min-semi-balanced systems of
cardinality at most $|\calB|$. The vector $\theta_{\calD}$ for a min-balanced system $\calD$ on $N$ with $|\calD|=2$
can also be written as a convex combination of other vectors. Nevertheless, the difference is that
the summands in the combination correspond to set systems of higher cardinality.

\begin{lem}\label{lem.two-balance-decompose}~\rm
Assume $|N|\geq 3$. If $\calD$ is a min-balanced system on $N$ with $|\calD|=2$ then
$\theta_{\calD}$ is a convex combination of $\theta_{\calS}$ and $\theta_{\calT}$ for purely min-semi-balanced
systems $\calS$ and $\calT$ on $N$.
\end{lem}

\begin{proof}
We have $\calD=\{Y,Z\}$ where $Y\cap Z=\emptyset$, $Y\cup Z=N$, $\emptyset\neq Z$ and $|Y|\geq 2$.
Thus, there exists a set $\emptyset\neq R\subset Y$ and one can put
$\calS:=\{\, Z,R,Z\cup R\,\}$ and $\calT:=\{\, Z\cup R, Y, R\,\}$.
The equalities $\chi_{Z}+\chi_{R}-\chi_{Z\cup R}=\0$
and $\chi_{Z\cup R}+\chi_{Y}-\chi_{R}=\chi_{N}$ together with affine/linear independence of involved vectors
allows one to observe using Lemma~\ref{lem.min-semi-bal} that both $\calS$ and $\calT$ is a min-semi-balanced system on $N$.
The respective affine semi-conic combinations are related as follows:
$$
\,[\,\frac{1}{2}\cdot\chi_{Y}+\frac{1}{2}\cdot\chi_{Z}\,] \,=\,
\frac{1}{2}\cdot [\,1\cdot\chi_{Z}+1\cdot\chi_{R}+(-1)\cdot\chi_{Z\cup R}\,]
+
\frac{1}{2}\cdot [\,1\cdot\chi_{Z\cup R}+1\cdot\chi_{Y}+(-1)\cdot\chi_{R}\,]\,,
$$
which equality implies using \eqref{eq.theta} that $\theta_{\calD}=\frac{1}{2}\cdot\theta_{\calS} +\frac{1}{2}\cdot\theta_{\calT}$.
\end{proof}

The previous two results allow one to derive the following conclusion.

\begin{cor}\label{cor.balanced-superfluous}~\rm
If $|N|\geq 3$ and $\calB$ is a min-balanced system on $N$ then the inequality given by $\theta=\theta_{\calB}$
is {\sf not\/} facet-defining for the exact cone $\calE(N)$.
\end{cor}

Recall that the fact that balanced systems provide superfluous inequalities has already been shown in
\cite[\S\,5,\,Theorem~5.4]{LBH12}; we give our short proof for the sake of completeness.

\begin{proof}
By combining Corollary~\ref{cor.balanc-decomp} with Lemma~\ref{lem.two-balance-decompose}
observe that every vector $\theta_{\calB}$ for a min-balanced system $\calB$ on $N$ is a non-trivial convex
combination of vectors $\theta_{\calS}$ for purely min-semi-balanced systems $\calS$ on $N$.
Then apply Corollary~\ref{cor.facet-char} to observe that the inequality \eqref{eq.theta-ineq}
with $\theta=\theta_{\calB}$ is not facet-defining for $\calE(N)$.
\end{proof}

\section{Indecomposable semi-balanced systems}\label{sec.indecomposable-min}
Nevertheless, even purely min-semi-balanced systems can induce superfluous inequalities.
We give a simple sufficient (combinatorial) condition for that.

\begin{defin}[decomposition, indecomposable system]~\label{def.decomposition}\rm \\
Assume $|N|\geq 3$. Given a purely min-semi-balanced set system $\calS$ on $N$ and $E\subseteq N$
with $E\not\in(\calS\cup\{\emptyset,N\})$ we say that $E$ yields
a {\em  decomposition\/} of\/ $\calS$ if $E$ is an exceptional set within $\calS\cup\{E\}$
(see Definition~\ref{def.semi-bal}).
A purely min-semi-balanced set system on $N$ will be called {\em indecomposable\/}
if it has no  decomposition.
\end{defin}

Recall from Section~\ref{ssec.complementarity} that a set is exceptional within
a non-trivial set system iff its complement is exceptional within its complementary
system. This implies that $\calS$ has a decomposition iff its complementary
system $\calS^{\refl}$ has a  decomposition. Therefore, $\calS$ is
indecomposable iff the same holds for $\calS^{\refl}$.
Our main result follows from the following lemma; its technical proof is moved to Appendix~\ref{app.lem-main-result}.

\begin{lem}\label{lem.main-result}~\rm
Let $\calS$ be a purely min-semi-balanced system on $N$, $|N|\geq 3$, with exceptional set $T\in\calS$
and $\calW:= \caP\setminus (\{\emptyset,N\}\cup\calS)$.
Then the following conditions are equivalent:
\begin{description}
\item[(a)] $\theta_{\calS}\not\in \ext(\Delta)$, (see Definition~\ref{def.aux-cone})
\item[(b)] there exists a convex combination $\theta_{\calS}=\sum_{D\in\calW\cup\{T\}} \alpha_{D}\cdot\theta^{D}$
where $\alpha_{T}<1$ and $\theta^{D}\in\tilde{\Theta}^{N}_{D}$ whenever $\alpha_{D}>0$,
\item[(c)] the set $\Delta(\calS):=\conv(\bigcup_{D\in\calW} \tilde{\Theta}^{N}_{D}) \cap \{\,\theta\in {\dv R}^{\caP}\,:\
\theta(W)\geq 0 ~\mbox{for $W\in\calW$}\,\}$ is non-empty,
\item[(d)] there exists $E\in\calW$ such that $E$ is exceptional in $\calS\cup\{E\}$, (= $\calS$ has
a decomposition)
\item[(e)] there exists $E\in\calW$ such that $E$ is exceptional in $(\calS\setminus\{T\})\cup\{E\}$,
\item[(f)] a min-semi-balanced system $\calD$ on $N$ exists such that $\calD\setminus\calS=\{E\}$
for some $E\in\calW$, the set $E$ is exceptional within $\calD$, and $T\not\in\calD$,
\item[(g)] a min-semi-balanced system $\calD$ on $N$ exists with an exceptional set $E$ and $\calD\setminus\calS=\{E\}$.
\end{description}
\end{lem}
\medskip

Now, we are ready to state our main result.

\begin{thm}\label{thm.main-result}
Given $|N|\geq 3$ the inequality $0\leq \langle \theta,\m\rangle$ for \mbox{$\m\in {\dv R}^{\caP}$}
with a coefficient vector $\theta\in {\dv R}^{\caP}$, where $\sum_{S\subseteq N} \theta (S)=0$ and
$\sum_{L:\,\emptyset\neq L\subset N} \theta(L)=-1$, is facet-defining for $m\in\calE (N)$ iff
$\theta=\theta_{\calS}$ for an indecomposable min-semi-balanced system $\calS$ on $N$.
\end{thm}

\begin{proof}
By Corollary~\ref{cor.facet-char}, facet-defining inequalities $0\leq\langle\theta,\m\rangle$ for $\m\in {\calE}(N)$
correspond to vertices $\theta$ of the polytope $\Delta$. Each $\theta\in\ext (\Delta)$ must be a vertex of
$\tilde{\Theta}^{N}_{D}$ for some $\emptyset\neq D\subseteq N$ (see Section~\ref{ssec.polyhedrals}).
By 
Lemma~\ref{lem.extr-theta}, every vertex $\theta\in\ext(\tilde{\Theta}^{N}_{D})$
has the form $\theta=\theta_{\calS}$ for a min-semi-balanced system $\calS$ on $N$.
Corollary~\ref{cor.balanced-superfluous} excludes the case that $\calS$ is min-balanced.
In case of a purely min-semi-balanced system $\calS$ one applies Lemma~\ref{lem.main-result}, the equivalence of negations
$\neg\mbox{\bf (a)}\Leftrightarrow\neg\mbox{\bf (d)}$, which says that $\theta=\theta_{\calS}\in\ext (\Delta)$
iff $\calS$ is indecomposable.
\end{proof}

Note in this context that, by Lemma~\ref{lem.main-result}{\bf (f)}, a purely min-semi-balanced system
$\calS$ on $N$ is indecomposable iff there is {\sf no\/} min-semi-balanced system $\calD$ with $T\not\in\calD$,
$\calD\setminus\calS=\{E\}$, and $E$ exceptional within $\calD$.
Thus, provided one has all purely min-semi-balanced systems on $N$ at disposal, the indecomposable
ones among them can be determined by this criterion.

\section{Relation of exact and totally balanced games}\label{sec.relation}
In this section we deal with the relation of the cone ${\cal E}(N)$ of exact games and the cone ${\cal T}(N)$ of
totally balanced games. In \cite[\S\,6]{KS19} a conjecture has been raised about
what are the facets of ${\cal E}(N)$, which is equivalent to the condition that
\begin{quote}
a game $\m$ over $N$ is exact iff both $\m$ and its anti-dual $\am$ are totally balanced.
\end{quote}
We give a counter-example to the conjecture in case $|N|=6$. On the other hand, we show that every
originally conjectured inequality from \cite[\S\,6]{KS19} is indeed facet-defining for
$\calE (N)$ whenever $|N|\geq 3$.

\subsection{Counter-example to a former conjecture}\label{ssec.counter-example}
Here we present a counter-example to the conjecture from \cite[\S\,6]{KS19}.
Despite we found our counter-example for $|N|=6$ computationally, by the method described in later
Remark~\ref{rem.procedure}, the reader need not repeat those computations to check its validity.
The values for coalitions in our counter-example $\m$ are given in Table~\ref{tab.counter-example}.
\medskip

\begin{table}[h]
$$
\begin{tabular}{|c||c|c|c|c|c|c|c|c|c|c|c|c|} \hline \hline
Coalition & $\emptyset$ & $a$ & $b$ & $c$ & $d$ & $e$ & $f$ & $ab$\\ \hline
Value & $0$ & $0$ & $0$ & $0$ & $0$ & $0$ & $0$ & $0$ \\
{\em tight vector}  & [1] & [7] & [1] & [2] & [1] & [11] & [1] & [7] \\
\hline \hline
Coalition & $ac$ & $ad$ & $ae$ & $af$ & $bc$ & $bd$ & $be$ & $bf$ \\ \hline
Value  & 0 & 0 & 0 & 0 & 0 & 0 & 4 & 0 \\
{\em tight vector}  & [15] & [7] & [11] & [12] & [17] & [1] & [1] & [1] \\
\hline \hline
Coalition & $cd$ & $ce$ & $cf$ & $de$ & $df$ & $ef$ & $abc$ & $abd$ \\ \hline
Value & 0 & 4 & 8 & 0 & 0 & 0 & 0 & 0 \\

{\em tight vector}
 & [2] & [18] & [12] & [11] & [1] & [12] & [17] & [7] \\
\hline  \hline
Coalition & $abe$ & $abf$ & $acd$ & $ace$ & $acf$ & $ade$ & $adf$ & $aef$ \\ \hline
Value & 4 & 0 & 0 & 8 & 8 & 0 & 0 & 0 \\
{\em tight vector}
& [8] & [16] & [15] & [2] & [12] & [11] & [13] & [12] \\
\hline  \hline
Coalition & $bcd$ & $bce$ & $bcf$ & $bde$ & $bdf$ & $bef$ & $cde$ & $cdf$ \\ \hline
Value  & 0 & 4 & 12 & 4 & 0 & 4 & 6 & 8 \\
{\em tight vector}
& [17] & [18] & [1] & [1] & [1] & [1] & [2] & [13] \\
\hline  \hline
Coalition & $cef$ & $def$ & $abcd$ & $abce$ & $abcf$ & $abde$ & $abdf$ & $abef$\\ \hline
Value & 8 & 0 & 0 & 8 & 12 & 4 & 4 & 4 \\
{\em tight vector}
 & [12] & [14] & [17] & [5] & [7] & [8] & [1] & [12] \\
 \hline  \hline
Coalition & $acde$ & $acdf$ & $acef$ & $adef$ & $bcde$ & $bcdf$ & $bcef$ & $bdef$\\ \hline
Value & 8 & 8 & 8 & 0 & 8 & 12 & 12 & 4 \\
{\em tight vector}
& [2] & [13] & [12] & [14] & [2] & [1] & [12] & [1] \\
\hline  \hline
Coalition & $cdef$ & $abcde$ & $abcdf$ & $abcef$ & $abdef$ & $acdef$ & $bcdef$ & $abcdef$\\ \hline
Value & 16 & 8 & 12 & 12 & 4 & 16 & 16 & 20 \\
{\em tight vector}
& [1] & [9] & [7] & [12] & [14] & [11] & [1] & [1] \\
\hline \hline
\end{tabular}
$$
\caption{Our counter-example $\m$ over $N=\{a,b,c,d,e,f\}$.}\label{tab.counter-example}
\end{table}

To show that $\m\in{\cal T}(N)$ for the game $\m$ from Table~\ref{tab.counter-example}
we provide its min-representation in Table~\ref{tab.ccore-min-repres} (see Section~\ref{ssec.game-concepts}). It is boring but the reader can verify manually that it is indeed a min-representation of $\m$. As a hint we indicate already in Table~\ref{tab.counter-example} at least one vector from Table~\ref{tab.ccore-min-repres} which is tight
for the respective coalition.
The first 17 vectors in Table~\ref{tab.ccore-min-repres}  are the vertices of the core of $\m$. None of them is tight for sets $\{c,e\}$ and $\{b,c,e\}$. It implies that there is no element in the core of $\m$ which is tight for
one of these sets. Thus, $\m\not\in\calE(N)$ (see again Section~\ref{ssec.game-concepts}).
\medskip

\begin{table}[h]
$$
\begin{tabular}{|c||c|c|c|c|c|c|} \hline
vector identifier & $a$ & $b$ & $c$ & $d$ & $e$ & $f$\\ \hline
[1] & 4 & 0 & 12 & 0 & 4 & 0 \\     \hline
[2] & 2 & 2 & 0 & 0 & 6 & 10 \\    \hline
[3] & 3 & 1 & 2 & 2 & 3 & 9 \\      \hline
[4] & 2 & 2 & 4 & 0 & 2 & 10 \\     \hline
[5] & 2 & 0 & 2 & 2 & 4 & 10 \\    \hline
[6] & 4 & 0 & 4 & 0 & 4 & 8 \\      \hline
[7] & 0 & 0 & 8 & 0 & 8 & 4 \\     \hline
[8] & 0 & 0 & 12 & 0 & 4 & 4 \\   \hline
[9] & 0 & 0 & 4 & 0 & 4 & 12 \\   \hline
[10] & 0 & 0 & 4 & 4 & 4 & 8 \\   \hline
[11] & 0 & 4 & 8 & 0 & 0 & 8 \\   \hline
[12] & 0 & 4 & 8 & 8 & 0 & 0 \\   \hline
[13] & 0 & 4 & 8 & 0 & 8 & 0 \\   \hline
[14] & 0 & 4 & 16 & 0 & 0 & 0 \\ \hline
[15] & 0 & 4 & 0 & 0 & 8 & 8 \\   \hline
[16] & 0 & 0 & 12 & 4 & 4 & 0 \\  \hline
[17] & 0 & 0 & 0 & 0 & 8 & 12 \\ \hline \hline
[18] & $20$ & $0$ & $0$ & $20$ & $4$ & $20$\\ \hline
\end{tabular}
$$
\caption{The min-representation of our counter-example $\m$.}\label{tab.ccore-min-repres}
\end{table}

The anti-dual $\am$ of our game $\m$ is given in Table~\ref{tab.anti-dual}, its min-representation in
Table~\ref{tab.anti-dual-min-repres}. The core of the anti-dual has also 17
vertices, presented in the first 17 lines of the table. None of them is tight for sets $\{a,d,f\}$ and
$\{a,b,d,f\}$; in particular, $\am\not\in\calE (N)$, which fact also follows from a former observation that
$\m\not\in\calE(N)$. A min-representation of $\am$ can be obtained by adding one additional vector. The reader
can verify manually that Table~\ref{tab.anti-dual-min-repres} indeed provides a min-representation of $\am$; for each coalition, we indicate in Table~\ref{tab.anti-dual} at least one vector from
Table~\ref{tab.anti-dual-min-repres} which is tight for it.

\begin{table}[h]
$$
\begin{tabular}{|c||c|c|c|c|c|c|c|c|c|c|c|c|} \hline \hline
Coalition & $\emptyset$ & $a$ & $b$ & $c$ & $d$ & $e$ & $f$ & $ab$\\ \hline
Value & 0 & -4 & -4 & -16 & -8 & -8 & -12 & -4 \\
{\em tight vector}
& [0] & [3] & [1] & [17] & [16] & [6] & [4] & [1] \\
\hline \hline
Coalition & $ac$ & $ad$ & $ae$ & $af$ & $bc$ & $bd$ & $be$ & $bf$ \\ \hline
Value & -16 & -8 & -8 & -12 & -20 & -12 & -12 & -12\\
{\em tight vector}
& [15] & [16] & [3] & [2] & [17] & [16] & [8] & [1] \\
\hline \hline
Coalition & $cd$ & $ce$ & $cf$ & $de$ & $df$ & $ef$ & $abc$ & $abd$ \\ \hline
Value & -16 & -16 & -16 & -8 & -12 & -20 & -20 & -12 \\
{\em tight vector}
& [13] & [5] & [1] & [6] & [4] & [9] & [17] & [16] \\
\hline  \hline
Coalition & $abe$ & $abf$ & $acd$ & $ace$ & $acf$ & $ade$ & $adf$ & $aef$ \\ \hline
Value & -12 & -14 & -16 & -20 & -16 & -8 & -16 & -20 \\
{\em tight vector}
& [8] & [2] & [13] & [15] & [1] & [3] & [18] & [9] \\
\hline  \hline
Coalition & $bcd$ & $bce$ & $bcf$ & $bde$ & $bdf$ & $bef$ & $cde$ & $cdf$ \\ \hline
Value & -20 & -20 & -20 & -12 & -12 & -20 & -20 & -16 \\
{\em tight vector}
& [16] & [14] & [1] & [8] & [1] & [8] & [13] & [1] \\
\hline  \hline
Coalition & $cef$ & $def$ & $abcd$ & $abce$ & $abcf$ & $abde$ & $abdf$ & $abef$\\ \hline
Value & -20 & -20 & -20 & -20 & -20 & -12 & -16 & -20 \\
{\em tight vector}
& [4] & [9] & [16] & [14] & [1] & [8] & [18] & [7] \\
\hline  \hline
Coalition & $acde$ & $acdf$ & $acef$ & $adef$ & $bcde$ & $bcdf$ & $bcef$ & $bdef$\\ \hline
Value & -20 & -16 & -20 & -20 & -20 & -20 & -20 & -20 \\
{\em tight vector}
& [13] & [1] & [3] & [9] & [13] & [1] & [1] & [8] \\
\hline  \hline
Coalition & $cdef$ & $abcde$ & $abcdf$ & $abcef$ & $abdef$ & $acdef$ & $bcdef$ & $abcdef$\\ \hline
Value & -20 & -20 & -20 & -20 & -20 & -20 & -20 & -20 \\
{\em tight vector}
& [4] & [13] & [1] & [1] & [7] & [3] & [1] & [1] \\
\hline \hline
\end{tabular}
$$
\caption{The anti-dual $\am$ of our counter-example $\m$.}\label{tab.anti-dual}
\end{table}

\begin{table}[h]
$$
\begin{tabular}{|c||c|c|c|c|c|c|} \hline
vector identifier & $a$ & $b$ & $c$ & $d$ & $e$ & $f$\\ \hline
[1] & 0 & -4 & -8 & 0 & 0 & -8 \\ \hline
[2] & -2 & -2 & -4 & 0 & -2 & -10 \\ \hline
[3] & -4 & 0 & -4 & 0 & -4 & -8 \\ \hline
[4] & 0 & 0 & -4 & 0 & -4 & -12 \\ \hline
[5] & 0 & -0 & -12 & 0 & -4 & -4 \\ \hline
[6] & 0 & -0 & -8 & 0 & -8 & -4 \\ \hline
[7] & -2 & -2 & 0 & 0 & -6 & -10 \\ \hline
[8] & 0 & -4 & 0 & 0 & -8 & -8 \\ \hline
[9] & 0 & 0 & 0 & 0 & -8 & -12 \\ \hline
[10] & -2 & 0 & -2 & -2 & -4 & -10 \\ \hline
[11] & 0 & 0 & -4 & -4 & -4 & -8 \\ \hline
[12] &  -3 & -1 & -2 & -2 & -3 & -9 \\ \hline
[13] & 0 & -0 & -12 & -4 & -4 & 0 \\ \hline
[14] & 0 & -4 & -8 & 0 & -8 & 0 \\ \hline
[15] & -4 & 0 & -12 & 0 & -4 & 0 \\ \hline
[16] & 0 & -4 & -8 & -8 & 0 & 0 \\ \hline
[17] & 0 & -4 & -16 & 0 & 0 & 0 \\ \hline  \hline
[18] & -4 & 0 & 0 & -4 & 0 & -8\\ \hline
\end{tabular}
$$
\caption{The min-representation of the anti-dual $\am$.}\label{tab.anti-dual-min-repres}
\end{table}

\begin{remark}\label{rem.procedure}\rm
This is to describe the way we found our counterexample. Our method was based on the characterization
of $\calE (N)$ from Section~\ref{sec.char-exact}. We have succeeded to compute the extreme rays of all the cones $\Theta^{N}_{D}$, $\emptyset\neq D\subseteq N$, in case $|N|=6$. Thus, we got a finite set of linear inequalities characterizing the cone $\calE (N)$ in this case, although a pretty big one. Additionally, on basis of the results from \cite{SKV19}, we were able to get a complete list $\calL$ of coefficient vectors for the conjectured facet-defining inequalities. Checking of the validity of the conjecture in case $|N|=6$ was, therefore, reduced to checking of whether, for every $\emptyset\neq D\subseteq N$, every generator of an extreme ray of $\Theta^{N}_{D}$ is in the conic hull of $\calL$.

This has appeared not to be the case: we found an extreme ray of $\Theta^{N}_{D}$ for $|D|=2$
which is not in the conic hull of ${\cal L}$. Specifically, it was an element of $\hat{\theta}\in \Theta^{N}_{D}$
for $N=\{a,b,c,d,e,f\}$ and $D=\{c,e\}$ defined as follows (we write $ce$ instead of $\{c,e\}$ here):
\begin{eqnarray*}
\lefteqn{\hspace*{-3mm}\hat{\theta}(\emptyset)=+1, \quad \hat{\theta}(ce)=+4,\quad \hat{\theta}(abcdef)=+3,}\\[0.2ex]
&& \hat{\theta}(be)=-1, \quad \hat{\theta}(ace)=-3, \quad \hat{\theta}(bcf)=-1, \quad \hat{\theta}(bcde)=-1, \quad \hat{\theta}(cdef)=-2\,.
\end{eqnarray*}
On basis of that objective vector $\hat{\theta}$ we found a game $\m$ that satisfies $\langle \hat{\theta},\m\rangle< 0$
while $\langle \theta^{\prime},\m\rangle\geq 0$ for any $\theta^{\prime}\in {\cal L}$,
which is just the game $\m$ presented in Table~\ref{tab.counter-example}.
Note, however, that the vector $\hat{\theta}$ has appeared {\sf not} to yield a facet-defining
inequality for $\m\in {\cal E}(N)$.

Note in this context that an alternative idea of computing the extreme rays of the cone $\{\m\in {\cal G}(N)\,:\ \langle \theta,\m\rangle\geq 0~~\mbox{for $\theta\in\calL$}\}$ has appeared to be computationally infeasible. This is because the number of the extreme rays of\/ ${\cal E}(N)$ grows very rapidly with $|N|$ and we were not able to compute them even in case $|N|=5$.
\end{remark}

\subsection{Facets shared with the cone of totally balanced games}\label{ssec.sub-balanced}
This is to relate our new concept of an indecomposable min-semi-balanced system
to earlier concepts and results from \cite{KS19}, where facet-defining inequalities
for the cone $\mathcal{T}(N)$ of totally balanced games were characterized.
The following is a simplified equivalent definition of a central concept from \cite[\S\,4]{KS19};
the equivalence of the original definition and the later simplified version of it was shown in \cite[\S\,2.3]{SKV19}.

\begin{defin}[reducible and irreducible balanced set system]~\label{def.reducible}\rm \\
A min-balanced set system $\calB$ on a finite set $M$, $|M|\geq 2$,
will be called {\em reducible\/} if there exists a set $\emptyset\neq E\subset M$ such that $\chi_{E}$ is a conic combination of $\{\, \chi_{S}\,:\ S\in\calB \,~\&~\, S\subset E\,\}$.\\
A min-balanced set system $\calB$ on $M$ which is not reducible is called {\em irreducible}.
\end{defin}

It was shown in \cite[Lemma~2.1]{KS19} that a balanced system $\calB$ on $M$ is minimal iff
vectors $\{\chi_{S}\,:\ S\in\calB\}$ are linearly independent. In particular, for a min-balanced set system $\calB$ on $M$, there is a unique linear combination $\sum_{S\in\calB} \gamma_{S}\cdot\chi_{S}$ yielding $\chi_{M}$ and this unique combination has all coefficients strictly positive: $\gamma_{S}>0$ for $S\in\calB$. Standard interpretation of a min-balanced system $\calB$ is then in terms of the assigned linear inequality
$$
\m(M)\geq \sum_{S\in\calB}\, \gamma_{S}\cdot\m(S)\,,
$$
which can be interpreted as an inequality for any game $\m$ on a superset $N$ of $M$. Recall that the main result from \cite[Theorem~5.1]{KS19} says that the facet-defining inequalities for the cone $\mathcal{T}(N)$ of totally balanced games on $N$, $|N|\geq 2$, are just the inequalities assigned to irreducible min-balanced systems on subsets $M\subseteq N$, $|M|\geq 2$.

Recall from Lemma~\ref{lem.clasify} that every min-balanced set system $\calB$ on {\sf proper\/} subset $M\subset N$, $|M|\geq 2$, corresponds to a min-semi-balanced system $\calS:=\calB\cup\{M\}$. The unique affine combination
$\sum_{S\in\calS} \lambda_{S}\cdot\chi_{S}$ yielding a constant vector $\0\in {\dv R}^{N}$
is semi-conic with $\lambda_{M}<0$ and can be written as $\sum_{S\in\calB} \lambda_{S}\cdot\chi_{S}=(-\lambda_{M})\cdot\chi_{M}$.
Therefore, the induced inequality \eqref{eq.theta-ineq} for $\m\in {\dv R}^{\caP}$, $\m(\emptyset)=0$, is equivalent to (= is a positive multiple of) the standard inequality assigned to the min-balanced system $\calB$.

The next result says that every irreducible min-balanced system on a proper subset $M$ of $N$ yields
an indecomposable min-semi-balanced  system on $N$. Thus, by Theorem~\ref{thm.main-result},
it implies that the respective inequalities are facet-defining for the exact cone ${\cal E}(N)$.
\smallskip

\begin{lem}\label{lem.sub-balanced}~\rm
Given a min-semi-balanced set system of the form $\calS:=\calB\cup\{M\}$ on $N$, $|N|\geq 3$,
where $\calB$ is a min-balanced set system on $M\subset N$, $|M|\geq 2$, the next two conditions
are equivalent:
\begin{description}
\item[(i)] ~$\calB$ is irreducible,
\item[(ii)] $\calS$ is indecomposable.
\end{description}
\end{lem}

\begin{proof}
We prove the equivalence of negations of those conditions.
\smallskip

To show $\neg \mbox{\bf (i)}\Rightarrow \neg \mbox{\bf (ii)}$ assume that $\calB$ is a reducible min-balanced system on $M\subset N$. Consider $\emptyset\neq E\subset M$ and a conic combination
$\sum_{S\in\calB:S\subset E} \lambda_{S}\cdot\chi_{S}=\chi_{E}$ 
and re-write that as a semi-conic combination $\sum_{S\in\calB:S\subset E} \lambda_{S}\cdot\chi_{S}+(-1)\cdot\chi_{E}=\0$ in the space ${\dv R}^{N}$. Then put $\lambda_{S}:=0$
for remaining $S\in\calS=\calB\cup\{M\}$, and, because $\0$ is a constant vector in ${\dv R}^{N}$, conclude that $E$
is exceptional within $\calS\cup\{E\}$. This means that $\calS$ has a  decomposition.
\smallskip

To verify $\neg \mbox{\bf (ii)}\Rightarrow \neg \mbox{\bf (i)}$ assume that $\calS=\calB\cup\{M\}$
has a  decomposition, that is,
a linear combination $\sum_{S\in\calS} \lambda_{S}\cdot\chi_{S}+\lambda_{E}\cdot\chi_{E}=r\cdot\chi_{N}$
yielding a constant vector in ${\dv R}^{N}$ exists with $\lambda_{S}\geq 0$ for $S\in\calS$
and with $\lambda_{E}<0$ for some $\emptyset\neq E\subset N$, $E\not\in\calS$.
Since it is a semi-conic combination, by Lemma~\ref{lem.semi-comb},
one has $r\geq 0$. In case $E\setminus M\neq\emptyset$ one can choose $i\in E\setminus M$ and
its substitution to the equality gives a contradictory conclusion $r=\lambda_{E}<0$.
Thus, as $E\not \in\calS=\calB\cup\{M\}$, $E\subset M\subset N$ and the choice $j\in N\setminus M$ and its substitution to the equality gives $r=0$, that is, $\sum_{S\in\calB\,\cup\{M\}} \lambda_{S}\cdot\chi_{S}+\lambda_{E}\cdot\chi_{E}=\sf0$.
Hence, the non-negativity of the coefficients implies that, for any $S\in\calS$
such that $S\setminus E\neq\emptyset$ (including $S=M$) one necessarily has $\lambda_{S}=0$.
Thus, one can write it in the form $(-\lambda_{E})\cdot\chi_{E}=\sum_{S\in\calB: S\subset E} \lambda_{S}\cdot\chi_{S}$,
which easily implies that $\calB$ is reducible.
\end{proof}

Hence, by Lemma~\ref{lem.sub-balanced}, every facet-defining inequality for the totally
balanced cone ${\cal T}(N)$ that corresponds to a strict subset $M\subset N$ is also
facet-defining for the exact cone ${\cal E}(N)$. These are facet-defining inequalities
for both cones. Note that, by Lemma~\ref{lem.conjug}, also conjugate inequalities to these
inequalities are facet-defining for ${\cal E}(N)$. These two classes of inequalities induced
by irreducible min-balanced systems on $M\subset N$ were originally conjectured in \cite[\S\,6]{KS19}
to be all facet-defining inequalities for ${\cal E}(N)$.

\section{Conclusions}\label{sec.coclusion}
The main achievement in this paper is the observation that a linear inequality
for games over $N$, $|N|\geq 3$, is facet-defining for the exact cone $\calE(N)$
iff it corresponds to (uniquely determined) {\em indecomposable min-semi-balanced\/}
set system on $N$ (Theorem~\ref{thm.main-result}). At first we got
these inequalities in case $|N|=6$ by computation and later we confirmed the conjecture
that the correspondence holds in general. Because of the computation, we know what are the numbers of
facets of $\calE(N)$ in cases $2\leq |N|\leq 6$; they are shown in Table~\ref{tab.facet-numbers}.

\begin{table}[h]
$$
\begin{tabular}{|l|c|c|c|c|c|} \hline
Number of players & $n=2$ & $n=3$ & $n=4$ & $n=5$ & $n=6$\\ \hline
Number of facets & $1$ & $6$ & $44$ & $280$ & $7006$\\
Number of its permutational types & $1$ & $2$ & $6$ & $16$ & $53$\\ \hline
\end{tabular}
$$
\caption{Numbers of facets of $\calE(N)$ and of its permutational types for $n=|N|\leq 6$.}\label{tab.facet-numbers}
\end{table}

Note that the case $|N|=2$ is special in some sense: then the cone of exact games coincides with the cone of balanced games and
the only facet-defining inequality corresponds to the only min-balanced set system on $N$ with $|N|=2$.
In case $|N|=3$ the cone of exact games coincides with the cone of supermodular (= convex) games, while
for $|N|\geq 4$ these two cones already differ.
For $|N|\leq 5$ every facet-defining inequality for $\calE (N)$ corresponds to a min-semi-balanced system
$\calS$ satisfying either $\bigcup\calS\subset N$ or $\bigcap\calS\neq\emptyset$; however, this is not the case
in case $|N|=6$.
\smallskip

Min-semi-balanced systems and their induced inequalities can synoptically be described
by means of special pictures/diagrams (see Section~\ref{ssec.picto-grams}). A catalogue of
permutational types of {\em irreducible min-semi-balanced\/} set systems over $N$, $3\leq |N|\leq 6$,
obtained as a result of our computation, is available:
$$
\mbox{\sf http://gogo.utia.cas.cz/indecomposable-min-semi-balanced-catalogue/}\,;
$$
thus, it implicitly provides an overview of all facet-defining inequalities in these cases.
\smallskip

Recall that the min-semi-balanced systems break into four basic classes (Section~\ref{ssec.classification})
and the pictorial representatives reflect this classification. They also reflect the fact that
(indecomposable) min-semi-balanced systems are closed under {\em complementarity\/} transform (Section~\ref{ssec.complementarity}).
Further relevant observation is that min-semi-balanced systems on $N$ can be obtained on basis of min-balanced systems on $N$
(Section~\ref{ssec.get-purely}), which suggests that one can possibly get all such systems for $|N|\geq 7$.
\smallskip

The second main result in this paper is an example of a game $\m$ over $N$, $|N|=6$, such that
both $\m$ and its anti-dual $\am$ are totally balanced while $\m$ is {\sf not\/} exact (Section~\ref{ssec.counter-example}).
%

\subsubsection*{Acknowledgements}
This research has been supported by the grant GA\v{C}R n.\ 19-04579S.
We are grateful to our colleague Tom\'{a}\v{s} Kroupa for a consultation with him
which helped us to overcome an obstacle on our way to compute facets of the exact cone in case $|N|=6$.

\appendix

\section{Proof of Lemma~\ref{lem.min-semi-bal}}\label{app.lem-2}
For reader's convenience we recall the result.
\bigskip

\noindent{\bfseries Lemma~\ref{lem.min-semi-bal}:} ~~
Given $|N|\geq 2$, let $\emptyset\neq\calS\subseteq {\caP}\setminus \{\emptyset,N\}$ be a non-trivial set system on $N$. Then the following conditions on $\calS$ are equivalent:
\begin{description}
\item[(a)] $\calS$ is a minimal set system such that there is a constant vector in ${\dv R}^{N}$ which can be
written as a non-zero semi-conic combination of vectors $\{ \chi_{S}\,:\ S\in\calS\}$,
\item[(b)] $\calS$ is a minimal semi-balanced set system on $N$,
\item[(c)] $\calS$ is semi-balanced on $N$, the vectors $\{ \chi_{S}\,:\ S\in\calS\}$ are affinely independent
and in case $\bigcup\calS=N$ even linearly independent,
\item[(d)] there is only one affine combination of vectors $\{ \chi_{S}\,:\ S\in\calS\}$ yielding
a constant vector in ${\dv R}^{N}$ and this unique combination is semi-conic and has all coefficients non-zero,
\item[(e)] there exists unique affine semi-conic combination of vectors $\{ \chi_{S}\,:\ S\in\calS\}$
which is a constant vector in ${\dv R}^{N}$ and this unique combination has all coefficients non-zero.
\end{description}

\begin{proof}
To show {\bf (a)}\,$\Rightarrow${\bf (b)} assume a non-zero semi-conic combination $\sum_{S\in\calS} \lambda_{S}\cdot\chi_{S}$ which is a constant vector in ${\dv R}^{N}$ and put $\calS^{\prime}:=\{ S\in\calS\, :\ \lambda_{S}\neq 0\}$.
Because of minimality of $\calS$ in {\bf (a)} one has $\calS^{\prime}=\calS$, which implies that $\calS$ is semi-balanced. The rest is evident.
\smallskip

To show {\bf (b)}\,$\Rightarrow${\bf (a)} it is enough to verify the minimality of $\calS$ in {\bf (a)}.
Assume for a contradiction that a set system $\calC\subset\calS$ exists with a non-zero semi-conic combination $\sum_{S\in\calC} \mu_{S}\cdot\chi_{S}$ yielding a constant vector in ${\dv R}^{N}$ and put $\calC^{\prime}:=\{ S\in\calC\, :\ \mu_{S}\neq 0\}$. Then $\calC^{\prime}$ is semi-balanced on $N$ and $\calC^{\prime}\subset\calS$ contradicts the minimality of
$\calS$ in {\bf (b)}.
\smallskip

To show {\bf (b)}\,$\Rightarrow${\bf (c)} let us fix a semi-conic combination $\sum_{S\in\calS} \lambda_{S}\cdot\chi_{S}$ yielding a constant vector in ${\dv R}^{N}$ with $\lambda_{S}\neq 0$ for $S\in\calS$. We first consider the case when it is {\sf not} a conic combination,
that is, there is a unique set $T\in\calS$ with $\lambda_{T}<0$; hence, $\lambda_{S}>0$ for $S\in\calS\setminus\{T\}$.

We then verify that $\{ \chi_{S}\,:\ S\in\calS\}$ are affinely independent. Assume for a contradiction that a non-zero linear combination $\sum_{S\in\calS} \mu_{S}\cdot\chi_{S}=\0$ with $\sum_{S\in\calS} \mu_{S}=0$ exists.  One can assume without loss of generality $\mu_{T}\leq 0$ for otherwise one can multiply the linear combination by $(-1)$. Moreover, $\sum_{S\in\calS} \mu_{S}=0$
implies the existence of $L\in\calS$ with $\mu_{L}<0$. In fact, there are at least two such sets:
otherwise $\sum_{S\in\calS} \mu_{S}\cdot\chi_{S}=\0$ is a semi-conic combination and, by Lemma~\ref{lem.semi-comb},
one has $\sum_{S\in\calS} \mu_{S}>0$ which contradicts the assumption. In particular, there is
$L\in\calS\setminus\{T\}$ with $\mu_{L}<0$.
For any $\varepsilon\geq 0$ and $S\in\calS$ we put $\lambda^{\varepsilon}_{S}:= \lambda_{S}+\varepsilon\cdot\mu_{S}$ and observe that $\sum_{S\in\calS} \lambda^{\varepsilon}_{S}\cdot\chi_{S}$ yields the same constant vector,
$\lambda^{\varepsilon}_{T}<0$, while $\lambda^{\varepsilon}_{S}>0$ for $S\in\calS\setminus\{T\}$ and
sufficiently small $\varepsilon$. Since $\lambda^{\varepsilon}_{L}$ tends to $-\infty$ with increasing $\varepsilon$
there is a maximal $\varepsilon^{\ast}>0$ such that $\lambda^{\varepsilon^{\ast}}_{S}\geq 0$ for all $S\in\calS\setminus\{T\}$. There must be $K\in\calS\setminus\{T\}$ with $\lambda^{\varepsilon^{\ast}}_{K}=0$;
then the system $\calC :=\{S\in\calS\,: \ \lambda^{\varepsilon^{\ast}}_{S}\neq 0\}\not\ni K$ is semi-balanced
on $N$, which contradicts the minimality of $\calS$.

The second step is to show that if $\bigcup\calS=N$ then $\{ \chi_{S}\,:\ S\in\calS\}$ are linearly independent. Note that the case of a conic combination $\sum_{S\in\calS} \lambda_{S}\cdot\chi_{S}$ is involved: this follows from Lemma~\ref{lem.semi-comb} because then $r>0$ enforces $\bigcup\calS=N$. In the
sequel, let $T$ denote a contingent set $T\in\calS$ with $\lambda_{T}<0$, which, however, need not exist. Recall that $\lambda_{S}>0$ for $S\in\calS\setminus\{T\}$.

Assume for a contradiction that a non-zero linear combination $\sum_{S\in\calS} \mu_{S}\cdot\chi_{S}=\0$ exists.
We then show that there exists such a linear combination which, additionally, satisfies $\mu_{T}\leq 0$ if $T$ exists and $\mu_{L}<0$ for at least one $L\in \calS\setminus\{T\}$. This is easy in case $T$ does not exist or $\mu_{T}=0$ because possible multiplication of the linear combination by $(-1)$ reaches the goal. In case $\mu_{T}\neq 0$ possible multiplication
ensures $\mu_{T}<0$. To show then the existence of $L\in \calS\setminus\{T\}$ with $\mu_{L}<0$ assume for
a contradiction the opposite, which means that $\sum_{S\in\calS} \mu_{S}\cdot\chi_{S}=\0$ is a semi-conic
combination and $\calC := \{S\in\calS\,:\ \mu_{S}\neq 0\}$ is a semi-balanced system. The minimality of $\calS$
then implies $\calC=\calS$ and, thus, $\mu_{S}>0$ for $S\in \calS\setminus\{T\}$. Since $T\subset N$ there
exists $i\in N\setminus T$ and $\bigcup\calS=N$ implies the existence of $K\in\calS$ with $i\in K$.
Since [\,$i\in S\in\calS ~\Rightarrow~ \mu_{S}>0$\,], this leads
to a contradictory conclusion $0=\sum_{S\in \calS} \mu_{S}\cdot\chi_{S}(i)=
\sum_{S\in \calS: i\in S} \mu_{S}\geq \mu_{K}>0$.

Finally, having a linear combination $\sum_{S\in\calS} \mu_{S}\cdot\chi_{S}=\0$ with $\mu_{T}\leq 0$ and
$\mu_{L}<0$ for some $L\in \calS\setminus\{T\}$, one can repeat the construction used in the previous case
(of affine independence) to get a contradiction with the minimality of $\calS$.
\smallskip

To show {\bf (c)}\,$\Rightarrow${\bf (d)} use Lemma~\ref{lem.semi-comb} to
obtain an affine semi-conic combination $\sum_{S\in\calS} \lambda_{S}\cdot\chi_{S}$
with all coefficients non-zero yielding a constant vector $\rho=[r,\ldots,r]\in {\dv R}^{N}$ with $r\in [0,1]$.
Let us fix this affine combination. Assume that $\sum_{S\in\calS} \sigma_{S}\cdot\chi_{S}$ is an affine combination yielding a constant vector $\varsigma=[s,\ldots,s]\in {\dv R}^{N}$. It is enough to show that these two combinations coincide. To this end we distinguish two cases.

In case $\bigcup\calS\subset N$ choose $i\in N\setminus\bigcup\calS$ and have $\sum_{S\in\calS} \lambda_{S}\cdot\chi_{S}(i)=0=\sum_{S\in\calS} \sigma_{S}\cdot\chi_{S}(i)$.
This implies that both $\rho=\0$ and $\varsigma=\0$.
By subtracting we get $\sum_{S\in\calS} (\lambda_{S}-\sigma_{S})\cdot\chi_{S}=\0$ with $\sum_{S\in\calS} (\lambda_{S}-\sigma_{S}) =0$ and  by affine independence $\sigma_{S}=\lambda_{S}$ for all $S\in\calS$.

In case $\bigcup\calS=N$ the linear independence of vectors $\{ \chi_{S}\,:\ S\in\calS\}$ implies
both $r\neq 0$ and $s\neq 0$ because solely their zero linear combination yields the vector $\0$.
Hence, we have $\sum_{S\in\calS} (r^{-1}\cdot\lambda_{S})\cdot\chi_{S}=\chi_{N}$ and
$\sum_{S\in\calS} (s^{-1}\cdot\sigma_{S})\cdot\chi_{S}=\chi_{N}$. By
subtracting we get $\sum_{S\in\calS} (r^{-1}\cdot\lambda_{S}-s^{-1}\cdot\sigma_{S})\cdot\chi_{S}=\0$
and by linear independence $r^{-1}\cdot\lambda_{S}=s^{-1}\cdot\sigma_{S}$ for all $S\in\calS$.
Thus, $\lambda_{S}=r\cdot s^{-1}\cdot\sigma_{S}$ for $S\in\calS$ and, because both combinations
are affine, by summing over $S\in\calS$ one derives $1=r\cdot s^{-1}$. Hence, $s=r$ and $\sigma_{S}=\lambda_{S}$ for all $S\in\calS$.
\smallskip

The implication {\bf (d)}\,$\Rightarrow${\bf (e)} is evident.
\smallskip

To show {\bf (e)}\,$\Rightarrow${\bf (b)} consider the affine semi-conic combination
$\sum_{S\in\calS} \lambda_{S}\cdot\chi_{S}$ with all coefficients non-zero  yielding a constant vector in ${\dv R}^{N}$.
Its existence implies that $\calS$ is semi-balanced. To verify the minimality of $\calS$ assume
for a contradiction that $\calC\subset\calS$ exists which is semi-balanced on $N$.
By Lemma~\ref{lem.semi-comb} applied to $\calC$ there exists an affine semi-conic combination
$\sum_{S\in\calC} \sigma_{S}\cdot\chi_{S}$ yielding a constant vector in ${\dv R}^{N}$. We extend it
by putting $\sigma_{S}=0$ for $S\in\calS\setminus\calC$. Thus, we get
two different affine semi-conic combinations of $\{ \chi_{S}\,:\ S\in\calS\}$
yielding a constant vector in ${\dv R}^{N}$, which contradicts the assumption.
\end{proof}

\section{Proof of Lemma~\ref{lem.theta}}\label{app.lem-theta}
For reader's convenience we recall what is claimed; see Definition~\ref{def.aux-cone} for notation.
\bigskip

\noindent{\bfseries Lemma~\ref{lem.theta}:} ~~
Given $|N|\geq 2$, every set $\tilde{\Theta}^{N}_{D}$, where $\emptyset\neq D\subseteq N$, is a bounded polyhedron.
Every vector $\theta\in\Theta^{N}_{D}$ satisfies both $\theta(N)\geq 0$ and $\theta(\emptyset)\geq 0$ and
every {\sf non-zero\/} vector $\theta\in\Theta^{N}_{D}$ satisfies $\theta(N)+\theta(\emptyset)>0$.
Given $\m\in {\dv R}^{\caP}$ with $\m(\emptyset)=0$, one has
$$
\hspace*{30mm}\m\in {\cal E}(N) ~~\Leftrightarrow~~ [\,\forall\, \theta\in \bigcup_{\emptyset\neq D\subseteq N}\tilde{\Theta}^{N}_{D}\qquad \langle \theta,\m\rangle\geq 0\,]\,.
\hspace*{32mm}\hfill \eqref{eq.exact-char}
$$

\begin{proof}
Note that $\Theta^{N}_{N}\subseteq \Theta^{N}_{D}$ if $D\subset N$; thus, assume without loss of generality $\emptyset\neq D\subset N$.
To show $\theta(N)\geq 0$ for $\theta\in \Theta^{N}_{D}$ choose $i\in N\setminus D$ and write
$\theta(N)=-\sum_{L\subset N:i\in L} \theta(L)\geq 0$. Note that,
for any $j\in N$, $\sum_{S\subseteq N\setminus\{j\}} \theta(S)=
\sum_{S\subseteq N} \theta(S)-\sum_{L\subseteq N:j\in L} \theta(L) = 0-0 =0$.
Hence, to show $\theta (\emptyset)\geq 0$ for $\theta\in \Theta^{N}_{D}$ take $j\in D$
and write $\theta(\emptyset)=-\sum_{\emptyset\neq S\subseteq N\setminus\{j\}} \theta(S)\geq 0$.

Thus, we have observed that every $\theta\in \Theta^{N}_{D}$ satisfies both $\theta(N)\geq 0$ and $\theta(\emptyset )\geq 0$.
In particular, if $\theta\in \tilde{\Theta}^{N}_{D}$ then $\theta(N)+\theta(\emptyset)=1$ gives both $0\leq \theta(\emptyset)\leq 1$ and $0\leq \theta(N)\leq 1$.

The next step is to show that, for $\theta\in \tilde{\Theta}^{N}_{D}$, if $S\subset N$ and $S\setminus D\neq\emptyset$ then $0\geq \theta(S)\geq -1$. Indeed, because of $S\not\in\{\emptyset,D,N\}$, the choice of $i\in S\setminus D$ gives
$$
0\geq \theta (S)\geq \sum_{L\subset N:i\in L} \theta (L)=-\theta(N)+ \underbrace{\sum_{L\subseteq N:i\in L} \theta (L)}_{=0} =-\theta(N)\geq -1\,.
$$
To observe that $0\geq \theta(S)\geq -1$ for $\theta\in \tilde{\Theta}^{N}_{D}$ whenever $\emptyset\neq S\subset N$ and $D\setminus S\neq\emptyset$ introduce a vector $\theta^{\refl}\in{\dv R}^{\caP}$ by $\theta^{\refl}(L):=\theta (N\setminus L)$ for $L\subseteq N$. It is easy to observe that $\theta^{\refl}\in\tilde{\Theta}^{N}_{N\setminus D}$: to this end write for any $i\in N$
\begin{eqnarray*}
\sum_{L\subseteq N:i\in L} \theta^{\refl}(L)=\sum_{L\subseteq N:i\in L} \theta(N\setminus L)=\sum_{S\subseteq N:i\not\in S} \theta(S) =\sum_{S\subseteq N} \theta(S)-\sum_{S\subseteq N:i\in S} \theta(S)=0-0=0.
\end{eqnarray*}
Thus, because of $\emptyset\neq D\setminus S= (N\setminus S)\setminus(N\setminus D)$, one has $0\geq \theta^{\refl}(N\setminus S)\geq -1$ by the previous observation applied to $\theta^{\refl}$, which, however, means $0\geq \theta(S)\geq -1$.

Altogether, we have $0\geq \theta(S)\geq -1$ for $\theta\in \tilde{\Theta}^{N}_{D}$ and $S\in\caP\setminus \{\emptyset,D,N\}$,
which implies  $0\geq \sum_{S: S\not\in \{\emptyset,D,N\}} \theta(S)\geq3-2^{|N|}$. Taking into consideration that
$$
\sum_{S: S\not\in \{\emptyset,D,N\}} \theta(S)=
-\theta(\emptyset)-\theta(D)-\theta(N)+\underbrace{\sum_{L\subseteq N} \theta (L)}_{=0}= -\theta(N)-\theta(\emptyset)-\theta(D)=-1-\theta(D)
$$
one gets $2^{|N|}-4\geq \theta(D)\geq -1$. In particular,
$\tilde{\Theta}^{N}_{D}$ is bounded
and so is $\tilde{\Theta}^{N}_{N}$.
\smallskip

The fact that, for any $\emptyset\neq D\subseteq N$, every non-zero vector $\theta\in\Theta^{N}_{D}$ satisfies $\theta(N)+\theta(\emptyset)>0$
follows directly from \cite[Lemma~5.1]{KS19}. In particular, every non-zero $\theta\in \Theta^{N}_{D}$ is a positive multiple of a vector
$\tilde{\theta}\in \tilde{\Theta}^{N}_{D}$.
Further important fact, which follows from \cite[Lemma~5.3]{KS19}, is that a game $\m$ is exact, that is, $\m\in {\cal E}(N)$, iff $[\,\forall\:\emptyset\neq D\subseteq N ~~\forall\, \theta\in \Theta^{N}_{D}\quad \langle \theta,\m\rangle\geq 0\,]$. Putting these two observations together gives \eqref{eq.exact-char}.
\end{proof}


\section{Proof of Lemma~\ref{lem.extr-theta}}\label{app.lem-extr-theta}
For reader's convenience we recall what is claimed; see Definition~\ref{def.aux-cone} for notation.
\bigskip

\noindent{\bfseries Lemma~\ref{lem.extr-theta}:} ~~
Given $|N|\geq 2$ and $\emptyset\neq D\subseteq N$, every vertex of $\tilde{\Theta}^{N}_{D}$ has either the form
$\theta_{\calB}$, where $\calB$ is a min-balanced set system on $N$, or the form $\theta_{\calS}$, where $\calS$ is a min-semi-balanced system on $N$ having $D$ as the exceptional set.\\
Conversely, in case $\emptyset\neq D\subset N$, every vector $\theta_{\calS}$, where
$\calS$ is a min-semi-balanced system on $N$ having $D$ as the exceptional set, is a vertex of
$\tilde{\Theta}^{N}_{D}$: $\theta_{\calS}\in \ext (\tilde{\Theta}^{N}_{D})$.
\medskip

\begin{proof}
Given a vertex $\theta\in \ext (\tilde{\Theta}^{N}_{D})$, we first observe that there exists a min-semi-balanced system $\calS$ on $N$ such that $\theta=\theta_{\calS}$. To this end we put
$$
\calS ~:=~ \{S\subseteq N\,:\ \emptyset\neq S\subset N ~~\&~~ \theta(S)\neq 0\,\}\qquad \mbox{and~~ $\lambda_{S}=-\theta(S)$\quad for $S\in\calS$.}
$$
Thus, $\sum_{S\in\calS} \lambda_{S}\cdot\chi_{S}$ yields a constant vector $\theta(N)\cdot\chi_{N}$ and $\calS$ is semi-balanced on $N$.
To evidence that $\calS$ is min-semi-balanced
we use the condition {\bf (c)} in Lemma~\ref{lem.min-semi-bal}. To verify affine independence of the vectors $\{ \chi_{S}\,:\ S\in\calS\}$ assume
for a contradiction that there is a non-zero linear combination $\sum_{S\in\calS} \sigma_{S}\cdot\chi_{S}=\0$ with $\sum_{S\in\calS} \sigma_{S}=0$.
Then we put $\sigma_{L}:=0$ for remaining $L\subseteq N$ and $\theta^{\varepsilon}(S):=\theta(S)+\varepsilon\cdot\sigma_{S}$ for any $S\subseteq N$
and $\varepsilon\in {\dv R}$.
Then one has $\theta^{\varepsilon}\in \tilde{\Theta}^{N}_{D}$ whenever $|\varepsilon|$ is small; thus, the relation
$\theta=\frac{1}{2}\cdot\theta^{\varepsilon}+\frac{1}{2}\cdot\theta^{-\varepsilon}$ then contradicts the assumption of the extremity of $\theta$ in $\tilde{\Theta}^{N}_{D}$, because of $\sigma_{L}\neq 0$ for some $L\in\calS$.

To verify linear independence of the vectors $\{ \chi_{S}\,:\ S\in\calS\}$ in case $\bigcup\calS=N$ we first
realize that one has $\theta(N)>0$ then. Indeed, if $D=N$ then $\theta(S)<0$ for any $S\in\calS$ and
if $D\neq N$ then we choose $i\in N\setminus D$ and have $\theta(S)<0$ for any $S\in\calS$ with $i\in S$;
this allows us to use $\sum_{L\subseteq N:\, i\in L} \theta (L)=0$ for $i\in N$ to derive $\theta(N)>0$.

Assume for a contradiction that there is a non-zero linear combination $\sum_{S\in\calS} \sigma_{S}\cdot\chi_{S}=\0$
and put $\varsigma:=\sum_{S\in\calS} \sigma_{S}$. The case $\varsigma=0$ leads to a contradiction as shown in the
case of affine independence. Thus, consider $\varsigma\neq 0$, put $\sigma_{\emptyset}:=-\varsigma$,
$\sigma_{L}:=0$ for remaining $L\subseteq N$, and $\theta^{\varepsilon}(S):=(1-\varepsilon\cdot\varsigma)^{-1}\cdot (\theta(S)+\varepsilon\cdot\sigma_{S})$ for any $S\subseteq N$ and $\varepsilon\in {\dv R}$, $\varepsilon\neq \varsigma^{-1}$.
One has $\theta^{\varepsilon}\in \tilde{\Theta}^{N}_{D}$ for small $|\varepsilon|$. Moreover, $\theta=\frac{1-\varepsilon\cdot\varsigma}{2}\cdot\theta^{\varepsilon}+
\frac{1+\varepsilon\cdot\varsigma}{2}\cdot\theta^{-\varepsilon}$.
Because of $\varsigma\neq 0$ we have $\theta^{\varepsilon}(N)=(1-\varepsilon\cdot\varsigma)^{-1}\cdot\theta (N)\neq\theta(N)$ if $\varepsilon\neq 0$. Hence, we get a contradiction with the assumption of extremity
of $\theta$ in $\tilde{\Theta}^{N}_{D}$.

We have thus shown that the set system $\calS$ is min-semi-balanced on $N$ and one clearly has $\theta=\theta_{\calS}$.
If $\calS=\calB$ is min-balanced then the first option $\theta=\theta_{\calB}$ occurs.
If $\calS$ has an exceptional set $T\in\calS$ then $\theta(T)=\theta_{\calS}(T)>0$ which forces $T=D$. Thus,
$D$ is the exceptional set in $\calS$ in this case.
\smallskip

To verify the second claim assume that $\calS$ is a min-semi-balanced set system on $N$ and that
$D$ is the exceptional set within $\calS$. It is straightforward to evidence that $\theta_{\calS}\in \tilde{\Theta}^{N}_{D}$.
To show that $\theta_{\calS}\in\ext (\tilde{\Theta}^{N}_{D})$ assume for a contradiction that there is a non-trivial convex combination $\theta_{\calS}=\alpha\cdot\theta^{0}+ (1-\alpha)\cdot\theta^{1}$ with $\theta^{0},\theta^{1}\in \tilde{\Theta}^{N}_{D}$,
$\alpha\in (0,1)$, and $\theta^{0}\neq\theta^{1}$.
We know that $\theta_{\calS}(S)=0$ any $S\subseteq N$, $S\not\in\{\emptyset,N\}\cup\calS$ and observe that
$\theta^{0}(S)=0=\theta^{1}(S)$ for any such $S\subseteq N$ as well. Indeed, the
inclusion $\tilde{\Theta}^{N}_{D}\subseteq \Theta^{N}_{D}$ implies $\theta^{i}(L)\leq 0$ for $L\subseteq N$,
$L\not\in\{\emptyset,D,N\}$, which forces $\theta^{i}(S)=0$ for $S\subseteq N$, $S\not\in\{\emptyset,N\}\cup\calS$ .

For every $\varepsilon\geq 0$ we put $\theta^{\varepsilon} := \theta^{0} +\varepsilon\cdot (\theta^{1}-\theta^{0})$.
Note that $\theta^{\varepsilon}(S)=0$ for any $S\subseteq N$, $S\not\in\{\emptyset,N\}\cup\calS$.
The fact that $\tilde{\Theta}^{N}_{D}$ is bounded (see Lemma~\ref{lem.theta}) implies the existence of
$\varepsilon^{\diamond}:=\max\, \{\varepsilon\geq 1\,:\ \theta^{\varepsilon}\in\tilde{\Theta}^{N}_{D}\,\}$.
There exists $L\in\calS$ with $\theta^{\varepsilon^{\diamond}}(L)=0$ because otherwise one gets a contradiction
with maximality of $\varepsilon^{\diamond}$.
We put $\calC:=\{ S\in\calS\,:\ \theta^{\varepsilon^{\diamond}}(S)\neq 0\,\}$ and observe, using the fact
$\theta^{\varepsilon^{\diamond}}\in\tilde{\Theta}^{N}_{D}$, that $\calC$ is a semi-balanced system on $N$.
Then $L\in\calS\setminus\calC$ contradicts the minimality of $\calS$.
Thus, there is no non-trivial convex combination of $\theta^{0},\theta^{1}\in \tilde{\Theta}^{N}_{D}$
yielding $\theta_{\calS}$, which means that $\theta_{\calS}$ is a vertex of $\tilde{\Theta}^{N}_{D}$.
\end{proof}

\section{Proof of Lemma~\ref{lem.main-result}}\label{app.lem-main-result}
A key induction step in the proof is based on the following auxiliary observations.

\begin{lem}\label{lem.induction-step}~\rm
Given $|N|\geq 3$, let $\calS$ be a purely min-semi-balanced system on $N$.
Denote $\calW :=\caP\setminus (\{\emptyset,N\}\cup\calS)$ and introduce the next special
polytopes (see Definition~\ref{def.aux-cone}):
\begin{eqnarray*}
\Sigma(\calS) &:=& \conv (\bigcup_{D\in\calW} \tilde{\Theta}^{N}_{D})\,,\\
\Delta(\calS) &:=& \Sigma(\calS)\,\cap\, \{\theta\in {\dv R}^{\caP}\,:\ \theta(W)\geq 0 \quad\mbox{for any $W\in\calW$}\,\}\,.
\end{eqnarray*}
Moreover, for any $\calZ\subseteq\calW$, we put:
\begin{eqnarray*}
\Sigma^{\calZ}(\calS) &:=& \Sigma(\calS)\,\cap\, \{\theta\in {\dv R}^{\caP}\,:\ \theta(Z)= 0 \quad\mbox{for any $Z\in\calZ$}\,\}\,,\\
\Delta^{\calZ}(\calS) &:=& \Sigma^{\calZ}(\calS)\,\cap\, \Delta(\calS)\,.
\end{eqnarray*}
Then, for any $\calZ\subseteq\calW$, one has:
\begin{description}
\item[(i)] $\Sigma^{\calW}(\calS)=\emptyset$,
\item[(ii)] if $\emptyset\neq \Delta^{\calZ}(\calS)$ and $\Sigma^{\calZ}(\calS)\setminus \Delta^{\calZ}(\calS)\neq\emptyset$
~then\, $\emptyset\neq \Delta^{\calZ\cup\{D\}}(\calS)$ for some $D\in\calW\setminus\calZ$,
\item[(iii)] if $\theta\in \ext (\Sigma^{\calZ}(\calS))$ then $\theta\in \bigcup_{D\in\calW} \tilde{\Theta}^{N}_{D}$.
\end{description}
\end{lem}
\smallskip

\begin{proof}
To show \mbox{\bf (i)} assume for a contradiction that $\theta\in\Sigma^{\calW}(\calS)$ exists. Thus,
$\theta(W)=0$ for any $W\in\calW$ forces $\calB:=\{ S\subset N\,:\ S\neq\emptyset ~\&~ \theta(S)\neq 0\,\}\subseteq\calS$.
Since, however, every $\eta\in \bigcup_{D\in\calW} \tilde{\Theta}^{N}_{D}$ satisfies $\eta(S)\leq 0$ for any $S\in\calS$,
one has $\theta(S)\leq 0$ for any $S\in\calS$ and $\theta\in \Sigma(\calS)$. This implies, using the
equality constraints $\sum_{L\subseteq N:\,i\in N} \theta (L)=0$ for $i\in N$ and $\theta\in\Sigma(\calS)$
(use Definition~\ref{def.aux-cone}), that $\calB$ is balanced on $N$. The minimality of $\calS$ then implies that $\calB=\calS$,
which contradicts the assumption that $\calS$ is {\sf not\/} balanced on $N$.
\smallskip

To show \mbox{\bf (ii)} assume the existence of $\eta\in\Delta^{\calZ}(\calS)$. By (i) one has $\calW\setminus\calZ\neq\emptyset$.
If there is $D\in\calW\setminus\calZ$ with $\eta (D)=0$ then $\eta\in\Delta^{\calZ\cup\{D\}}(\calS)$ and we are done.
Thus, assume that $\eta(W)\neq 0$ for any $W\in\calW\setminus\calZ$, which implies, by definition of $\Delta(\calS)$,
that $\eta(W)> 0$ for any $W\in\calW\setminus\calZ$. The second assumption in (ii) means that there exists
$\theta\in\Sigma^{\calZ}(\calS)\setminus \Delta^{\calZ}(\calS)$; thus, necessarily $\eta\neq\theta$.
Put $\theta^{\alpha}:=(1-\alpha)\cdot\eta+\alpha\cdot\theta$ for $0\leq\alpha\leq 1$; the convexity of $\Sigma^{\calZ}(\calS)$
gives $\theta^{\alpha}\in\Sigma^{\calZ}(\calS)$. Take $\beta:=\max\, \{\,\alpha\geq 0\,:\ \theta^{\alpha}\in \Delta^{\calZ}(\calS)\,\}$
and $\theta^{\prime}:=\theta^{\beta}$. Note that there exists $D\in\calW\setminus\calZ$ with $\theta^{\prime}(D)=\theta^{\beta}(D)=0$
as otherwise one has $\theta^{\beta}(W)>0$ for any $W\in\calW\setminus\calZ$, which contradicts the maximality of $\beta$.
Thus, the facts $\theta^{\prime}\in \Delta^{\calZ}(\calS)$ and $\theta^{\prime}(D)=0$ imply together that
$\theta^{\prime}\in \Delta^{\calZ\cup\{D\}}(\calS)$.
\smallskip

The condition \mbox{\bf (iii)} can be verified by induction on $|\calZ|$. In case $\calZ=\emptyset$
this follows directly from the definition of $\Sigma(\calS)$ using basic facts from polyhedral geometry
recalled in Section~\ref{ssec.polyhedrals}. To verify the induction step assume the claim is true
for some $\calZ\subset\calW$, take $Z\in\calW\setminus\calZ$ and evidence the claim for $\calZ\cup\{Z\}$.
Note that the polytope ${\sf Q}:=\Sigma^{\calZ\cup\{Z\}}(\calS)$ is the intersection of
the polytope ${\sf P}:=\Sigma^{\calZ\cup\{Z\}}(\calS)$ with the
hyperplane ${\sf H}:=\{\,\theta\,:\ \theta(Z)=0\}$. 
The characterization of vertices of ${\sf Q}={\sf P}\cap {\sf H}$ recalled in Section~\ref{ssec.polyhedrals}
says that any vertex $\theta\in \ext ({\sf Q})$ is {\sf either\/} a vertex of ${\sf P}$, in which
case one has $\theta\in \bigcup_{D\in\calW} \tilde{\Theta}^{N}_{D}$ by the induction assumption, {\sf or\/}
there is an edge $[\eta,\sigma]$ of ${\sf P}$ such that $\theta\in\,]\eta,\sigma[$ and $[\eta,\sigma]\cap{\sf H}=\{\theta\}$.
Since $\eta,\sigma\in\ext ({\sf P})$, by the induction hypothesis, one has
$\eta,\sigma\in \bigcup_{D\in\calW} \tilde{\Theta}^{N}_{D}$.
As $[\eta,\sigma]\cap{\sf H}=\{\theta\}$ one has $\eta(Z)\neq 0\neq\sigma(Z)$ and
can assume without loss of generality that $\eta(Z)>0$ and $\sigma(Z)<0$. This forces
that $\eta\in\tilde{\Theta}^{N}_{Z}$; assume that $\sigma\in \tilde{\Theta}^{N}_{E}$ for some $E\in\calW$
(possibly $E=Z$). It makes no problem to observe that these facts imply that $\theta\in \tilde{\Theta}^{N}_{E}$
(see Definition~\ref{def.aux-cone}). Hence, $\theta\in \bigcup_{D\in\calW} \tilde{\Theta}^{N}_{D}$ and
the induction step has been verified.
\end{proof}

For reader's convenience we recall what is claimed.
\bigskip

\noindent{\bfseries Lemma~\ref{lem.main-result}:} ~~
Let $\calS$ be a purely min-semi-balanced system on $N$, $|N|\geq 3$, with exceptional set $T\in\calS$
and $\calW:= \caP\setminus (\{\emptyset,N\}\cup\calS)$.
Then the following conditions are equivalent:
\begin{description}
\item[(a)] $\theta_{\calS}\not\in \ext(\Delta)$, (see Definition~\ref{def.aux-cone})
\item[(b)] there exists a convex combination $\theta_{\calS}=\sum_{D\in\calW\cup\{T\}} \alpha_{D}\cdot\theta^{D}$
where $\alpha_{T}<1$ and $\theta^{D}\in\tilde{\Theta}^{N}_{D}$ whenever $\alpha_{D}>0$,
\item[(c)] the set $\Delta(\calS):=\conv(\bigcup_{D\in\calW} \tilde{\Theta}^{N}_{D}) \cap \{\,\theta\in {\dv R}^{\caP}\,:\
\theta(W)\geq 0 ~\mbox{for $W\in\calW$}\,\}$ is non-empty,
\item[(d)] there exists $E\in\calW$ such that $E$ is exceptional in $\calS\cup\{E\}$,
\item[(e)] there exists $E\in\calW$ such that $E$ is exceptional in $(\calS\setminus\{T\})\cup\{E\}$,
\item[(f)] a min-semi-balanced system $\calD$ on $N$ exists such that $\calD\setminus\calS=\{E\}$
for some $E\in\calW$, the set $E$ is exceptional within $\calD$, and $T\not\in\calD$,
\item[(g)] a min-semi-balanced system $\calD$ on $N$ exists with an exceptional set $E$ and $\calD\setminus\calS=\{E\}$.
\end{description}
\medskip

\begin{proof}
To show \mbox{\bf (a)}$\Rightarrow$\mbox{\bf (b)} assume that $\theta_{\calS}\in\,]\eta,\sigma[$
for $\eta,\sigma\in\Delta$, $\eta\neq\theta_{\calS}\neq\sigma$. By the definition of $\Delta$ a finite convex combination
yielding $\eta$ exists: $\eta=\sum_{j\in J} \beta_{j}\cdot \eta_{j}$, where $\beta_{j}>0$ and
$\eta_{j}\in\bigcup_{\emptyset\neq D\subseteq N} \tilde{\Theta}^{N}_{D}$ for $j\in J$.
Since $\eta\neq\theta_{\calS}$, there exists $j\in J$ with $\eta_{j}\neq\theta_{\calB}$.
An analogous convex combination exists for $\sigma$ and by combining them observe that
a finite convex combination yielding $\theta_{\calS}$ exists: $\theta_{\calS}=\sum_{i\in I} \alpha_{i}\cdot \theta_{i}$,
where $\alpha_{i}>0$ and $\theta_{i}\in\bigcup_{\emptyset\neq D\subseteq N} \tilde{\Theta}^{N}_{D}$ for $i\in I$,
where, moreover, $\theta_{i}\neq \theta_{\calS}$ for at least one $i\in I$. Assume without loss of generality
that $\theta_{i}$, $i\in I$, differ from each other and, if there is $i\in I$ with
$\theta_{i}=\theta_{\calS}$ then, by easy modification, one gets such a convex combination yielding $\theta_{\calS}$ where
$\theta_{i}\neq\theta_{\calS}$ for all $i\in I$.

Observe that one can even assume without loss of generality that $\theta_{i}(S)\leq 0$ for any $i\in I$ and $S\in\calS\setminus\{T\}$.
Indeed, if there exists $j\in I$ and $S\in\calS\setminus\{T\}$ with $\theta_{j}(S)>0$ then the convex combination
can be replaced by another convex combination $\theta_{\calS}=\hat{\alpha}_{j}\cdot \hat{\theta}_{j}+\sum_{i\in I\setminus \{j\}} \hat{\alpha}_{i}\cdot \theta_{i}$, where $\hat{\theta}_{j}(S)=0$. To observe that realize that $\alpha_{j}<1$, as otherwise
$\theta_{j}=\theta_{\calS}$. Thus, one can put $\alpha := 1-\alpha_{j}=\sum_{i\in I\setminus\{j\}} \alpha_{i}$ and
$$
\theta ~:=~ \sum_{i\in I\setminus\{j\}} \frac{\alpha_{i}}{\alpha}\cdot \theta_{i}\,,\quad
\mbox{which gives~ $\theta_{\calS}=\alpha_{j}\cdot\theta_{j}+ \alpha\cdot \theta$, that is, $\theta_{\calS}\in\,]\theta_{j},\theta[$.}
$$
On the other hand, since $\theta_{j}(S)>0$ and $\theta_{\calS}(S)<0$ there exists unique $\hat{\theta}_{j}\in\,]\theta_{j},\theta_{\calS}[$
such that $\hat{\theta}_{j}(S)=0$. Note that necessarily $\theta_{j}\in \tilde{\Theta}^{N}_{S}$ while
$\theta_{\calS}\in \tilde{\Theta}^{N}_{T}$ which facts together allow one to observe that $\hat{\theta}_{j}\in \tilde{\Theta}^{N}_{T}$.
The vectors $\theta_{j},\hat{\theta}_{j},\theta_{\calS}$ and $\theta$ are on the same line, which implies that
$\theta_{\calS}\in\,]\hat{\theta}_{j},\theta[$. Thus, $0<\gamma<1$ exists with $\theta_{\calS}=\gamma\cdot\hat{\theta}_{j}+
(1-\gamma)\cdot\theta$ and the substitution gets
$$
\theta_{\calS} ~=~ \gamma\cdot \hat{\theta}_{j} + \sum_{i\in I\setminus\{j\}} \underbrace{\frac{(1-\gamma)\cdot\alpha_{i}}{\alpha}}_{\hat{\alpha}_{i}}\,\cdot\, \theta_{i}\,,
$$
where it suffices to put $\hat{\alpha}_{j}:=\gamma$.

For any $i\in I$, one has $\theta_{i}\in \tilde{\Theta}^{N}_{D}$ for some $\emptyset\neq D\subseteq N$.
In case $D\in \calS\setminus\{T\}$ the above inequality $\theta_{i}(S)\leq 0$ for 
$S\in\calS\setminus\{T\}$ gives $\theta_{i}\in \tilde{\Theta}^{N}_{N}$
and the inclusion $\tilde{\Theta}^{N}_{N}\subseteq \tilde{\Theta}^{N}_{T}$ allows one to conclude that
$\theta_{i}\in \tilde{\Theta}^{N}_{D}$ for some $D\in\calW\cup\{T\}$.
To summarize that: there exists a convex combination  $\theta_{\calS}=\sum_{i\in I} \alpha_{i}\cdot \theta_{i}$,
where $\alpha_{i}>0$, $\theta_{i}\neq\theta_{\calS}$ and $\theta_{i}\in\bigcup_{D\in\calW\cup\{T\}} \tilde{\Theta}^{N}_{D}$ for
all $i\in I$.

Further observation is that there exists $j\in I$ and $W\in\calW$ such that $\theta_{j}(W)>0$. To show that assume
for a contradiction the converse, that is, $\theta_{i}(W)\leq 0$ for any $i\in I$ and $W\in\calW$. That basically means, that,
for any $i\in I$, the condition $\theta_{i}\in \tilde{\Theta}^{N}_{W}$ implies
$\theta_{i}\in \tilde{\Theta}^{N}_{N}\subseteq \tilde{\Theta}^{N}_{T}$. In particular, one would have
$\theta_{i}\in \tilde{\Theta}^{N}_{T}$, $\theta_{i}\neq\theta_{\calS}$ for any $i\in N$, which contradicts
the fact $\theta_{\calS}\in \ext (\tilde{\Theta}^{N}_{T})$ claimed by (the second claim in) Lemma~\ref{lem.extr-theta}.

Thus, any $i\in I$ can be assigned to some $D\in\calW\cup\{T\}$ such that $\theta_{i}\in\tilde{\Theta}^{N}_{D}$;
let us fix that choice and write $i\mapsto D$ to denote that. We already know that there exists
$j\in I$ and $W\in\calW$ with $\theta_{j}(W)>0$, which necessitates $j\mapsto W$.
For any $D\in \calW\cup\{T\}$ we put
$$
\alpha_{D} ~:=~ \sum_{i\in I:\, i\,\mapsto D} \alpha_{i} \qquad \mbox{and, if $\alpha_{D}>0$,}\quad
\theta^{D} ~:=~ \sum_{i\in I:\, i\,\mapsto D} \frac{\alpha_{i}}{\alpha_{D}}\cdot \theta_{i}\in \tilde{\Theta}^{N}_{D}\,.
$$
Thus, one has $\theta_{\calS}=\sum_{D\in\calW\cup\{T\}} \alpha_{D}\cdot\theta^{D}$, where
$\theta^{D}$ can be chosen arbitrarily in case $\alpha_{D}=0$. It is a convex combination and the
existence of $j\in I$ and $W\in\calW$ with $j\mapsto W$ implies $\alpha_{W}>0$. Hence, $\alpha_{T}<1$,
which gives the condition \mbox{\bf (b)}.
\smallskip

To show \mbox{\bf (b)}$\Rightarrow$\mbox{\bf (c)} consider a convex combination
$\theta_{\calS}=\sum_{D\in\calW\cup\{T\}} \alpha_{D}\cdot\theta^{D}$ where $\alpha_{T}<1$ and
$\theta^{D}\in\tilde{\Theta}^{N}_{D}$ whenever $\alpha_{D}>0$.
%
The fact $\theta_{\calS}(T)>0$ forces both $\alpha_{T}>0$ and $\theta^{T}(T)>0$ because $\eta(T)\leq 0$ for any
$\eta\in \bigcup_{D\in\calW} \tilde{\Theta}^{N}_{D}$. In particular, one has $\theta^{T}\in \tilde{\Theta}^{N}_{T}$.
Let us put
$$
\theta ~:=~ \sum_{D\in\calW}\,\, \frac{\alpha_{D}}{1-\alpha_{T}}\,\cdot\, \theta^{D}\qquad \mbox{and observe that~
$\theta_{\calS}=(1-\alpha_{T})\cdot\theta+ \alpha_{T}\cdot\theta^{T}$.}
$$
By definition, $\theta\in \conv(\bigcup_{D\in\calW} \tilde{\Theta}^{N}_{D})$. The fact $\theta_{\calS}\in\,]\theta,\theta^{T}[$
together with $\theta_{\calS}(W)=0$ and  $\theta^{T}(W)\leq 0$ for any $W\in\calW$ forces
$\theta(W)\geq 0$ for any such $W$. In particular, $\theta\in\Delta(\calS)$ and the condition {\bf (c)} has been verified.
\smallskip

To show \mbox{\bf (c)}$\Rightarrow$\mbox{\bf (d)} we use Lemma~\ref{lem.induction-step};
the condition \mbox{\bf (c)} means, by notation from Lemma~\ref{lem.induction-step}, that
$\Delta^{\calZ}(\calS)\neq\emptyset$ for empty $\calZ=\emptyset$.
By inductive application of Lemma~\ref{lem.induction-step}(ii) we find $\calZ\subseteq\calW$ with
$\Sigma^{\calZ}(\calS)=\Delta^{\calZ}(\calS)\neq\emptyset$; indeed, by Lemma~\ref{lem.induction-step}(i)
one has $\Sigma^{\calW}(\calS)=\emptyset$, which ensures that inductive enlarging of $\calZ$
has to finish with some desired $\calZ\subset\calW$. Thus, $\emptyset\neq \ext (\Sigma^{\calZ}(\calS))\subseteq \Delta^{\calZ}(\calS)$
for some $\calZ\subset\calW$. One can take some $\theta\in \ext (\Sigma^{\calZ}(\calS))$ and, by
Lemma~\ref{lem.induction-step}(iii), there exists $E\in\calW$ such that
$\theta\in \tilde{\Theta}^{N}_{E}$. Hence, $\theta(W)\leq 0$ for any $W\in\calW\setminus\{E\}$ while
$\theta\in\Delta(\calS)$ says $\theta(W)\geq 0$ for any $W\in\calW$.
That together means $\theta(W)=0$ for any $W\in \calW\setminus\{E\}$ and $\theta(E)\geq 0$.
One cannot have $\theta(E)=0$ for otherwise $\theta\in\Sigma^{\calW}(\calS)$ contradicts the
claim in Lemma~\ref{lem.induction-step}(i). Thus, necessarily $\theta(E)>0$; note also that
$\theta(S)\leq 0$ for any $S\in\calS$.
The equality constraints $\sum_{L\subseteq N:\,i\in N} \theta (L)=0$ for any $i\in N$ and $\theta\in \tilde{\Theta}^{N}_{E}$
allow one to conclude that $\sum_{S\in\calS\cup\{E\}} -\theta(S)\cdot\chi_{S}=\theta(N)\cdot\chi_{N}$ is
a semi-conic combination yielding a constant vector in ${\dv R}^{N}$. Thus, by definition,
$E$ is exceptional in $\calS\cup\{E\}$ and the condition \mbox{\bf (d)} has been verified.
\smallskip

To show \mbox{\bf (d)}$\Rightarrow$\mbox{\bf (e)} assume that $E\in\calW$ is exceptional within $\calS\cup\{E\}$.
This means that there exists a linear combination $\sum_{S\in\calS\cup\{E\}} \nu_{S}\cdot\chi_{S}$
yielding a constant vector in ${\dv R}^{N}$ where $\nu_{E}<0$ and $\nu_{S}\geq 0$ for $S\in\calS$.
Observe that $E$ is exceptional within a smaller set system $(\calS\setminus\{T\})\cup\{E\}$.
In case $\nu_{T}=0$ we are done.
In case $\nu_{T}>0$ we apply Lemma~\ref{lem.min-semi-bal} to $\calS$ which says that there exists (unique)
affine combination $\sum_{S\in\calS} \lambda_{S}\cdot\chi_{S}$ yielding a constant vector in ${\dv R}^{N}$
with both $\lambda_{T}<0$ and $\lambda_{S}> 0$ for $S\in\calS\setminus\{T\}$.
Then we put $\lambda_{E}:=0$ and
$$
\kappa_{S} ~:=~ \frac{\nu_{T}}{\nu_{T}-\lambda_{T}}\cdot \lambda_{S} + \frac{-\lambda_{T}}{\nu_{T}-\lambda_{T}}\cdot \nu_{S}
\qquad \mbox{for $S\in\calS\cup\{E\}$}.
$$
Hence, $\sum_{S\in\calS\cup\{E\}} \kappa_{S}\cdot\chi_{S}$ yields a constant vector in ${\dv R}^{N}$,
$\kappa_{E}<0$, and $\kappa_{T}=0$. Thus, $(\calS\setminus\{T\})\cup\{E\}$ is a semi-balanced system on $N$
and $E$ is an exceptional set within it.
\smallskip

To show \mbox{\bf (e)}$\Rightarrow$\mbox{\bf (f)} we fix a semi-conic combination $\sum_{S\in \calS\cup\{E\}} \kappa_{S}\cdot\chi_{S}$ yielding a constant vector in ${\dv R}^{N}$ with $\kappa_{E}<0$ and $\kappa_{T}=0$.
In particular, one can choose a min-semi-balanced system $\calD\subseteq(\calS\setminus\{T\})\cup\{E\}$.
One has $E\in\calD$ as otherwise $\calD\subset\calS$ contradicts the minimality of $\calS$.
By Lemma~\ref{lem.min-semi-bal}, there exists (unique) affine combination
$\sum_{S\in\calD} \sigma_{S}\cdot\chi_{S}$ yielding a constant vector in ${\dv R}^{N}$
which is semi-conic and has all coefficients non-zero. Assume for a contradiction that $\sigma_{E}>0$. Then we put $\sigma_{S}:=0$
for $S\in \calS\setminus\calD$ and
$$
\mu_{S} ~:=~ \frac{\sigma_{E}}{\sigma_{E}-\kappa_{E}}\cdot \kappa_{S} + \frac{-\kappa_{E}}{\sigma_{E}-\kappa_{E}}\cdot \sigma_{S}
\qquad \mbox{for $S\in\calS\cup\{E\}$}.
$$
Hence, $\sum_{S\in\calS\cup\{E\}} \mu_{S}\cdot\chi_{S}$ yields a constant vector in ${\dv R}^{N}$ and $\mu_{E}=0=\mu_{T}$. Since it is a semi-conic combination $\calT:=\{ S\in\calS\,:\ \mu_{S}\neq 0\}\subset\calS$ is a semi-balanced system on $N$, which
contradicts the minimality of $\calS$. As $\sigma_{E}\neq 0$, one necessarily has $\sigma_{E}<0$ and
the set $E$ is exceptional within $\calD$. Thus, the condition \mbox{\bf (f)} has been verified.
\smallskip

The implication \mbox{\bf (f)}$\Rightarrow$\mbox{\bf (g)} is evident.
\smallskip

To show \mbox{\bf (g)}$\Rightarrow$\mbox{\bf (a)} assume that $\calD$ is the min-semi-balanced system
with an exceptional set $E\in\calW$ such that $\calD\setminus\calS=\{E\}$. Note that, by \eqref{eq.theta},
one has both $\theta_{\calD}(E)>0$ and $\theta_{\calD}(W)=0$ for any $W\in\calW\setminus\{E\}$.
Then we put
$$
\theta^{\varepsilon} ~:=~ (1+\varepsilon)\cdot\theta_{\calS} + (-\varepsilon)\cdot\theta_{\calD}
~=~ \theta_{\calS} +\varepsilon\cdot (\theta_{\calS}-\theta_{\calD})\qquad
\mbox{for every $\varepsilon\geq 0$.}
$$
It makes no problem to observe that, for small $\varepsilon>0$, one has $\theta^{\varepsilon}\in\tilde{\Theta}^{N}_{T}\subseteq\Delta$.
Because of $\theta_{\calS}\in\, ]\theta_{\calD},\theta^{\varepsilon}[$ and $\theta_{\calD}\in\tilde{\Theta}^{N}_{E}\subseteq\Delta$
one gets $\theta_{\calS}\not\in\ext(\Delta)$ and \mbox{\bf (a)} has been verified.
\end{proof}

\end{document}